\documentclass[reqno]{conm-p-l} 
\usepackage{amssymb}

\makeatletter 

 \def\activeat#1{\csname @#1\endcsname}
 \def\def@#1{\expandafter\def\csname @#1\endcsname}
 {\catcode`\@=\active \gdef@{\activeat}}

\let\ssize\scriptstyle
\newdimen\ex@	\ex@.2326ex

 \def\requalfill{\cleaders\hbox{$\mkern-2mu\mathord=\mkern-2mu$}\hfill
  \mkern-6mu\mathord=$}
 \def\eqfill{$\m@th\mathord=\mkern-6mu\requalfill}
 \def\deffill{\hbox{$:=$}$\m@th\mkern-6mu\requalfill}
 \def\fiberbox{\hbox{$\vcenter{\hrule\hbox{\vrule\kern1ex
     \vbox{\kern1.2ex}\vrule}\hrule}$}}

 \newdimen\arrwd 
  \newdimen\minCDarrwd \minCDarrwd=2.5pc
   \setbox\z@\hbox{$\rightarrow\,$} \minCDarrwd=\wd\z@
 	
 \def\findarrwd#1#2#3{\arrwd=#3%
  \setbox\z@\hbox{$\ssize\;{#1}\;\;$}%
  \setbox\@ne\hbox{$\ssize\;{#2}\;\;$}%
  \ifdim\wd\z@>\arrwd \arrwd=\wd\z@\fi
  \ifdim\wd\@ne>\arrwd \arrwd=\wd\@ne\fi}
 \newdimen\arrowsp\arrowsp=0.375em  	
 \def\findCDarrwd#1#2{\findarrwd{#1}{#2}{\minCDarrwd}
    \advance\arrwd by 2\arrowsp}
 \newdimen\minarrwd 
 \setbox\z@\hbox{$\longrightarrow$} \minarrwd=\wd\z@

 \def\harrow#1#2#3#4{{\minarrwd=#1\minarrwd
   \findarrwd{#2}{#3}{\minarrwd}\kern\arrowsp
    \mathrel{\mathop{\hbox to\arrwd{#4}}\limits^{#2}_{#3}}\kern\arrowsp}}

 \def@]#1>#2>#3>{\harrow{#1}{#2}{#3}\rightarrowfill}
 \def@>#1>#2>{\harrow1{#1}{#2}\rightarrowfill}
 \def@<#1<#2<{\harrow1{#1}{#2}\leftarrowfill}
 \def@={\harrow1{}{}\eqfill}
 \def@:#1={\harrow1{}{}\deffill}
 \def@ N#1N#2N{\vCDarrow{#1}{#2}\UpDownarrow}
 \def\UpDownarrow{\uparrow\,\Big\downarrow}

\def\hookrightarrowfill{\hbox{$\lhook\joinrel$}\rightarrowfill}
\def@{hk>}#1>#2>{\harrow1{#1}{#2}\hookrightarrowfill}
\def@ c>#1>#2>{\harrow1{#1}{#2}\hookrightarrowfill}
\def\hookleftarrowfill{\leftarrowfill\hbox{$\joinrel\rhook$}}
\def@{hk<}#1<#2<{\harrow1{#1}{#2}\hookleftarrowfill}

 \def@.{\ifodd\row\relax\harrow1{}{}\hfill
   \else\vCDarrow{}{}.\fi}
 \def@|{\vCDarrow{}{}\Vert}
 \def@ V#1V#2V{\vCDarrow{#1}{#2}\downarrow}
 \def@ A#1A#2A{\vCDarrow{#1}{#2}\uparrow}
 \def@(#1){\arrwd=\csname col\the\col\endcsname\relax
   \hbox to 0pt{\hbox to \arrwd{\hss$\vcenter{\hbox{$#1$}}$\hss}\hss}}

 \def\squash#1{\setbox\z@=\hbox{$#1$}\finsm@@sh}
\def\finsm@@sh{\ifnum\row>1\ht\z@\z@\fi \dp\z@\z@ \box\z@}

 \newcount\row \newcount\col \newcount\numcol \newcount\arrspan
 \newdimen\vrtxhalfwd  \newbox\tempbox

 \def\innernewdimen{\alloc@1\dimen\dimendef\insc@unt}
 \def\measureinit{\col=1\vrtxhalfwd=0pt\arrspan=1\arrwd=0pt 
   \setbox\tempbox=\hbox\bgroup$}
 \def\setinit{\col=1\hbox\bgroup$\ifodd\row
   \kern\csname col1\endcsname
   \kern-\csname row\the\row col1\endcsname\fi}
 \def\findvrtxhalfsum{$\egroup
  \expandafter\innernewdimen\csname row\the\row col\the\col\endcsname
  \global\csname row\the\row col\the\col\endcsname=\vrtxhalfwd
  \vrtxhalfwd=0.5\wd\tempbox
  \global\advance\csname row\the\row col\the\col\endcsname by \vrtxhalfwd 
  \advance\arrwd by \csname row\the\row col\the\col\endcsname
  \divide\arrwd by \arrspan
  \loop\ifnum\col>\numcol \numcol=\col%
     \expandafter\innernewdimen \csname col\the\col\endcsname
     \global\csname col\the\col\endcsname=\arrwd
   \else \ifdim\arrwd >\csname col\the\col\endcsname
      \global\csname col\the\col\endcsname=\arrwd\fi\fi
   \advance\arrspan by -1 %
   \ifnum\arrspan>0 \repeat}
 \def\setCDarrow#1#2#3#4{\advance\col by 1 \arrspan=#1 
    \arrwd= -\csname row\the\row col\the\col\endcsname\relax
    \loop\advance\arrwd by \csname col\the\col\endcsname
     \ifnum\arrspan>1 \advance\col by 1 \advance\arrspan by -1%
     \repeat
    \squash{\mathop{
     \hbox to\arrwd{\kern\arrowsp#4\kern\arrowsp}}\limits^{#2}_{#3}}}
 \def\measureCDarrow#1#2#3#4{\findvrtxhalfsum\advance\col by 1%
   \arrspan=#1\findCDarrwd{#2}{#3}%
    \setbox\tempbox=\hbox\bgroup$}
 \def\vCDarrow#1#2#3{\kern\csname col\the\col\endcsname
    \hbox to 0pt{\hss$\vcenter{\llap{$\ssize#1$}}%
     \Big#3\vcenter{\rlap{$\ssize#2$}}$\hss}\advance\col by 1}

 \def\setCD{\def\harrow{\setCDarrow}%
  \def\\{$\egroup\advance\row by 1\setinit}
  \m@th\lineskip3\ex@\lineskiplimit3\ex@ \row=1\setinit}
 \def\endsetCD{$\egroup}
 \def\drop#1\\{\findvrtxhalfsum\advance\row by 2 \measureinit}
 \def\measure{\bgroup
  \def\harrow{\measureCDarrow}%
  \def\\##1{\ifx##1\endmeasure\endmeasure\else\expandafter\drop\fi}%
  \row=1\numcol=0\measureinit}
 \def\endmeasure{\findvrtxhalfsum\egroup}

 \newcount\savedcount	
 \def\LCD#1\end{\savedcount=\count11
   \measure#1\endmeasure
   \vcenter{\setCD#1\endsetCD\kern\medskipamount}%
   \global\count11=\savedcount\end}

 \newenvironment{CD}{\let\at=@\catcode`\@=\active\LCD}{\catcode`\@=12\relax}

\makeatother 

\newcommand{\emdash}{\unskip\penalty10000\thinspace
	---\penalty-500\thinspace\ignorespaces}

\DeclareMathOperator{\Br}{Br}
\DeclareMathOperator{\Div}{Div}
\DeclareMathOperator{\IDiv}{\bf Div}
\DeclareMathOperator{\et}{\text{\rm(\'et)}}
\DeclareMathOperator{\fppf}{(fppf)}
\DeclareMathOperator{\fpqc}{(fpqc)}
\DeclareMathOperator{\Hilb}{Hilb}
\DeclareMathOperator{\IHilb}{\bf Hilb}
\DeclareMathOperator{\length}{\rm length}
\DeclareMathOperator{\IPic}{\bf Pic}
 \def\IPicz{\IPic^{\smash0}}
 \def\IPicm{\IPic^{\smash m}}
 \def\IPict{\IPic^{\smash\tau}}
 \def\IPicp{\IPic^{\smash\phi}}
\DeclareMathOperator{\IQuot}{\bf Quot}
\DeclareMathOperator{\Hom}{Hom}
\DeclareMathOperator{\SHom}{\it Hom}
\DeclareMathOperator{\LinSys}{LinSys}
\DeclareMathOperator{\Pic}{Pic}
\DeclareMathOperator{\Proj}{Proj}
\DeclareMathOperator{\cok}{cok}
\DeclareMathOperator{\rank}{rank}
\DeclareMathOperator{\red}{red}
\DeclareMathOperator{\reg}{reg}
\DeclareMathOperator{\Spec}{Spec}
\DeclareMathOperator{\Supp}{Supp}
\DeclareMathOperator{\Tor}{Tor}
\DeclareMathOperator{\zar}{(zar)}

\let\?=\overline
\let\:=\colon
\let\bb=\mathbb
\let\I=\mathbf
\let\into=\hookrightarrow
\let\mc=\mathcal
\let\onto=\twoheadrightarrow
\let\ox=\otimes
\let\To=\Longrightarrow
\let\tu=\textup
\let\x=\times
\let\ve=\varepsilon
\let\vf=\varphi
\let\vt=\vartheta
\let\wh=\widehat
\let\wt=\widetilde

\newcommand{\Gm}{\bb G_\tu m}
\def\risom{\buildrel\sim\over{\smashedlongrightarrow}}
 \def\smashedlongrightarrow{\setbox0=\hbox{$\longrightarrow$}\ht0=1.25pt\box0}
\newcommand{\uH}{\tu H}
\def\bigomega{\mbox{\LARGE $\omega$}}

\theoremstyle{plain}
 \newtheorem{thm}{Theorem}[section] 
 \newtheorem{cor}[thm]{Corollary}
 \newtheorem{lem}[thm]{Lemma} 
 \newtheorem{prp}[thm]{Proposition}

\newtheoremstyle{sbs}
  {\smallskipamount}
  {\smallskipamount}
  {\it}
  {\parindent}
  {\rm}
  {.}
  {.5em}
  {}
\theoremstyle{sbs}
 \newtheorem{sbsthm}{Theorem}[thm]
 \newtheorem{sbscor}[sbsthm]{Corollary}

\theoremstyle{definition}
 \newtheorem{dfn}[thm]{Definition}
 \newtheorem{eg}[thm]{Example}
 \newtheorem{ex}[thm]{Exercise}
 \newtheorem{rmk}[thm]{Remark}
 \newtheorem{sbs}[thm]{} 

\makeatletter
 \@addtoreset{equation}{thm}
\makeatother
\renewcommand{\theequation}{\arabic{section}.\arabic{equation}}

\newenvironment{ans}[1]
  {\renewcommand{\theequation}{A.\ref{ex:#1}.\arabic{equation}}
\begin{proof}[Answer \ref{ex:#1}]}
  {\end{proof}}

\def\mylistparam
 {\setlength{\topsep}{\smallskipamount}%
  \setlength{\itemsep}{0pt}\setlength{\parsep}{0pt}%
  \setlength{\itemindent}{0pt}%
  \setlength{\leftmargin}{\parindent}%
  \setbox0=\hbox{1.\kern0.5em}\addtolength{\leftmargin}{\wd0}%
 }%
\renewenvironment{enumerate}%
 {\begin{list}{\llap{\rm (\roman{enumi})}}%
  {\usecounter{enumi}\mylistparam}%
  }%
 {\end{list}}

\begin{document}

\title
  {The Picard scheme} 

\author[S. L. Kleiman]{Steven L. Kleiman}
 \address
 {Math Department, MIT \\
 77 Massachusetts Avenue\\
 Cambridge, MA 02139, USA}
 \email{kleiman@math.mit.edu}

 \thanks{The author is grateful to Allen Altman, Ethan Cotterill,
Eduardo Esteves, Rebecca Lehman, Jean-Pierre Serre, and Angelo Vistoli
for their many valuable comments on earlier versions of these notes.}

\begin{abstract}
 We develop in detail most of the theory of the Picard scheme that
Grothendieck sketched in two Bourbaki talks and in commentaries on them.
Also, we review in brief much of the rest of the theory developed by
Grothendieck and by others.  But we begin with a historical
introduction.
 \end{abstract}
 
\maketitle 
\tableofcontents

\section{Introduction}\label{sc:intro}

On any ringed space $X$, the isomorphism classes of invertible sheaves
form a group; it is denoted by $\Pic(X)$, and called the (absolute) {\it
Picard group}.  Suppose $X$ is a ``projective variety''; in other words,
$X$ is an integral scheme that is projective over an algebraically
closed field $k$.  Then, as is proved in these notes, the group
$\Pic(X)$ underlies a natural $k$-scheme, which is a disjoint union of
quasi-projective schemes, and the operations of multiplying and of
inverting are given by $k$-maps.  This scheme is denoted by
$\IPic_{X/k}$, and called the {\it Picard scheme}.  It is reduced in
characteristic zero, but not always in positive characteristic.  When
$X$ varies in an algebraic family, correspondingly, $\IPic_{X/k}$ does
too.

The Picard scheme was introduced in 1962 by Grothendieck. He sketched
his theory in two Bourbaki talks, nos.~232 and 236, which were reprinted
along with his commentaries in \cite{FGA}.  But Grothendieck advanced an
old subject, which was actively being developed by many others at the
time.  Nevertheless, Grothendieck's theory was revolutionary, both in
concept and in technique.

In order to appreciate Grothendieck's contribution fully, we have to
review the history of the Picard scheme.  Reviewing this history serves
as well to introduce and to motivate Grothendieck's theory.
Furthermore, the history is rich and fascinating, and it is a
significant part of the history of algebraic geometry.

So let us now review the history of the Picard scheme up to 1962.  We
need only summarize and elaborate on scattered parts of Brigaglia,
Ciliberto, and Pedrini's article \cite{BCP} and of the author's article
\cite{Kl04}.  Both articles give many precise references to the original
sources and to the secondary literature; so few references are given
here.

The Picard scheme has roots in the 1600s.  Over the course of that
century, the Calculus was developed, through the efforts of many
individuals, in order to design lenses, to aim cannons, to make clocks,
to hang cables, and so on.  Thus interest arose in the properties of
functions appearing as indefinite integrals.

Notably, in 1694, James Bernoulli analyzed the way rods bend, and was
led to introduce the ``lemniscate,'' a figure eight with equation
$(x^2+y^2)^2=a^2(x^2-y^2)$
 where $a$ is nonzero.  In polar coordinates, he found the arc length
$s$ to be given by
	$$s=\int_0^r\frac{a^2\,dr}{\sqrt{a^4-r^4}}.$$
He surmised that $s$ can not be expressed in terms of the elementary
functions.  Similar integrals had already arisen in attempts to rectify
elliptical orbits; so these integrals became known as ``elliptic
integrals.''

In 1698, James's brother, John, recalled there are algebraic relations
among the arguments of sums and differences of logarithms and of the
inverse trigonometric functions.  Then he showed that, similarly, given
two arcs from the origin on the cubical parabola $y=x^3$, their lengths
differ by the length of a certain third such arc.  And he posed the
problem of finding more cases of this phenomenon.

Sure enough, between 1714 and 1720, Fagnano found, in an ad hoc manner,
similar relations for the cords and arcs of ellipses, hyperbolas, and
lemniscates.  In turn, Fagnano's work led Euler in 1757 to discover the
``addition formula''
	$$\int_0^{x_1} \frac{dx}{\sqrt{1-x^4}}
	\pm\int_0^{x_2} \frac{dx}{\sqrt{1-x^4}}
	=\int_0^{x_3} \frac{dx}{\sqrt{1-x^4}}$$
 where the variables $x_1,\,x_2,\,x_3$ must satisfy the symmetric relation
\begin{multline*}
	x_1^4x_2^4x_3^4
	+2x_1^4x_2^2x_3^2+2x_1^2x_2^4x_3^2+2x_1^2x_2^2x_3^4\\
	+x_1^4+x_2^4+x_3^4
	-2x_1^2x_2^2-2x_1^2x_3^2-2x_2^2x_3^2=0.  \end{multline*}
 In 1759, Euler generalized this formula to some other elliptic
integrals.  Specifically, Euler found the sum or difference of two to be
equal to a certain third plus an elementary function.  Moreover, he
expressed regret that he could handle only square roots and fourth
powers, but not higher roots or powers.

In 1826, Abel made a great advance: he discovered an addition theorem of
sweeping generality.  It concerns certain algebraic integrals, which
soon came to be called ``Abelian integrals.''  They are of the following
form:
	$$\psi x := \int^x_{x_0} R(x,y)\,dx$$
 where $x$ is an independent complex variable, $R$ is a rational
function, and $y=y(x)$ is an integral algebraic function; that is,
$y$ is the implicit multivalued function defined by an irreducible
equation of the form
	$$ f(x,y):=y^n+f_1(x)y^{n-1}+\dots+f_n(x)=0$$
where the $f_i(x)$ are polynomials in $x$.

Let $p$ be the genus of the curve $f=0$, and let $h_1,\dotsc,h_\alpha$
be rational numbers.  Then Abel's addition theorem asserts that
	$$h_1\psi x_1+\cdots+h_\alpha\psi x_\alpha
	 = v+\psi x_1'+\cdots+\psi x_p'  $$
 where $v$ is an elementary function of the independent variables
$x_1,\dots,x_\alpha$ and where $x_1',\dots,x_p'$ are algebraic
functions of them.  More precisely, $v$ is a complex-linear combination
of one algebraic function of $x_1,\dotsc,x_\alpha$ and of logarithms of
others; moreover, $x_1',\dots,x_p'$ work for every choice
of\/ $\psi x$.  Lastly, $p$ is minimal: given algebraic functions
$x_1',\dots,x_{p-1}'$ of\/ $x_1,\dotsc,\penalty0 x_p$, there exists an
integral $\psi x$ such that, for any elementary function $v$,
	$$\psi x_1+\cdots+\psi x_p
	 \neq v+\psi x_1'+\cdots+\psi x_{p-1}'.$$

Abel finished his 61-page manuscript in Paris, and submitted it in
person on 30 October 1826 to the Royal Academy of Sciences, which
appointed Cauchy and Legendre as referees.  However, the Academy did not
publish it until 1841, long after Abel's death from tuberculosis on 6
April 1829.

Meanwhile, Abel feared his manuscript was lost forever.  So in Crelle's
Journal, {\bf 3} (1828), he summarized his general addition theorem
informally.  Then he treated in detail a major special case, that in
which $f(x,y):=y^2-\varphi(x)$ where $\varphi(x)$ is a square-free
polynomial of degree $d\ge1$.  In particular, Abel found
	$$p=\begin{cases}(d-1)/2, &\text{if $d$ is odd};\\
			  (d-2)/2, &\text{if $d$ is even}.
		\end{cases}$$
Thus, if $d\ge5$, then $p\ge2$, and so Euler's formula does not extend.

With Jacobi's help, Legendre came to appreciate the importance of this
case.  To it, Legendre devoted the third supplement to his long treatise
on elliptic integrals, which are recovered when $d=3,4$.
For $d\ge5$, the integrals share many of the same formal properties.  So
Legendre termed them ``ultra-elliptiques.''

Legendre sent a copy of the supplement to Crelle for review on 24 March
1832, and Crelle asked Jacobi to review it.  Jacobi translated
``ultra-elliptiques'' by ``hyperelliptischen,'' and the prefix ``hyper''
has stuck.  In his cover letter, Legendre praised Abel's addition
theorem, calling it, in the immortal words of Horace's Ode 3, XXX.1, ``a
monument more lasting than bronze'' (monumentum aere perennius).  In his
review, Jacobi said that the theorem would be a most noble monument were
it to acquire the name {\bf Abel's Theorem.}  And it did!

Jacobi was inspired to give, a few months later, the first of several
proofs of Abel's Theorem in the hyperelliptic case.  Furthermore, he
posed the famous problem, which became known as the ``Jacobi Inversion
Problem.''  He asked, ``what, in the general case, are those functions
whose inverses are Abelian integrals, and what does Abel's theorem show
about them?''

Jacobi solved the inversion problem when $d=5,6$.  Namely, he formed
 $$\psi x := \int^x_{x_0}\frac{dx}{\sqrt{\varphi(x)}}\text{ and }
      \psi_1 x := \int^x_{x_0} \frac{x\,dx}{\sqrt{\varphi(x)}},$$
 and he set
 $$\psi x + \psi y = u\text{ and } \psi_1 x + \psi_1 y = v.$$
 He showed $x+y$ and $xy$ are single-valued functions of $u$ and $v$
with four periods.

Some historians have felt Abel had this inversion in mind, but ran out
of time.  At any rate, in 1827, Abel had originated the idea of
inverting elliptic integrals, obtaining what became known as ``elliptic
functions.''  He and Jacobi studied them extensively.  Moreover, Jacobi
introduced ``theta functions'' as an aid in the study; they were
generalized by Riemann in 1857, and used to solve the inversion problem
in arbitrary genus.

Abel's paper on hyperelliptic integrals fills twelve pages.  Eight are
devoted to a computational proof of a key intermediate result.  A half
year later, in Crelle's Journal, {\bf 4} (1829), Abel published a 2-page
paper with a conceptual proof of this result for any Abelian integral
$\psi x$.  The result says that  $\psi x_1+\cdots+\psi x_\mu$ is
equal to an elementary function $v$ if $x_1,\dotsc,x_\mu$ are not
independent, but are the abscissas of the variable points of
intersection of the curve $f=0$ with a second plane curve that varies in
a linear system \emdash although this geometric formulation is Clebsch's.

In each of the first two papers, Abel addressed two more, intermediate
questions:  First, when is the sum $\psi x_1+\cdots+\psi x_\mu$
constant?  Second, what is the number $\alpha$ of $x_i$ that can vary
independently?  Remarkably, the answers involve the genus $p$.

For hyperelliptic integrals, Abel found $\psi x_1+\cdots+\psi x_\mu$ is
constant if $\psi x$ is a linear combination of the following $p$
integrals:
	$$\int^x_{x_0} \frac{dx}{\sqrt{\varphi(x)}},\
	  \int^x_{x_0} \frac{x\,dx}{\sqrt{\varphi(x)}},\dotsc,
	  \int^x_{x_0} \frac{x^{p-1}\,dx}{\sqrt{\varphi(x)}}.$$
 Here $x_1,\dotsc,x_\mu$ are the abscissas of the variable points of
intersection of the curve $y^2=\varphi(x)$ and the curve
$\theta_1(x)y=\theta_2(x)$, whose coefficients vary, but $\theta_2(x)$
and $\varphi(x)$ retain a fixed common factor $\varphi_1(x)$.
Furthermore, Abel found
\begin{equation}\label{equation:1ineq}
 	\mu-\alpha\ge p;
 \end{equation}
 equality does not always hold, but can be achieved, given $d$ and
$\alpha$, by choosing the degrees of $\theta_1(x)$ and $\varphi_1(x)$
appropriately.

  Suppose $\mu-\alpha=p$.  Then
 $$\psi x_1+\cdots+\psi x_\alpha=v-(\psi x_{\alpha+1}+\cdots+\psi
	x_{\alpha+p})$$
 where $x_{\alpha+1},\dotsc,x_{\alpha+p}$ are algebraic functions of
$x_1,\dotsc,x_{\alpha}$.  Similarly, given any
$x_1',\dotsc,x_{\alpha'}'$, we get
     $$\psi x_1'+\cdots+\psi x_{\alpha'}'+\psi x_{\alpha+1}+\cdots
	   +\psi x_{\alpha+p} =v'-(\psi x_1''+\cdots+\psi x_p'').$$
Subtract this formula from the one above, and set $V:=v-v'$.  The result
is
   $$\psi x_1+\cdots+\psi x_\alpha-\psi x_1'-\cdots-\psi x_{\alpha'}'
	=V+\psi x_1''+\cdots+\psi x_p'',$$
 namely, the addition theorem with $h_i=\pm1$.  This result is
essentially Abel's main theorem on hyperelliptic integrals.

In his Paris manuscript, Abel addressed the two intermediate questions
for an arbitrary $f$.  However, his computations are more involved, and
his results, less definitive.  He found constancy holds when $\psi x$ is
of the form
\begin{equation}\label{equation:2adj}
   \psi x := \int^x_{x_0} \frac{h(x,y)}{\partial f/\partial y}\,dx
 \end{equation}
 where $\deg h\le\deg f-3$.  Also, $h$ must satisfy certain linear
conditions; namely, $h$ must vanish suitably everywhere $\partial
f/\partial y$ does on the curve $f=0$, at finite distance and at
infinity.  Abel took the maximum number of independent $h$ as the genus
$p$.

Furthermore, Abel found that there exists an $i\ge0$ such that
 \begin{equation}\label{equation:RR}
	\mu-\alpha= p-i.
 \end{equation}
 This equation does not contradict (\ref{equation:1ineq}) as the two
 $\alpha$'s differ; in (\ref{equation:1ineq}), the linear system of
 intersections is incomplete, whereas in (\ref{equation:RR}), the
 system is complete.

Abel's ideas have been clarified and completed over the course of time
through the efforts of many.  Doubtless, Riemann made the greatest
contribution in his revolutionary 1857 paper on Abelian functions.  In
his thesis of 1851, he had developed a way of extending complex analysis
to a multivalued function $y$ of a single variable $x$ by viewing $y$ as
a single-valued function on an abstract multisheeted covering of the
$x$-plane, the ``Riemann surface'' of $y$.  In 1857, he treated the case
where the surface is compact, and showed this case is precisely the case
where $y$ is algebraic.

Riemann defined the genus $p$ topologically, essentially as half the
first Betti number of the surface.  However, the term ``genus'' is not
Riemann's, but Clebsch's.  Clebsch introduced it in 1865 to signal his
aim of using $p$ in order to classify algebraic curves.  And he showed
that every curve of genus 0 is birationally equivalent to a line, and
every curve of genus 1, to a nonsingular plane cubic.

Also in 1865, Clebsch gave an algebro-geometric formula for the genus
$p$ of a plane curve: if the curve has degree $d$ and, at worst,
$\delta$ nodes and $\kappa$ cusps, then
	      $$p = (d-1)(d-2)/2-\delta-\kappa.$$
 The next year, Clebsch and Gordan employed this formula to prove the
birational invariance of $p$; they determined how $d$, $\delta$, and
$\kappa$ change.

Plainly, birationally equivalent curves have homeomorphic Riemann
surfaces, and so the same genus $p$.  But Clebsch was no longer
satisfied in just showing the consequences of Riemann's work.  He now
wanted to establish the theory of Abelian integrals on the basis of the
algebraic theory of curves as developed by Cayley, Salmon, and
Sylvester.  At the time, Riemann's theory was strange and suspect; there
was, as yet, no theory of manifolds, and no proof of the Dirichlet
principle. Clebsch's efforts led to a sea change in algebraic geometry,
which turned toward the study of birational invariants.

Riemann defined an integral to be of the ``first kind'' if it is finite
everywhere.  He proved these integrals form a vector space of dimension $p$.
Furthermore, each can be expressed in the form (\ref{equation:2adj})
provided the curve $f=0$ has, at worst, double points; if so, the linear
conditions on $h$ just require $h$ to vanish at these double points.  In
1874, Brill and M. Noether generalized this result to ordinary $m$-fold
points: $h$ must vanish to order $m-1$.  They termed such $h$
``adjoints.''  Meanwhile, starting with Kronecker in 1858 and Noether in
1871 and continuing through Muhly and Zariski in 1938, many algebraic
geometers developed corresponding ways of reducing the singularities of
a given plane curve by means of birational transformations.

Euler noted the integral $\int_0^xdx/\sqrt{1-x^4}$ has a ``modulus of
multivaluedness'' like that of the inverse trigonometric functions.
Abel noted an arbitrary Abelian integral has a similar ambiguity, but
viewed it as a sort of constant of integration, and avoided it by
keeping the domain small.  Riemann clarified the issue completely.  He
proved every integral of the first kind has $2p$ ``periods,'' which are
numbers that generate all possible changes in the value of the integral
arising from changes in the path of integration.

Riemann, in effect, did as follows.  He fixed a basis
$\psi_1x,\dotsc,\psi_px$ of the integrals of the first kind, and he
fixed a homology basis of $2p$ paths.  Then, inside the vector space
$\mathbb{C}^{\,p}$, he formed the lattice $\mathbb{L}$ generated by the
$2p$ corresponding $p$-vectors of periods.  And he proved the quotient
is a $p$-dimensional complex torus
	$$J:=\mathbb{C}^{\,p}/\mathbb{L}.$$
 Later $J$ was termed the ``Jacobian'' to honor Jacobi's work on
inversion.

Let $C$ be the curve $f=0$, or better, the associated Riemann surface.
Let $C^{(\mu)}$ be its $\mu$-fold symmetric product.  Riemann, in
effect, formed the following map:
 $$\Psi_\mu\colon C^{(\mu)}\to J\mbox{ given by }
  \Psi_\mu(x_1,\dotsc,x_\mu) = \biggl( \sum_{i=1}^\mu\psi_1x_i,
 	\dots, \sum_{i=1}^\mu \psi_px_i\biggr).
 $$
 This map $\Psi_\mu$ is rather important.  It has been called the
``Abel--Jacobi map'' and the ``Abel map.''  The latter name is
historically more correct and shorter, so better.

Riemann, in effect, studied the fibers of the Abel map $\Psi_\mu$.
He proved that, if two divisors $x_1+\dotsb+x_\mu$  and
$x_1'+\dotsb+x_\mu'$  are linearly equivalent, then 
	$$\Psi_\mu(x_1,\dotsc,x_\mu)
	 = \Psi_\mu(x_1',\dotsc,x_\mu').$$
 Riemann called this result ``Abel's Addition Theorem,'' and cited
Jacobi's 1832 proof of it in the hyperelliptic case.

The converse of this result holds too.  But Abel did not recognize
it, and it lies, at best, between the lines of Riemann's paper.  The
converse was first explicitly stated by Clebsch in 1864, and first
proved in full generality some time later by Weierstrass.  In 1913, in
Weyl's celebrated book on Riemann surfaces, Weyl combined the result and
its converse under the heading of Abel's Theorem.  Ever since then, most
mathematicians have done the same, even though Weyl explained it is not
historically correct to do so.

Together, the above result and its converse imply the fiber
$\Psi_\mu^{-1}\Psi_\mu(x_1,\dotsc,x_\mu)$ is the complete linear system
determined by $x_1+\dotsb+x_\mu$.  Its dimension is just Abel's
$\alpha$, the number of $x_i$ that can vary independently in the system.
Furthermore, in effect, Riemann rediscovered Abel's formula
(\ref{equation:RR}), and in 1864, Roch identified $i$ as the number of
independent adjoints vanishing on $x_1,\dotsc,x_\mu$.  In 1874, Brill
and Noether, inspired by Clebsch, gave the first algebro-geometric
treatment of the formula, whose statement they named the ``Riemann--Roch
Theorem.''

Finally, Riemann treated the inversion problem.  In effect, he proved
that the Abel map $\Psi_p$ is biholomorphic on a certain saturated
Zariski open subset $U\subset C^{(p)}$; namely, $U$ is the complement of
the image of $C^{(p-1)}$ in $C^{(p)}$ under the map
	$$(x_2,\dotsc,x_p)\mapsto (x_0,x_2,\dotsc,x_p).$$
 The inverse map $\Psi_pU\to U$ can be expressed using the coordinate
functions on $C^{(p)}$, so in terms of functions on $\Psi_pU$.  Since
$J:= \mathbb{C}^{\,p}/\mathbb{L}$, these functions can be lifted to an
open subset of $\mathbb{C}^{\,p}$, and then continued to meromorphic
functions on $\mathbb{C}^{\,p}$ with $2p$ periods.  Riemann termed these
special functions ``Abelian functions.''

Two years later, in 1859, Riemann proved that every meromorphic function
$F$ in $p$ variables has at most $2p$ independent period vectors.  Those
$F$ with exactly $2p$ soon became known as ``Abelian functions.''  They were
studied by many, including Weierstrass, Frobenius, and Poincar\'e.  In
particular, in 1869, Weierstrass observed that not every $F$ comes from a
curve.

Form the set $K$ of all ``Abelian functions'' whose group of periods
contains a given lattice $\mathbb{L}$ in $\mathbb{C}^{\,p}$ 
of rank $2p$.  It turns out  $K$
is a field of transcendence degree $p$ over $\mathbb C$.  Hence $K$ is
the field of rational functions on a $p$-dimensional projective
algebraic variety $A$, which is parameterized on a Zariski open set by
$p$ of them.  There are many such $A$, and all were called ``Abelian
varieties'' at first.  In 1919, Lefschetz proved there is a
distinguished $A$, whose underlying set can be identified with
$\mathbb{C}^{\,p}/\mathbb{L}$ in a natural way, and he restricted the
term ``Abelian variety'' to it.

Not only do the points of an Abelian variety $A$ form a group, but the
operations of adding and of inverting are given by polynomials.  Thus
$A$ is a complete algebraic group, or an ``Abelian variety'' in Weil's
sense of 1948.  Weil proved each such abstract Abelian variety is
commutative.  Earlier, in 1889, Picard had, in effect, proved this
commutativity in the case of a surface.

Every connected projective algebraic group is parameterized globally by
Abel\-ian functions with a common lattice of periods.  This fact was
proved by Picard for surfaces in 1889 and, assuming the group is
commutative, in any dimension in 1895.  His proof was completed at
certain points of analysis in 1903 by Painlev\'e.  Thus the two
definitions of Abel\-ian variety agree, Lefschetz's and Weil's; however,
Weil worked in arbitrary characteristic.

In the case of $C$ above, its Jacobian $J$ is thus an Abelian variety.
Moreover, $J$ is the quotient of $C^{(\mu)}$ for any $\mu\ge p$ by
linear equivalence.  So $J$ and $\Psi_\mu$ are defined by integrals, but
given by polynomials!  And addition on $J$ corresponds to addition of
divisors.  Therefore, $J$ and $\Psi_\mu$ can be constructed
algebro-geometrically just by forming the quotient.  Severi attributed
this construction to Castelnuovo.

In 1905, Castelnuovo generalized the construction to surfaces.  To set
the stage, he reviewed the case of curves, calling it very well known
(notissimo).  His work is a milestone in the history of irregular
surfaces, which began in 1868.

In 1868, Clebsch generalized Abel's formula (\ref{equation:2adj}) to a
surface $f(x,y,z)=0$ with ordinary singularities and no point at
infinity on the $z$-axis; in other words, $f=0$ is a general projection
of a smooth surface.  Clebsch showed that every double integral of the
first kind is of the form
	$$\int\!\!\int\frac{h(x,y,z)}{\partial f/\partial z}\,dx\,dy$$
 where $\deg h = \deg f -4$ and $h$ vanishes when $\partial f/\partial
z$ does.  The number of independent integrals became known as the
``geometric genus'' and denoted by $p_g$.

In 1870, M. Noether found an algebro-geometric proof that $p_g$ is a
birational invariant, as conjectured by Clebsch.  In 1871, Cayley found
a formula for the expected number of independent $h$, later called the
``arithmetic genus'' and denoted by $p_a$.  He observed that, if $f=0$
is a ruled surface over a base curve of genus $p$, then $p_a=-p$, but
$p_g=0$.  Later in 1871, Zeuthen used Cayley's formula to give an
algebro-geometric proof that $p_a$ too is a birational invariant.

In 1875, Noether explained the unexpected discrepancy between $p_g$ and
$p_a$: the vanishing conditions on $h$ need not be independent; in any
case, $p_g\ge p_a$.  He noted that, if the surface $f=0$ in 3-space is
smooth or rational, then $p_g=p_a$.  It was expected that equality
usually holds; so when it did, $f=0$ was termed ``regular.''  The
difference $p_g-p_a$ gives a quantitative measure of the failure of
$f=0$ to be regular; so $p_g-p_a$ was termed its ``irregularity.''

In 1884, Picard studied, on the surface $f=0$, simple integrals
	$$\int P(x,y,z)\,dx + \int Q(x,y,z)\,dy$$
 that are closed, or $\partial P/\partial y = \partial Q/\partial x$;
they became known as ``Picard integrals.''  And $q$ was used to denote
the number of independent Picard integrals of the first kind.  Picard
noted that, if $f=0$ is smooth, then $q=0$.  In 1894, Humbert proved
that, if $q=0$, then every algebraic system of curves is contained in a
linear system.

Inspired by Humbert's result, in 1896, Castelnuovo proved that, if
$p_g-p_a=0$, then every algebraic system of curves on $f=0$ is contained
in a linear system under a certain restriction.  In 1899, Enriques
removed the restriction.  For a modern version of this result and of its
converse, which together characterize regular surfaces, see
Exercises~\ref{ex:q=0} and \ref{ex:Enriques}.

In 1897, Castelnuovo fixed a linear system of curves on the surface
$f=0$, and studied the ``characteristic'' linear system cut out on a
general member by the other members.  Let $\delta$ be the amount, termed
the ``deficiency,'' by which the dimension of the characteristic system
falls short of the dimension of its complete linear system.  Castelnuovo
proved that
  $$\delta\le p_g-p_a,$$
 and equality holds for the system cut out by the surfaces of suitably
high degree.

In February 1904, Severi extended Castelnuovo's work.  Severi took a
complete algebraic system of curves, without repetition, on the surface
$f=0$, say with parameter space $\Sigma$ of dimension $R$.  Let
$\sigma\in \Sigma$ be a general point, and $C_\sigma$ the corresponding
curve. Say $C_\sigma$ moves in a complete linear system of dimension
$r$.  Now, to each tangent direction at $\sigma\in \Sigma$, Severi
associated, in an injective fashion, a member of the complete
characteristic system on $C_\sigma$.  Thus he got an $R$-dimensional
``characteristic'' linear subsystem.  The complete system has dimension
$r+\delta$.  Hence, 
 $$R\le r+p_g-p_a.$$

A few months later, Enriques and, shortly afterward, Severi gave proofs
that, if $C_\sigma$ is sufficiently positive, then equality holds above;
in other words, then the characteristic system of $\Sigma$ is complete.
Both proofs turned out to have serious gaps, as Severi himself observed
in 1921.  Meanwhile, in 1910, Poincar\'e gave an analytic construction
of a family with $R=p_g-p_a$ and $r=0$.  It follows formally, by means
of the Riemann--Roch Theorem for surfaces, that whenever $C_\sigma$ is
sufficiently positive, the characteristic system is complete.  After
Severi's criticism, it became a major open problem to find a purely
algebro-geometric treatment of this issue.  But see
Corollary~\ref{cor:Poincare} and Remark~\ref{rmk:charsys}: the solution
finally came forty years later with Grothendieck's systematic use of
nilpotents!

In mid January 1905, Severi proved that $p_g-p_a\ge q$ and that
$p_g-p_a=b-q$ where $b$ is the number of independent Picard integrals of
the second kind, which is equal to the first Betti number.
Simultaneously and independently, Picard too proved that $p_g-p_a=b-q$.

A week later and more fully that May and June, Castelnuovo took the last
step in this direction.  He fixed $C_\sigma$ sufficiently positive, and
formed the quotient, $P$ say, of $\Sigma$ modulo linear equivalence.  So
$P$ is projective, and
	$\dim P=p_g-p_a$.
 Furthermore, since two sufficiently positive curves sum to a third, it
follows that $P$ is independent of the choice of $C_\sigma$, and is a
commutative group variety.  Hence $P$ is an Abelian variety by the
general theorem of Picard, completed by Painlev\'e, mentioned above.
Hence $P$ is parameterized by $\dim P$ Abelian functions.  Castelnuovo
proved they induce independent Picard integrals on $f=0$.  Therefore,
$p_g-p_a\le q$.  Thus Castelnuovo obtained the Fundamental Theorem of
Irregular Surfaces: \begin{equation}\label{equation:FTIS}
 	\dim P=p_g-p_a=q =b/2.
 \end{equation}
 For a modern discussion of the result, see Remark~\ref{rmk:Igusa} and
Exercise~\ref{ex:q=0}.  In 1905, the term ``Abelian variety'' was not
yet in use, and so, naturally enough, Castelnuovo termed $P$ the
``Picard variety'' of the surface $f=0$.

Picard also studied Picard integrals of the third kind on the surface
$f=0$.  In 1901, he proved that there is a smallest integer $\varrho$
such that any $\varrho+1$ curves are the logarithmic curves of some such
integral.  On the basis this result, in 1905, Severi proved, in effect,
that $\varrho$ is the rank of the group of all curves modulo algebraic
equivalence.  In 1952, N\'eron proved the result in arbitrary
characteristic.  So the group is now called the ``N\'eron--Severi
group,'' and  $\varrho$ is called the ``Picard number.''

In 1908 and 1910, Severi studied, in effect, the torsion subgroup of the
N\'eron--Severi group, notably proving it is finite.  In 1957, Matsusaka
proved this finiteness in arbitrary characteristic.  However, there is
no special name for this subgroup or for its order.  For more about them
and $\varrho$, see Corollary (\ref{cor:torgp}) and Remark
(\ref{rmk:Ptaufin}).

The impetus to work in arbitrary characteristic came from developments
in number theory.  In 1921, E. Artin developed, in his thesis, an
analogue of the Riemann Hypothesis, in effect, for a hyperelliptic curve
over a prime field of odd characteristic.  In 1929, F. K. Schmidt
generalized Artin's work to all curves over all finite fields, and recast
it in the geometric style of Dedekind and Weber.  In 1882, they had
viewed a curve as the set of discrete valuation rings in a finitely
generated field of transcendence degree 1 over $\mathbb C$, and they had
given an abstract algebraic treatment of the Riemann--Roch Theorem.
Schmidt observed that their treatment works with little change in
arbitrary characteristic, and he used the Riemann--Roch Theorem to prove
that Artin's zeta function satisfies a natural functional equation.

In 1936, Hasse proved Artin's Riemann hypothesis in genus 1 using an
analogue of the theory of elliptic functions.  Then he and Deuring noted
that to extend the proof to higher genus would require developing a
theory of correspondences between curves analogous to that developed by
Hurwitz and others.  This work inspired Weil to study the fixed points
of the Frobenius correspondence, and led to his announcement in 1940 and
to his two great proofs in 1948 of Artin's Riemann hypothesis for the
zeta function of an arbitrary curve and also to his proof of the
integrality of his analogue of Artin's $L$-functions of 1923 and 1930.

First, in 1946, Weil carefully rebuilt the foundation of algebraic
geometry from scratch.  Following in the footsteps of E. Noether, van
der Waerden, and Schmidt, Weil took a variable coefficient field of
arbitrary characteristic inside a fixed algebraically closed coordinate
field of infinite transcendence degree.  Then he formed ``abstract''
varieties by patching pieces of projective varieties, and said when
these varieties are ``complete.''  Finally, he developed a calculus of
cycles.

In 1948, Weil published two exciting monographs.  In the first, he
reproved the Riemann--Roch theorem for (smooth complete) curves, a 
theorem he regarded as fundamental (see \cite[I, p.~562, top; II,
p.~541, top]{We79}).  Then he developed an elementary theory of
correspondences between curves, which included Castel\-nuovo's theorem of
1906 on the positive definiteness of the equivalence defect of a
correspondence.  Of course, Castelnuovo's proof was set over $\mathbb
C$, but ``its translation into abstract terms was essentially a routine
matter once the necessary techniques had been created,'' as Weil put it
in his 1954 ICM talk.  Finally, Weil derived the Riemann hypothesis.

In the second monograph, Weil established the abstract theory of Abelian
varieties.  He constructed the Jacobian $J$ of a curve $C$ of genus $p$
by patching together copies of an open subset of the symmetric product
$C^{(p)}$.  Then taking a prime $l$ different from the characteristic,
he constructed, out of the points on $J$ of order $l^n$ for all $n\ge1$,
an $l$-adic representation of the ring of correspondences, equivalent to
the representation on the first cohomology group of $C$.  Finally, he
proved the positive definiteness of the trace of this representation,
reproved the Riemann hypothesis for the zeta function, and completed the
proof of his analogue of Artin's conjectured integrality for
$L$-functions of number fields.

Weil left open two questions: Do a curve and its Jacobian have the same
coefficient field?  Is every Abelian variety projective?  Both questions
were soon answered in the affirmative by Chow and Matsusaka.  However,
there has remained some general interest in constructing nonprojective
varieties and in finding criteria for projectivity.  Furthermore, Weil
was led in 1956 to study the general question of descent of the
coefficient field, and this work in turn inspired Grothendieck's general
descent theory, which he sketched in \cite[no.~236]{FGA}.

In 1949, Weil published his celebrated conjectures about the zeta
function of a variety of arbitrary dimension.  Weil did not explicitly
explain these conjectures in terms of a hypothetical cohomology theory,
but such an explanation lies between the lines of his paper.
Furthermore, it was credited to him explicitly in Serre's 1956 ``Mexico
paper'' \cite[p.~24]{Sr56} and in Grothendieck's 1958 ICM talk.

In his talk, Grothendieck announced that he had found a new approach to
developing the desired ``Weil cohomology.''  He wrote: ``it was
suggested to me by the connections between sheaf-theoretic cohomology
and cohomology of Galois groups on the one hand, and the classification
of unramified coverings of a variety on the other \dots, and by Serre's
idea that a `reasonable' algebraic principal fiber space \dots, if it is
not locally trivial, should become locally trivial on some covering
unramified over a given point.''  This is the announcement of
Grothendieck topology.

In 1960, Grothendieck and Dieudonn\'e \cite[p.~6]{EGAI} listed the
titles of the chapters they planned to write.  The last one, Chapter
XIII, is entitled ``Cohomologie de Weil.''  The next-to-last is
entitled ``Sch\'emas ab\'eliens et sch\'emas de Picard.''  Earlier, at
the end of his 1958 ICM talk, Grothendieck had listed five open
problems; the fifth is to construct the Picard scheme.

In 1950, Weil published a remarkable note on Abelian varieties.  For
each complete normal variety $X$ of any dimension in any characteristic,
he said there ought to be two associated Abelian varieties, the
``Picard'' variety $P$ and the ``Albanese'' variety $A$, with certain
properties, discussed just below.  He explained he had complete proofs
for smooth complex $X$, and ``sketches'' in general.  Soon all was
proved.

Weil's sketches rest on two criteria for linear equivalence, developed
in 1906 by Severi and reformulated in the 1950 note by Weil.  He
announced proofs of them in 1952, and published the details in 1954.
For some more information, see \cite[p.~120]{Za35} and
Remark~\ref{rmk:RamSam}.  In Weil's commentaries on his '54 paper, he
wrote: ``Ever since 1949, I considered the construction of an algebraic
theory of the Picard variety as the task of greatest urgency in abstract
algebraic geometry.''

The properties are these.  First, $P$ parameterizes the linear
equivalence classes of divisors on $X$.  And there exists a map $X\to A$
that is ``universal'' in the sense that every map from $X$ to an Abelian
variety factors through it.  In his commentaries on the note, Weil
explained that $P$ had been introduced and named by Castelnuovo; so
historically speaking, it would be justified to name $P$ after him, but
it was better not to tamper with common usage.  By contrast, Weil chose
to name $A$ after Albanese in order to honor his work in 1934 viewing
$A$ as a quotient of a symmetric power of $X$, although $A$ had been
introduced and studied in 1913 by Severi.

Second, if $X$ is an Abelian variety, then $X$ is equal to the Picard
variety of $P$; so each of $X$ and $P$ is the ``dual'' of the other.  If
$X$ is arbitrary, then $A$ and $P$ are dual Abelian varieties; in fact,
the universal map $X\to A$ induces the canonical isomorphism from the
Picard variety of $A$ onto $P$.  If $X$ is a curve, then both $A$ and
$P$ coincide with the Jacobian of $X$, and the universal map $X\to A$ is
just the Abel map; in other words, the Jacobian is ``autodual.''  This
autoduality can be viewed as an algebro-geometric statement of Abel's
theorem and its converse for integrals of the first kind.  For some more
information, see Remarks~\ref{rmk:Ablsch}--\ref{rmk:Jac}.

In 1951, Matsusaka gave the first algebraic construction of $P$.
However, he had to extend the ground field because he applied Weil's
results: one of the equivalence criteria, and the construction of the
Jacobian.  Both applications involve the ``generic curve,'' which is the
section of $X$ by a generic linear space of complementary dimension.  In
1952, Matsusaka gave a different construction; it does not require
extending the ground field, but requires $X$ to be smooth.

Both constructions are like Castelnuovo's in that they begin by
constructing a complete algebraic system of sufficiently positive
divisors, and then form the quotient modulo linear equivalence.  To
parameterize the divisors, Matsusaka used the theory of ``Chow
coordinates,'' which was developed by Chow and van der Waerden in 1938
and refined by Chow contemporaneously.  In 1952, Matsusaka also used
this theory to form the quotient.  In the same paper, he gave the first
construction of $A$, again using the Jacobian of the generic curve, but
he did not relate $A$ and $P$.

In 1954, Chow published a construction of the Jacobian similar to
Matsusaka's second construction of $P$.  Chow had announced it in 1949,
and both Weil and Matsusaka had referred to it in the meantime.  In
1955, Chow constructed $A$ and $P$ by a new procedure; he took the
``image'' and the ``trace'' of the Jacobian of a generic curve.
Moreover, he showed that the universal map $X\to A$ induces an
isomorphism from the Picard variety of $A$ onto $P$.

In a course at the University of Chicago, 1954--55, Weil gave a more
complete and elegant treatment, based on the ``see-saw principle,'' which
he adapted from Severi, and on his own Theorem of the Square and Theorem
of the cube.  This treatment became the core of Lang's 1959 book,
``Abelian.Varieties.''  The idea is to construct $A$ first using the
generic curve, and then to construct $P$ as a quotient of $A$ modulo a
finite subgroup; thus there is no need for Chow coordinates.

In 1959 and 1960, Nishi and Cartier independently established the
duality between $A$ and $P$ in full generality.

Between 1952 and 1957, Rosenlicht published a remarkable series of
papers, which grew out of his 1950 Harvard thesis.  It was supervised by
Zariski, who had studied Abelian functions and algebraic geometry with
Castelnuovo, Enriques, and Severi in Rome from 1921 to 1927.  Notably
Rosenlicht generalized to a curve with arbitrary singularities the
notions of linear equivalence and differentials of the first kind.  Then
he constructed a ``generalized Jacobian'' over $\mathbb{C}$ by
integrating and in arbitrary characteristic by patching.  It is not an
Abelian variety, but an extension of the Jacobian of the normalized
curve by an affine algebraic group.

Rosenlicht cited Severi's 1947 monograph, ``Funzioni Quasi Abeliane,''
where the generalized Jacobian was discussed for the first time, but
only for curves with double points.  In turn, Severi traced the history
of the corresponding theory of quasi-Abelian functions back to Klein,
Picard, Poincar\'e, and Lefschetz.

In 1956, Igusa established the compatibility of specializing a curve
with specializing its generalized Jacobian in arbitrary characteristic
when the general curve is smooth and the special curve has at most one
node.  Igusa explained that, in 1952, N\'eron had studied the total
space of such a family of Jacobians, but had not explicitly analyzed the
special fiber.  Igusa's approach is, in spirit, like Castelnuovo's,
Chow's, and Matsusaka's before him and Grothendieck's after him.  But
Grothendieck went considerably further: he proved compatibility with
specialization for a family of varieties of arbitrary dimension with
arbitrary singularities, both in equicharacteristic and in mixed.

In 1960, Chevalley constructed a Picard variety for any normal variety
$X$ using locally principal divisorial cycles.  Cartier had already
focused on these cycles in his 1958 Paris thesis.  But Chevalley said he
would call them simply ``divisors,'' and we follow suit, although they
are now commonly called ``Cartier divisors.''

First, Chevalley constructed a ``strict'' Albanese variety; it is
universal for regular maps (morphisms) into Abelian varieties.  Then he
took its Picard variety to be that of $X$.  He noted his Picard and
Albanese varieties need not be equal to those of a desingularization of
$X$.  By contrast, Weil's $P$ and $A$ are birational invariants, and his
universal map $X\to A$ is a rational map, which is defined wherever $X$
is smooth.  In 1962, Seshadri generalized Chevalley's construction to a
variety with arbitrary singularities, thus recovering Rosenlicht's
generalized Jacobian.

Back in 1924, van der Waerden initiated the project of rebuilding the
whole foundation of algebraic geometry on the basis of commutative
algebra.  His goal was to develop a rigorous theory of Schubert's
enumerative geometry, as called for by Hilbert's fifteenth problem.  Van
der Waerden drew on Elimination Theory, Ideal Theory, and Field Theory
as developed in the schools of Kronecker, of Dedekind, and of Hilbert.
Van der Waerden originated, notably, the algebraic notion of
specialization as a replacement for the topological notion of
continuity.

In 1934, as Zariski wrote his book~\cite{Za35}, he lost confidence in
the clarity, precision, and completeness of the algebraic geometry of
his Italian teachers.  He spent a couple of years studying the algebra
of E. Noether and Krull, and aimed to reduce singularities rigorously.
He introduced three algebraic tools: normalization, valuation theory,
and completion.  He worked extensively with the rings obtained by
localizing affine coordinate domains at arbitrary primes over arbitrary
fields.  And, in 1944, he put a topology on the set of all valuation
rings in a field of algebraic functions, and used the property that any
open covering has a finite subcovering.

In 1949, Weil saw that the ``Zariski topology'' can be put on his
abstract varieties, simplifying the old exposition and suggesting the
construction of new objects, such as locally trivial fiber spaces.  In
his paper of 1950 on Abelian varieties, he noted that line bundles
correspond to linear equivalence classes of divisors, and predicted that
line bundles would play a role in the theory of quasi-Abelian functions.

In 1955, Serre provided abstract algebraic geometry with a very powerful
new tool: sheaf cohomology.  Given a variety equipped with the Zariski
topology, he assembled the local rings into the stalks of a
``structure'' sheaf.  Then he developed a cohomology theory of coherent
sheaves, analogous to the one that he and Cartan and Kodaira and Spencer
had just developed in complex analytic geometry, and had so successfully
applied to establish and to generalize previous work on complex
algebraic varieties.

About the same time, a general theory of abstract algebraic geometry was
developed by Chevalley.  He did not use sheaves and cohomology, but did
work with what he called ``schemes,'' obtained by patching ``affine
schemes''; an affine scheme is the set of rings obtained by localizing a
finitely generated domain over an arbitrary field.  Nevertheless, he
soon returned to a more traditional theory of ``varieties'' when he
worked on the theory of algebraic groups.

In January 1954, Chevalley lectured on schemes at Kyoto University.  His
lectures inspired Nagata in \cite{Na56} to generalize the theory by
replacing the coefficient field by a Dedekind domain.  But Nagata used
Zariski's term ``model,'' not Chevalley's term ``scheme.''  Earlier, at
the 1950 ICM, Weil had recalled Kronecker's dream of an algebraic
geometry over the integers; however, Nagata did not cite Weil's talk,
and likely was not motivated by it.

In the fall of 1955, Chevalley lectured on schemes over fields at the
S\'eminaire Cartan--Chevalley, and Grothendieck was there.  By February
1956 (see \cite[p.~32]{CS01}), he was patching the spectra of arbitrary
Noetherian rings, and studying the cohomology of Cartier's
``quasi-coherent'' sheaves.  There is good reason for the added generality:
nilpotents allow better handling of higher-order infinitesimal
deformations, of inseparability in positive characteristic, and of
passage to the completion; quasi-coherent sheaves have the technical
convenience of coherent sheaves without their cumbersome finiteness.  By
October 1958 (see \cite[p.~63]{CS01}), Grothendieck and Dieudonn\'e had
begun the gigantic program of writing EGA\emdash rebuilding once again
the foundation of algebraic geometry in order to provide a more flexible
framework, more powerful methods, and a more refined theory.

Also in 1958, there appeared two other papers, which discuss objects
similar to Chevalley's schemes: K\"ahler published a 400 page
foundational monograph, which introduced general base changes, and Chow
and Igusa published a four-page note, which proved the K\"unneth Formula
for coherent sheaves.  The two works are mentioned briefly in
 \cite[p.~101]{CS01}, in \cite[p.~8]{EGAI} and in \cite[p.~6]{EGAG};
 however, they seem to have had little or no influence on Grothendieck.

Finally, in 1961--1962, Grothendieck constructed the Hilbert scheme
and the Picard scheme.  The construction is a technical masterpiece,
showing the tremendous power of Grothendieck's new tools.  In
particular, a central role is played by the theory of flatness.  It was
introduced by Serre in 1956 as a formal device for use in comparing
algebraic functions and analytic functions.  Then Grothendieck developed
the theory extensively, for he recognized that flatness is the technical
condition that best expresses continuity across a family.

Grothendieck \cite[p.~221-1]{FGA} saw Hilbert schemes as ``destined to
replace'' Chow coordinates.  However, he \cite[p.~221-7]{FGA} had to
appeal to the theory of Chow coordinates for a key finiteness result: in
projective space, the subvarieties of given degree form a bounded
family.  A few years later, Mumford \cite[Lects.~14, 15]{Mm66} gave a
simple direct proof of this finiteness; his proof introduced an
important new tool, now known as ``Castelnuovo--Mumford'' regularity.
In a slightly modified form, this tool plays a central role in the
proofs of the finiteness theorems for the Picard scheme, which are
addressed in Chapter 6 below.

In spirit, Grothendieck's construction of the Picard scheme is 
like Castelnuovo's and Matsusaka's.  He began with the component
$\Sigma$ of the Hilbert scheme determined by a sufficiently positive
divisor.  Then he formed the quotient; in fact, he did so twice for
diversity.  First, he used ``quasi-sections''; second, and more
elegantly, he used the Hilbert scheme of $\Sigma$.

Grothendieck's definition of the two schemes is yet a greater
contribution than his construction of them.  He defined them by their
functors of points.  These schemes are universal parameter spaces; so
they receive a map from a scheme $T$ just when $T$ parameterizes a
family of subschemes or of invertible sheaves, respectively, and this
map is unique.

  What is a universal family?  The answer seems obvious when we use
schemes.  But Chow coordinates parameterize positive cycles, not
subschemes.  And even in the analytic theory of the Picard variety $P$,
there was some question about the sense in which $P$ parameterizes
divisor classes.  Indeed, the American Journal of Math., {\bf 74}
(1952), contains three papers on $P$.  First, Igusa constructed $P$, but
left universality unsettled.  Then Weil and Chow settled it with
different arguments.

A functor of points, or ``representable functor,'' is not an arbitrary
contravariant functor from schemes to sets.  It is determined locally,
so is a sheaf.  But it suffices to represent a sheaf locally, as the
patching is determined.  Thus to construct the Hilbert scheme, the first
step is to check that the Hilbert functor is, in fact, a sheaf for the
Zariski topology, that is, a ``Zariski sheaf.''

The naive Picard functor is not a Zariski sheaf.  So the first step is
to localize it, or form the associated sheaf.  This time, the Zariski
topology is not fine enough.  However, a representable functor is an
fpqc sheaf by a main theorem of Grothendieck's descent theory.  In
practice, it is enough to localize for the \'etale Grothendieck topology
or for the fppf, and these localizations are more convenient to work
with.  The localizations of the Picard functor are discussed in Chapter
2.

The next step is to cover the localized Picard functor by representable
Zariski open subfunctors.  This step is elementary.  But it is
technically involved, more so than any other argument in these notes.
It is carried out in the proof of the main result,
Theorem~\ref{th:main}.  Each subfunctor is represented by a quotient of
an open subscheme of the Hilbert scheme.  Thus the Picard scheme is
constructed.

In sum, Grothendieck's method of representable functors is like
Descartes's method of coordinate axes: simple, yet powerful.  Here is
one hallmark of genius!


In the notes that follow, our primary aim is to develop in detail most
of Grothendieck's original theory of the Picard scheme basically by
filling in his sketch in \cite{FGA}.  Our secondary aim is to review in
brief much of the rest of the theory developed by Grothendieck and by
others.  We review the secondary material in a series of scattered
remarks.  The remarks refer to each other and to the primary discussion,
but the primary discussion never refers to the remarks.  So the remarks
may be safely ignored in a first reading.

Notably, the primary discussion does not develop Grothendieck's
method of ``relative representability.''  Indeed, the details would take
us too far afield.  On the other hand, were we to use the method, we
could obtain certain existence theorems and finiteness theorems in
greater generality by reducing to the cases that we do handle.
Consequently, in Sections~4--6, a number of results just assume the
Picard scheme exists, rather than assume hypotheses guaranteeing it
does.  However, we do discuss the method and its applications in
Remark~\ref{rk:exist} and in other remarks.

These notes also contain many exercises, which call for working
examples and constructing proofs.  Unlike the remarks, these exercises
are an integral part of the primary discussion, which not only is
enhanced by them, but also is based in part on them.  Furthermore, the
exercises are designed to foster comprehension.  The answers involve no
new concepts or techniques.  The exercises are meant to be easy; if a
part seems to be hard, then some review and reflection may be in
order.  However, all the answers are worked out in detail in
Appendix~\ref{sc:A}.

These notes assume familiarity with the basic algebraic geometry
developed in Chapters II and III of Hartshorne's popular
textbook~\cite{Ha83}, and assume familiarity, but to a lesser extent,
with the foundational material developed in Grothendieck and
Dieudonn\'e's monumental reference books \cite{EGAI} to \cite{EGAG}.  In
addition, these notes assume familiarity with basic Grothendieck
topology, descent theory, and Hilbert-scheme theory, such as that
explained on pp.~129--147, 199--201, and 215--221 in Bosch,
L\"utkebohmert, and Raynaud's welcome survey book~\cite{BLR}; this
material and more was introduced by Grothendieck in three Bourbaki
talks, nos.~190, 212, and 221, which were reprinted in \cite{FGA} and
are still worth reading.  Of course, when specialized results are used
below, precise references are provided.

Throughout these notes, we work only with locally Noetherian schemes,
just as Grothendieck did in \cite{FGA}.  Shortly afterward,
Grothendieck promoted the elimination of this restriction, through a
limiting process that reduces the general case to the Noetherian case.
Ever since, it has been common to make this reduction.  However, the
process is elementary and straightforward.  Using it here would only be
distracting.

Given two locally Noetherian schemes lying over a third, the (fibered)
product is not necessarily locally Noetherian.  Consequently, there are
minor technical difficulties in working with the fpqc topology on the
category of locally Noetherian schemes.  On the other hand, in practice,
there is no need to use the fpqc topology.  Therefore, its use has been
eliminated from these notes.

Throughout, we work with a  {\bf fixed} separated map of finite type
	$$f\:X\to S.$$
  For convenience, when given an $S$-scheme $T$, we set
	$$X_T:=X\x_ST$$
 and denote the projection by $f_T\:X_T\to T$.  Also, when given a
$T$-scheme $T'$ and given quasi-coherent sheaves $\mc N$ on $T$ and $\mc
M$ on $X_T$, we denote the pullback sheaves by $\mc N|T'$ or $\mc
N_{T'}$ and by $\mc M|X_{T'}$ or $\mc M_{T'}$.

Given an $S$-scheme $P$, we call an $S$-map $T\to P$ a $T$-{\it point\/}
of $P$, and we denote the set of all $T$-points by $P(T)$.  As $T$
varies, the sets $P(T)$ form a contra\-variant functor on the category
of $S$-schemes, called the {\it functor of points\/} of $P$.

Section 2 introduces and compares the four common relative Picard
functors, the likely candidates for the functor of points of the Picard
scheme.  They are simply the functor $T\mapsto\Pic(X_T)\big/\Pic(T)$ and
its associated sheaves in the Zariski topology, the \'etale topology,
and the fppf topology.  Section 3 treats relative effective (Cartier)
divisors on $X/S$ and the relation of linear equivalence.  We prove
these divisors are parameterized by an open subscheme of the Hilbert
scheme of $X/S$.  Furthermore, we study the subscheme parameterizing
the divisors whose fibers are linearly equivalent, and prove it is of
the form $\I P(\mc Q)$ where $\mc Q$ is a certain coherent sheaf on $S$.

Section 4 begins the study of the Picard scheme $\IPic_{X/S}$ itself.
Notably, we prove Grothendieck's main theorem: $\IPic_{X/S}$ exists when
$X/S$ is projective and flat and its geometric fibers are integral.
Then we work out Mumford's example showing the necessity of the
integrality hypothesis.  Section 5 studies $\IPicz_{X/S}$, which is the
union of the connected components of the identity element of the fibers
of $\IPicz_{X/S}$.  In particular, we compute the tangent space at the
identity of each fiber.  It is remarkable how much we can prove {\it
formally\/} about $\IPicz_{X/S}$.  Section 6 proves the two deeper
finiteness theorems.  They concern $\IPict_{X/S}$, which consists of the
points with a multiple in $\IPicz_{X/S}$, and $\IPicp_{X/S}$, whose
points represent the invertible sheaves with a given Hilbert polynomial
$\phi$.

Finally, there are two appendices.  Appendix A contains detailed answers
to all the exercises.  Appendix B develops basic divisorial intersection
theory, which is used freely throughout Section 6.  The treatment is
short, simple, and elementary.

\section{The several Picard functors}\label{sc:Pfctrs}
\renewcommand{\theequation}{\thethm.\arabic{equation}}
Our first job is to identify a likely candidate for the functor of
points of the Picard scheme.  In fact, there are several reasonable such
Picard functors, and each one is more likely to be representable than
the preceding.  In this section, they all are formally introduced and
compared.

\begin{dfn}\label{df:aPf}
 The {\it absolute Picard functor} $\Pic_X$ is the functor from  the
 category of (locally Noetherian) $S$-schemes $T$ to the category of
abelian groups defined by the formula
	$$\Pic_X(T):=\Pic(X_T).$$
\end{dfn}

The absolute Picard functor  is a ``prepared presheaf'' in this
sense: given any family of  $S$-schemes $T_i$, we have
	$$\Pic_X\bigl(\textstyle\coprod T_i\bigr)=\prod\Pic_X(T_i).$$
 Hence, given a covering family $\{T_i\to T\}$ in the Zariski topology,
the \'etale topology, or any other Grothendieck topology on the category
of $S$-schemes, there is no harm, when we consider the sheaf associated
to $\Pic_X$, in working simply with the single map $T'\to T$ where
$T':=\coprod T_i$.  And doing so lightens the notation, making for
easier reading. Therefore, we do so throughout, calling $T'\to T$ simply
a {\it covering\/} in the given topology.

The absolute Picard functor $\Pic_X$ is never a separated presheaf in
the Zariski topology.  Indeed, take an $S$-scheme $T$ that carries an
invertible sheaf $\mc N$ such that $f^*_T\mc N$ is nontrivial.  (For
example, take $T:=\I P^1_X$ and $\mc N:=\mc O_{\!T}(1)$.  Then the diagonal
map $X\to X\x X$ induces a section $g$ of $f_T$; that is, $gf_T=1_T$.
Hence $f^*_T\mc N$ is nontrivial.)  Now, there exists a Zariski covering
$T'\to T$ such that the pullback $\mc N|T'$ is trivial; here, $T'$ is
simply the disjoint union of the subsets in a suitable ordinary open
covering of $T$.  So the pullback $f^*_{T'}\mc N|X_{T'}$ is trivial too.
Thus $\Pic_X(T)$ has a nonzero element whose restriction is zero in
$\Pic_X(T')$.

According to descent theory, every representable functor is a sheaf in
the Zariski topology \emdash in fact, in the \'etale and fppf
topologies as well.  Therefore, the absolute Picard functor $\Pic_X$ is never
representable.  So, in the hope of obtaining a representable functor
that differs as little as possible from $\Pic_X$, we now define a
sequence of three successively more promising relative functors.

\begin{dfn}\label{df:Pfs}
  The {\it relative Picard functor} $\Pic_{X/S}$ is defined by
	$$\Pic_{X/S}(T):=\Pic(X_T)\big/\Pic(T).$$
 Denote its associated sheaves in the Zariski, \'etale, and fppf
topologies  by
	$$\Pic_{(X/S)\zar},\ \Pic_{(X/S)\et},\
	\Pic_{(X/S)\fppf}.$$
 \end{dfn}

We now have a sequence of five Picard functors, and each one maps
naturally into the next.  So each one maps naturally into any of its
successors.  If the latter functor is one of the three just displayed,
then it is a sheaf in the indicated topology; in fact, it is the sheaf
associated to any one of its predecessors, and the map between them is
the natural map from a presheaf to its associated sheaf, as is easy to
check.  In particular, the three displayed sheaves are the sheaves
associated to the absolute Picard functor $\Pic_X$, as well as to the
relative Picard functor $\Pic_{X/S}$.

Since $\Pic_{X/S}$ is not a priori a sheaf, it is remarkable that it is
representable so often in practice.

Note that, for every $S$-scheme $T$, each $T$-point of
$\Pic_{(X/S)\fppf}$, or element of $\Pic_{(X/S)\fpqc}(T)$, is
represented by an invertible sheaf $\mc L'$ on $X_{T'}$ for some
fppf-covering $T'\to T$.    Moreover, there must
be an fppf-covering $T''\to T'\x_T T'$ such that the two pullbacks of
$\mc L'$ to $X_{T''}$ are isomorphic \emdash or to put it informally, the
restrictions of $\mc L'$ must agree on a covering of the overlaps.

Furthermore, a second such sheaf $\mc L_1$ on $X_{T_1}$ represents the
same $T$-point if and only if there is an fppf-covering $T_1'\to T_1\x_T
T'$ such that the pullbacks of $\mc L'$ and $\mc L_1$ to $X_{T_1'}$ are
isomorphic.  Technically, this condition includes the preceding one,
which concerns the case where $T_1=T'$ and $\mc L_1=\mc L'$ since $\mc
L'$ must represent the same $T$-point as itself.  Of course, similar
considerations apply to the Picard functors for the Zariski and \'etale
topologies as well.

\begin{ex}\label{ex:Alr}
 Given an $S$-scheme $T$ of the form $T=\Spec(A)$ where $A$ is a local
ring, show that the natural maps are isomorphisms
	$$\Pic_X(A)\risom\Pic_{X/S}(A)\risom\Pic_{(X/S)\zar}(A)$$
 where $\Pic_X(A):=\Pic_X(T)$, where $\Pic_{X/S}(A):=\Pic_{X/S}(T)$,
and so forth.

Assume $A$ is Artin local with algebraically closed residue field.  Show
	$$\Pic_X(A)\risom\Pic_{(X/S)\et}(A).$$

Assume $A$ is an algebraically closed  field $k$.  Show
	$$\Pic_X(k)\risom\Pic_{(X/S)\fppf}(k).$$
  \end{ex}

\begin{ex}\label{ex:Pfs}
 Show that the natural map
	$$\Pic_{(X/S)\zar}\to \Pic_{(X/S)\et}$$
 need not be an isomorphism.  Specifically, take $X$ to be the following
 curve  in the real projective plane:
	$$X:u^2+v^2+w^2=0 \text{ in }\I P^2_\bb R.$$
 Then $X$ has no $\bb R$-points, but  over the complex numbers $\bb C$,
 there is an isomorphism
	$$\vf\:X_{\bb C}\risom\I P^1_{\bb C}.$$
 Show that $\vf^*\mc O(1)$ defines an element of $\Pic_{(X/\bb R)\et}(\bb
R)$, which is not in the image of $\Pic_{(X/\bb R)\zar}(\bb R)$.
 \end{ex}

The main result of this section is the following comparison theorem.
It identifies two useful conditions: the first guarantees that the first
three relative functors can be viewed as subfunctors of the fourth;
together, the two conditions guarantee that all four functors coincide.
The second condition has three successively weaker forms.  Before we can
prove the theorem, we must develop some theory.

\begin{thm}[Comparison]\label{th:cmp}
 Assume $\mc O_{\!S}\risom f_*\mc O_{\!X}$ holds universally; that is, for any
$S$-scheme $T$, the comorphism of $f_T$ is an isomorphism, $\mc O_{\!T}\risom
f_{\smash{T}*}\mc O_{\!X_T}$.
 \smallskip

\tu{1.}\enspace
 Then the natural maps are injections:
	$$\Pic_{X/S}\into \Pic_{(X/S)\zar}\into \Pic_{(X/S)\et}\into
	\Pic_{(X/S)\fppf}.$$

\tu{2.}\enspace All three maps are isomorphisms if also $f$ has a
section; the latter two maps are isomorphisms if also $f$ has a section
locally in the Zariski topology; and the last map is isomorphism if also
$f$ has a section locally in the \'etale topology.
 \end{thm}

\begin{ex}\label{ex:gpts}
 Assume $\mc O_{\!S}\risom f_*\mc O_{\!X}$ holds universally.  Using
Theorem~\ref{th:cmp}, show that its four functors  have the same
geometric points; in other words, for every algebraically closed field
$k$ containing the residue class field of a point of $S$, the $k$-points
of all four functors are, in a natural way, the same.  Show, in fact,
that these $k$-points are just the elements of $\Pic(X_k)$.

What if $\mc O_{\!S}\risom f_*\mc O_{\!X}$ does not necessarily hold
universally?
 \end{ex}

\begin{lem}\label{lm:fff}
 Assume $\mc O_{\!S}\risom f_*\mc O_{\!X}$.  Then the functor $\mc N\mapsto
f^*\mc N$ is fully-faithful from the category $\mc C$ of locally free
sheaves of finite rank on $S$ to that on $X$.  The essential image is
formed by the sheaves $\mc M$ on $X$ such that \tu{(i)} the image
$f_*\mc M$ is in $\mc C$ and \tu{(ii)} the natural map $f^*f_*\mc
M\to\mc M$ is an isomorphism.
 \end{lem}

\begin{proof} (Compare with \cite[21.13.2]{EGAIV4}.) For any $\mc
N$ in $\mc C$, there is a string of three natural isomorphisms
  \begin{equation}\label{eq:2a}
	\mc N\risom \mc N\ox f_*\mc O_{\!X}\risom
	  \mc N\ox f_*f^*\mc O_{\!S}\risom f_*f^*\mc N.
  \end{equation}
 The first isomorphism arises by tensor product with the comorphism of
$f$; this comorphism is an isomorphism by hypothesis.  The second
isomorphism arises from the identification $\mc O_{\!X}=f^*\mc O_{\!S}$.  The
third arises from the projection formula.

For any $\mc N'$ in $\mc C$, also $\SHom(\mc N,\,\mc N')$ is in $\mc C$.
Hence,
 (\ref{eq:2a}) yields an isomorphism
	$$\SHom(\mc N,\,\mc N')\risom f_*f^*\SHom(\mc N,\,\mc N').$$
 Now, since $\mc N$ and $\mc N'$ are locally free of finite rank, the
natural map
	$$f^*\SHom(\mc N,\,\mc N')\to \SHom(f^*\mc N,\,f^*\mc N')$$
 is an isomorphism locally, so globally.  Hence there is an isomorphism
of groups
	$$\Hom(\mc N,\,\mc N')\risom \Hom(f^*\mc N,\,f^*\mc N').$$
 In other words, $\mc N\mapsto f^*\mc N$ is fully-faithful.

Finally, the essential image consists of the sheaves $\mc M$ that are
isomorphic to those of the form $f^*\mc N$ for some $\mc N$ in $\mc C$.
Given such an $\mc M$ and $\mc N$, there is an isomorphism $f_*\mc
M\simeq \mc N$ owing to (\ref{eq:2a}).  Hence $f_*\mc M$ is in $\mc C$,
and $f^*f_*\mc M\to\mc M$ is an isomorphism locally, so globally; thus
(i) and (ii) hold.  Conversely, if (i) and (ii) hold, then $\mc M$ is,
by definition, in the essential image.
 \end{proof}

\begin{proof}[Proof of Part 1 of Theorem \ref{th:cmp}]
 Given $\lambda\in\Pic_{X/S}(T)$, represent $\lambda$ by an invertible
sheaf $\mc L$ on $X_T$.  Suppose $\lambda$ maps to $0$ in
$\Pic_{(X/S)\fppf}(T)$.  Then there exist an fppf covering $p\:T'\to T$
and an isomorphism $p_X^*\mc L\simeq f_{T'}^*\mc N'$ for some invertible
sheaf $\mc N'$ on $T'$. Hence Lemma~\ref{lm:fff} implies that
$f_{\smash{T'}*}p_X^*\mc L\simeq\mc N'$.  Now,  $p$ is flat, so
$p^*f_{\smash{T}*}\mc L\risom f_{\smash{T'}*}p_X^*\mc L$. So
$p^*f_{\smash{T}*}\mc L\simeq\mc N'$.  Hence $f_{\smash{T}*}\mc L$
is invertible and the natural map $f^*_Tf_{\smash{T}*}\mc L\to\mc L$ is
an isomorphism, as both statements hold after pullback via $p$, which
is faithfully flat.  Therefore, $\lambda=0$.  Thus $\Pic_{X/S}\into
\Pic_{(X/S)\fppf}$.

The rest is formal.  Indeed, take the last injection, and form the
associated sheaves in the Zariski topology. This operation is exact by
general (Grothendieck) topology, and $\Pic_{(X/S)\fppf}$ remains the
same, as it is already a Zariski sheaf.  Thus $\Pic_{(X/S)\zar}\into
\Pic_{(X/S)\fppf}$.  Similarly, $\Pic_{(X/S)\et}\into
\Pic_{(X/S)\fppf}$.

Alternatively, we can avoid the use of Lemma~\ref{lm:fff} by starting
from the fact that $\Pic_{(X/S)\fppf}$ is the sheaf associated to
$\Pic_X$, rather than to $\Pic_{X/S}$.  This way, we may assume $\mc
N'=\mc O_{T'}$.  Then $f_{\smash{T'}*}p_X^*\mc L\simeq\mc N'$ because
$\mc O_{\!S}\risom f_*\mc O_{\!X}$ holds universally. We now proceed
just as before.
 \end{proof}

\begin{dfn}\label{df:rgd}
 Assume $f$ has a section $g$, so $fg=1$.  Let $T$ be an $S$-scheme,
and $\mc L$ a sheaf on $X_T$.  Then a $g$-{\it rigidification} of $\mc
L$ is the choice of an isomorphism $u\:\mc O_{\!T}\risom g^*_T\mc L$,
assuming one exists.
 \end{dfn}

\begin{lem}\label{lm:idn}
 Assume $f$ has a section $g$, and let $T$ be an $S$-scheme.  Form the
group of isomorphism classes of pairs  $(\mc L,u)$ where $\mc L$ is an
invertible sheaf on $X_T$ and $u$ is a $g$-rigidification of $\mc
L$.  Then this group is carried isomorphically onto $\Pic_{X/S}(T)$ by
the homomorphism $\rho$ defined by $\rho(\mc L,u):=\mc L$. 
 \end{lem}
 \begin{proof}
 Given  $\lambda$ in $\Pic_{X/S}(T)$, represent  $\lambda$ by an
invertible sheaf $\mc M$ on $X_T$.  Set $\mc L:=\mc M\ox(f_T^*g^*_T\mc
M)^{-1}$. Then $\mc L$ too represents $\lambda$.  Also
$g^*_T\mc L =g_T^*\mc M\ox g_T^*f_T^*g_T^*\mc M^{-1}$.  Now,
$g_T^*f_T^*=1$ as $fg=1$.  Hence the natural isomorphism $g_T^*\mc
M\ox(g_T^*\mc M)^{-1}\risom \mc \mc O_{\!T}$ induces a $g$-rigidification of
$\mc L$.  Thus $\rho$ is surjective.

To prove $\rho$ is injective, let $(\mc L,u)$ represent an element of
its kernel. Then there exist an invertible sheaf $\mc N$ on $T$ and an
isomorphism $v\:\mc L\risom f^*_T\mc N$.  Set $w:=g_T^*v\circ u$, so
$w\:\mc O_{\!T}\risom g^*_T\mc L\risom\mc N$.  Now, a map of pairs is just a
map $w'$ of the first components such that $g^*_Tw'$ is compatible with
the two $g$-rigidifications.  So $v\:(\mc L,u)\risom(f^*_T\mc N,w)$ and
$f^*_Tw\:(\mc O_{\!X_T},1)\risom(f^*_T\mc N,w)$. Thus $\rho$ is injective.
\end{proof}

\begin{lem}\label{lm:aut}
 Assume $f$ has a section $g$, and assume $\mc O_{\!S}\risom f_*\mc O_{\!X}$
holds universally.  Let $T$ be an $S$-scheme, $\mc L$ an invertible
sheaf on $X_T$, and $u$ a $g$-rigidification of $\mc L$.  Then every
automorphism of the pair $(\mc L,u)$ is trivial.
 \end{lem}
 \begin{proof}
 An automorphism of $(\mc L,u)$ is just an automorphism $v\:\mc L\risom
\mc L$ such that $g_T^*v\circ u\:\mc O_{\!T}\risom g^*_T\mc L\risom g^*_T\mc
L$ is equal to $u$.  But then $g_T^*v=1$.  Now,
	$$v\in\Hom(\mc L,\mc L)=\uH^0(\SHom(\mc L,\mc L))
	=\uH^0(\mc O_{\!X_T})=\uH^0(\mc O_{T});$$
 the middle equation holds since the natural map $O_{\!X_T}\to\SHom(\mc
L,\mc L)$  is locally an isomorphism, so globally one, and the
last equation holds since $\mc O_{\!T}\risom f_*\mc O_{\!X_T}$.  But
$g_T^*v=1$.  Therefore,  $v=1$.
 \end{proof}

\begin{proof}[Proof of Part 2 of Theorem \ref{th:cmp}]
 Suppose $f$ has a section $g$.  Owing to Part 1, it suffices to prove
that every $\lambda\in \Pic_{(X/S)\fppf}(T)$ lies in $\Pic_{X/S}(T)$.
Represent $\lambda$ by a $\lambda'\in \Pic_{X/S}(T')$ where $T'\to
T$ is an fppf covering.  Then there is an fppf covering $T''\to
T'\x_TT'$ such that the two pullbacks of $\lambda'$ to $X_{T''}$ are
equal.  We may assume $T''\risom T'\x_T T'$ because $\Pic_{X/S}$ is separated
for the fppf topology, again owing to Part 1.

Owing to Lemma~\ref {lm:idn}, we may represent $\lambda'$ by a pair
$(\mc L',u')$ where $\mc L'$ is an invertible sheaf on $X_{T'}$ and $u'$
is a $g$-rigidification of $\mc L'$.  Furthermore, on $X_{T'\x T'}$, there
is an isomorphism $v'$ from the pullback of $(\mc L',u')$ via the first
projection onto the pullback via the second.

Consider the three projections $X_{T'\x T'\x T'} \to X_{T'\x T'}$.  Let
$v'_{ij}$ denote the pullback of $v'$ via the projection to the $i$th
and $j$th factors.  Then $v_{13}^{\prime-1}v'_{23}v'_{12}$ is an
automorphism of the pullback of $(\mc L',u')$ via the first projection
$X_{T'\x T'\x T'} \to X_{T'}$.  So, owing to Lemma~\ref {lm:aut}, this
automorphism is trivial.  Hence $(\mc L',u')$ descends to a pair
$(\mc L,u)$ on $X_T$.  Therefore, $\lambda$ lies in $\Pic_{X/S}(T)$.

The rest is formal.  Indeed, suppose that there is a Zariski covering
$T'\to T$ such that $f_{T'}$ has a section.  Then, by the above,
$\Pic_{(X/S)\zar}|T'\risom \Pic_{(X/S)\fppf}|T'$.  Hence, by general
(Grothendieck) topology, $\Pic_{(X/S)\zar}\risom \Pic_{(X/S)\fppf}$
since both source and target are sheaves in the Zariski topology and a map
of sheaves is an isomorphism if it is so locally.  Similarly,
$\Pic_{(X/S)\et}\risom \Pic_{(X/S)\fppf}$ if $f$ has a section locally
in the \'etale topology.
 \end{proof}

\begin{rmk}\label{rk:coh}
 There is another way to prove Theorem~\ref{th:cmp}.  This way is more
sophisticated, and  yields more information, which we won't need.
Here is the idea.

Recall \cite[Ex.~III, 4.5, p.~224]{Ha83} that, for any ringed space $R$,
there is a natural isomorphism
 \begin{equation}\label{eq:2b}
	\Pic(R) = \uH^1(R,\mc O_{\!R}^*).
  \end{equation}
 Now, given any $S$-scheme  $T$, form the presheaf $T'\mapsto
\uH^1(X_{T'},\mc O_{\!X_{T'}}^*)$ on $T$.  Its associated sheaf is
\cite[Prp.~III, 8.1, p.~250]{Ha83} simply $\tu R^1f_{\smash{T}*}\mc
O_{\!X_T}^*$.  Therefore,
  \begin{equation}\label{eq:2e}
   \Pic_{(X/S)\zar}(T)=\uH^0(T,\tu R^1f_{\smash{T}*}\mc O_{\!X_T}^{*}).
  \end{equation}

Consider the Leray spectral sequence \cite[Thm.~II, 4.17.1, p.~201]{Gd58}
    $$\tu E_2^{pq}:=\uH^p(T,\tu R^qf_{\smash{T}*}\mc O_{\!X_T}^*)
    \Longrightarrow \uH^{p+q}(X_T,\mc O_{\!X_T}^*),$$
 and form its exact sequence of terms of low degree \cite[Thm.~I, 4.5.1,
p.~82]{Gd58}
\begin{multline}\label{eq:2c}
 0\to\uH^1(T,f_{\smash{T}*}\mc O_{\!X_T}^*)\to\uH^1(X_T,\mc O_{\!X_T}^*)
 \to\uH^0(T,\tu R^1f_{\smash{T}*}\mc O_{\!X_T}^*)\\
 \to\uH^2(T,f_{\smash{T}*}\mc O_{\!X_T}^*)\to\uH^2(X_T,\mc O_{\!X_T}^*).
\end{multline}

If $\mc O_{\!S}\risom f_*\mc O_{\!X}$ holds universally, then $\uH^1(T,\mc
O_{\!T}^*)\risom\uH^1(T,f_{\smash{T}*}\mc O_{\!X_T}^*)$.  Hence the beginning
of (\ref{eq:2c}) becomes
	$$0\to\Pic(T)\to\Pic(X_T)\to\Pic_{(X/S)\zar}(T).$$
 Thus $\Pic_{X/S}\into \Pic_{(X/S)\zar}$.

If also $f$ has a section $g$, then $g$ induces, for each $p$, a left
inverse of the map $\uH^p(T,f_{\smash{T}*}\mc O_{\!X_T}^*)\to\uH^p(X_T,\mc
O_{\!X_T}^*)$ induced by $f$.  So the latter is injective.  Hence,
$\uH^1(X_T,\mc O_{\!X_T}^*)\to\uH^0(T,\tu R^1f_{\smash{T}*}\mc O_{\!X_T}^*)$
is surjective.  Thus $\Pic_{X/S}\risom \Pic_{(X/S)\zar}$.

The preceding argument works for the \'etale and fppf topologies
too with little change.  First of all, consider the functor
	$$\Gm(T):=\uH^0(T,\mc O_{\!T}^*).$$
 Let $u$ be an indeterminate.  Then $\Gm(T)$ is representable by the
$S$-scheme
	$$\Gm:=\Spec(\mc O_{\!S}[u,u^{-1}]).$$
 Indeed, giving an $S$-map $T\to\Gm$ is the same as giving an
$\uH^0(S,\mc O_{\!S})$-homo\-mor\-phism from $\uH^0(S,\mc
O_{\!S})[u,u^{-1}]$ to $\uH^0(T,\mc O_{\!T})$, and giving such a
homomorphism is the same as assigning to $u$ a unit in $\uH^0(T,\mc
O_{\!T})$.  Now, since $\Gm(T)$ is representable, it is a sheaf for the
\'etale and fppf topologies.

Grothendieck's generalization \cite[p.~190-16]{FGA} of Hilbert's Theorem
90 asserts the formula
	$$\Pic(T)=\uH^1(T,\Gm)$$
where the $\uH^1$ can be computed in either the \'etale or fppf
topology.  The proof is simple, and similar to the proof of
(\ref{eq:2b}).  The $\uH^1$ can be computed as a \v Cech group.  And,
essentially by definition, for a covering $T'\to T$, a \v Cech cocycle 
with values in $\Gm$ amounts to descent data on $\mc
O_{T'}$.  The data is effective by descent theory, and the resulting
sheaf on $T$ is invertible since the covering is faithfully flat.

In the present context, the exact sequence~(\ref{eq:2c}) becomes
\begin{multline}\label{eq:2d}
 0\to\uH^1(T,f_{\smash{T}*}\Gm)\to\uH^1(X_T,\Gm)
 \to\uH^0(T,\tu R^1f_{\smash{T}*}\Gm)\\
 \to\uH^2(T,f_{\smash{T}*}\Gm)\to\uH^2(X_T,\Gm).$$
\end{multline}
 Furthermore, the proof of (\ref{eq:2e}) yields, for example in the fppf
topology,
	$$\Pic_{(X/S)\fppf}(T)=\uH^0(T,\tu R^1f_{\smash{T}*}\Gm).$$
 
If $\mc O_{\!S}\risom f_*\mc O_{\!X}$ holds universally, then it follows from
the definitions that $f_{\smash{T}*}\Gm=\Gm$. Hence the beginning
of (\ref{eq:2d}) becomes
	$$0\to\Pic(T)\to\Pic(X_T)\to\Pic_{(X/S)\fppf}(T).$$
 Thus $\Pic_{X/S}\into \Pic_{(X/S)\fppf}$.  And $\Pic_{X/S}\risom
\Pic_{(X/S)\fppf}$ if also $f$ has a section.  As before, the rest of
Theorem~\ref{th:cmp} follows formally.

The \'etale group $\uH^2(T,\Gm)$ was studied extensively by Grothendieck
\cite[pp.~46--188]{Dix}.  He showed that it gives one of two significant
generalizations of the Brauer group of central simple algebras over a
field.  The other is the group of Azumaya algebras on $T$.  He denoted
the latter by $\Br(T)$ and the former by $\Br'(T)$. Hence
\cite[pp.~127--128]{Dix}, if $\mc O_{\!S}\risom f_*\mc O_{\!X}$ holds
universally, then (\ref{eq:2d}) becomes
\begin{equation*}\label{eq:2ee}
 0\to\Pic(T)\to\Pic(X_T)\to\Pic_{(X/S)\et}(T)\to\Br'(T)\to\Br'(X_T);
 \end{equation*}
 in particular, the obstruction to representing an element of
$\Pic_{(X/S)\et}(T)$ by an invertible sheaf on $X_T$ is given by an
element of $\Br'(T)$, which maps to $0$ in $\Br'(X_T)$.

Using the smoothness of $\Gm$ as an $S$-scheme, Grothendieck
\cite[p.~180]{Dix} proved  the natural homomorphisms are
isomorphisms from the \'etale groups $\uH^p(T,\Gm)$ to the
corresponding fppf groups.  Hence, if $\mc O_{\!S}\risom f_*\mc O_{\!X}$
holds universally, then it follows from (\ref{eq:2d}) via the Five Lemma
that,  whether $f$ has a section locally in the \'etale topology or not,
 \begin{equation*}
   \Pic_{(X/S)\et}\risom\Pic_{(X/S)\fppf}.
 \end{equation*}



Nevertheless, when a discussion is set in the greatest possible
generality, it is common to work with $\Pic_{(X/S)\fppf}$ and call it
{\it the\/} Picard functor.
  \end{rmk}

\section{Relative effective divisors}\label{sc:red}

Grothendieck constructed the Picard scheme by taking a suitable family
of effective divisors and forming the quotient modulo linear
equivalence.  This section develops the basic theory of these notions.

\begin{sbs}[Effective divisors] \label{sb:ediv}
 A closed subscheme $D\subset X$ is called an {\it effective (Cartier)
divisor\/} if its ideal $\mc I$ is invertible.  Given an $\mc
O_{\!X}$-module $\mc F$ and   $n\in\bb Z$, set
	$$\mc F(nD):=\mc F\ox\mc I^{\ox -n}.$$

In particular, $\mc O_{\!X}(-D)=\mc I$.  So the inclusion $\mc I\into\mc
O_{\!X}$ yields, via tensor product with $\mc O_{\!X}(D)$, an injection
$\mc O_{\!X}\into \mc O_{\!X}(D)$, which, in turn, corresponds to a
global section of $\mc O_{\!X}(D)$.  This section is not arbitrary since
it corresponds to an injection.  Sections corresponding to injections
are termed {\it regular}.

Conversely, given an arbitrary invertible sheaf $\mc L$ on $X$, let
$\uH^0(X,\mc L)_{\reg}$ denote the subset of $\uH^0(X,\mc L)$ consisting
of the regular sections, those corresponding to injections $\mc
L^{-1}\into\mc O_{\!X}$.  And let $|\mc L|$ denote the set of effective
divisors $D$ such that $\mc O_{\!X}(D)$ is, in some way, isomorphic to $\mc
L$.  For historical reasons, $|\mc L|$ is called the {\it complete
linear system\/} associated to $\mc L$ (but $|\mc L|$ needn't be a $\I P^n$ 
if $X$ isn't integral).
 \end{sbs}

\begin{ex}\label{ex:lsys}
 Under the conditions of (\ref{sb:ediv}), establish a canonical isomorphism
     $$\uH^0(X,\mc L)_{\reg}\big/\uH^0(X,\mc O_{\!X}^*)\risom|\mc L|.$$
  \end{ex}

\begin{dfn}\label{df:red}
 A {\it relative effective divisor\/} on $X/S$ is an effective divisor
$D\subset X$ that is  $S$-flat.
 \end{dfn}

\begin{lem}\label{lm:ctn}
 Let $D\subset X$ be a closed subscheme, $x\in D$ a point, and $s\in S$
its image.  Then the following statements are equivalent:
 \begin{enumerate}
  \item The subscheme $D\subset X$ is a relative effective divisor at
$x$ (that is, in a neighborhood of $x$).
  \item The schemes $X$ and $D$ are $S$-flat at $x$, and the fiber $D_s$
is an effective divisor on $X_s$ at $x$. 
  \item The scheme $X$ is $S$-flat at $x$, and the subscheme $D\subset
X$ is cut out at $x$ by one element that is regular (a nonzerodivisor)
on the fiber $X_s$.
 \end{enumerate}
\end{lem}

\begin{proof} 
 For convenience, set $A:=\mc O_{S,s}$ and denote its residue field by
$k$.  In addition, set $B:=\mc O_{\!X,x}$ and $C:=\mc O_{D,x}$.  Then
$B\ox_A k=\mc O_{\!X_s,x}$.

Assume (i), and let's prove (ii).  By hypothesis, $D$ is an effective
divisor at $x$.  So there is a regular element $b\in B$ that generates
the ideal of $D$.  Multiplication by $b$ defines a short exact sequence
	$$0\to B\to B\to C\to0.$$
 In turn, this sequence induces the following exact sequence:
 $$\Tor^A_1(B,k)\to\Tor^A_1(B,k)\to\Tor^A_1(C,k)\to B\ox k\to B\ox k.$$

By hypothesis, $D$ is $S$-flat at $x$; hence, $\Tor^A_1(C,k)=0$.  So
$B\ox k\to B\ox k$ is injective.  Its image is the ideal of $D_s$.  Thus
$D_s$ is an effective divisor.

Since $\Tor^A_1(C,k)=0$, the map $\Tor^A_1(B,k)\to\Tor^A_1(B,k)$ is
surjective.  This map is given by multiplication by $b$, and $b$ lies in
the maximal ideal of $B$.  Also, $\Tor^A_1(B,k)$ is a finitely
generated $B$-module.  Hence,  $\Tor^A_1(B,k)=0$ by Nakayama's lemma.
Therefore, by the local criterion of flatness \cite[Thm.~5.6,
p.~98]{SGA1} or \cite[Thm.~6.1, p.~73]{OB72}, $B$ is $A$-flat; in other
words, $X$ is $S$-flat at $x$.  Thus (ii) holds.

Assume (ii).  To prove (iii), denote the ideal of $D$ in $B$ by $I$, and
that of $D_s$ in $B\ox k$ by $I'$.  Take an element $b\in I$ whose image $b'$ 
in $B\ox k$ generates $I'$; such a $b$ exists because $D_s$ is an
effective divisor at $x$ by hypothesis.  For the same reason, $b'$ is
regular.  It remains to prove $b$ generates $I$. 

Consider the short exact sequence
	$$0\to I\to B\to C\to0.$$
 By hypothesis, $C$ is $A$-flat.  Hence the map $I\ox k\to B\ox k$ is
injective.  Its image is $I'$, which is generated by $b'$.  So the image
of $b$ in $I\ox k$ generates it.  Hence, by Nakayama's lemma, $b$
generates $I$.  Thus (iii) holds.
 
Assume (iii).  To prove (i), again denote the ideal of $D$ in $B$ by
$I$.  By hypothesis, $I$ is generated by an element $b$ whose image $b'$
in $B\ox k$ is regular.  We have to prove $b$ is regular and $C$ is
$A$-flat.

The exact sequence $0\to I\to B\to C\to0$ yields this one:
 \begin{equation}\label{eq:3a}
  \Tor^A_1(B,k)\to\Tor^A_1(C,k)\to I\ox k\to B\ox k.
  \end{equation}
 The last map is injective for the following reason.  Since $I=Bb$,
multiplication by $b$ induces a surjection $B\to I$, so a surjection
$B\ox k\to I\ox k$.  Consider the composition
	$$B\ox k\to I\ox k\to B\ox k.$$
 It is given by multiplication by $b'$, so is injective because $b'$ is
regular. Hence $B\ox k\risom I\ox k$.  Therefore, $I\ox k\to B\ox k$ is
injective.

By hypothesis, $B$ is $A$-flat.  So $\Tor^A_1(B,k)=0.$ Hence the
exactness of (\ref{eq:3a}) implies $\Tor^A_1(C,k)=0$.  Therefore, by the
local criterion, $C$ is $A$-flat.

Since $B$ and $C$ are $A$-flat and $0\to I\to B\to C\to0$ is exact, also
$I$ is $A$-flat.

Define $K$ by the exact sequence  $0\to K\to B\to I\to0$.  Then the
sequence
	$$0\to K\ox k\to B\ox k\to I\ox k\to0$$
 is exact since $I$ is $A$-flat.  But $B\ox k\risom I\ox k$. Hence $K\ox
k=0$.  Therefore, $K=0$ by Nakayama's lemma.  But $K$ is the kernel of
multiplication by $b$ on $B$.  So $b$ is regular.  Thus (i) holds.
 \end{proof}

\begin{ex}\label{ex:sum}
 Let $D$ and $E$ be relative effective divisors on $X/S$, and $D+E$
their sum. Show $D+E$ is a relative effective divisor too.
 \end{ex}

\begin{dfn}\label{df:div}
 Define a functor $\Div_{X/S}$ by the formula
    $$\Div_{X/S}(T):=\{\,\text{relative effective divisors $D$
     on $X_T/T$}\,\}.$$
 \end{dfn}

Note $\Div_{X/S}$ is indeed a functor.  Namely, given a relative
effective divisor $D$ on $X_T/T$ and an arbitrary $S$-map $p\:T'\to T$,
we have to see the $T'$-flat closed subscheme $D_{T'}\subset
X_{T'}$ is an effective divisor.  So let $\mc I$ denote the ideal of
$D$.  Since $D$ is $T$-flat, $p_{X_T}^*\mc I$ is equal to the ideal of
$D_{T'}$.  But, since $\mc I$ is invertible, so is $p_{X_T}^*\mc I$.
Thus $D_{T'}$ is a (relative) effective divisor.

\begin{thm}\label{th:repDiv}
 Assume $X/S$ is projective and flat.  Then $\Div_{X/S}$ is
representable by an open subscheme $\IDiv_{X/S}$ of the Hilbert scheme
$\IHilb_{X/S}$.
 \end{thm}

\begin{proof}
 Set $H:=\IHilb_{X/S}$, and let $W\subset X\x H$ be the universal
(closed) subscheme, and $q\:W\to H$ the projection.  Let $V$ denote the
set of points $w\in W$ at which $W$ is an effective divisor. Plainly $V$
is open in $W$.  Set $Z:=q(W-V)$.  Then $Z$ is closed because $q$ is
proper.  Set $U:=H-Z$.  Then $U$ is open, and $q^{-1}U$ is an effective
divisor in $X\x U$.  In fact, since $q$ is flat, $q^{-1}U$ is a relative
effective divisor in $X\x U/U$.

It remains to show that $U$ represents $\Div_{X/S}$.  So let $T$ be an
$S$-scheme, and $D\subset X_T/T$ a relative effective divisor.  By the
universal property of the pair $(H,W)$, there exists a unique map
$g\:T\to H$ such that $g_X^{-1}W=D$.  We have to show that $g$ factors
through $U$.

For each $t\in T$, the fiber $D_t$ is an effective divisor since it is
obtained by base change (or owing to Lemma~\ref{lm:ctn}).  But
$D_t=W_{g(t)}\ox k_t$ where $k_t$ is the residue field of $t$.  So
$W_{g(t)}$ too is a divisor, as a field extension is faithfully flat.
Hence, since $X\x H$ and $W$ are $H$-flat, $W$ is, by
Lemma~\ref{lm:ctn}, a relative effective divisor along the fiber over
$g(t)$.  Therefore, $g(t)\in U$.  So, since $U$ is open, $g$ factors
through $U$.
 \end{proof}

\begin{ex}\label{ex:DivC}
  Assume $f\:X\to S$ is flat and is projective Zariski locally over $S$.
Assume its fibers are curves, that is, of {\it pure} dimension 1.  Given
$m\ge1$, let $\Div_{X/S}^m$ be the functor whose $T$-points are the
relative effective divisors $D$ on $X_T/T$ with fibers $D_t$ of degree
$m$.

Show the $\Div_{X/S}^m$ are representable by open and closed subschemes
of finite type $\IDiv_{X/S}^m\subset\IDiv_{X/S}$, which are disjoint and
cover.

Let $X_0\subset X$ be the subscheme where $X/S$ is smooth.  Show
$X_0=\IDiv_{X/S}^1$.

Let  $X_0^m$ be the $m$-fold  $S$-fibered product.  Show there is a natural
$S$-map $$\alpha\:X_0^m\to\IDiv_{X/S}^m,$$ which is given on  $T$-points
by $\alpha(\Gamma_1,\dotsc,\Gamma_m)= \sum \Gamma_i$.
 \end{ex}

\begin{rmk}\label{rmk:symprod}
 Consider the map $\alpha$ of Exercise~\ref{ex:DivC}.  Plainly $\alpha$
is compatible with permuting the factors of $X_0^m$.  Hence $\alpha$
factors through the symmetric product $X_0^{(m)}$.  In fact, $\alpha$
induces an isomorphism
	$$X_0^{(m)}\risom \IDiv_{X_0/S}^m.$$
 This isomorphism is treated in detail in \cite[Prp.~6.3.9,
p.~437]{De73} and in outline in \cite[Prp.~3, p.~254]{BLR}.
 \end{rmk}

\begin{sbs}[The module $\mc Q$]\label{sb:Q}
 Assume $f\:X\to S$ is proper, and let $\mc F$ be a coherent $\mc
O_{\!X}$-module flat over $S$.  Recall from \cite[7.7.6]{EGAIII2} that there
exist a coherent $\mc O_{\!S}$-module $\mc Q$ and an isomorphism of functors
in the quasi-coherent $\mc O_{\!S}$-module $\mc N$:
 \begin{equation}\label{eq:Q}
	q\:\SHom(\mc Q, \mc N)\risom f_*(\mc F\ox f^*\mc N).
 \end{equation}
  The pair $(\mc Q, q)$ is unique, up to unique isomorphism, and by
\cite[7.7.9]{EGAIII2}, forming it commutes with changing the base, in
particular, with localizing.

Fix $s\in S$ and assume $S=\Spec(\mc O_{S,s})$.  Note that the following
conditions are equivalent:
 \begin{enumerate}
  \item the $\mc O_{\!S}$-module $\mc Q$ is free (or equivalently,
projective);
  \item the functor $\mc N\mapsto f_*(\mc F\ox f^*\mc N)$ is right
exact;
  \item for all $\mc N$, the natural map is an isomorphism, $f_*(\mc
F)\ox \mc N\risom f_*(\mc F\ox f^*\mc N)$;
  \item the natural map is a surjection,
 $\uH^0(X,\mc F)\ox k_s\onto\uH^0(X_s,\mc F_s)$.
 \end{enumerate}
 Indeed, the equivalence of (i)--(iii) is elementary and straightforward.
Moreover, (iv) is a special case of (iii).  Conversely, (iv) implies
(iii) by \cite[7.7.10]{EGAIII2} or \cite[Cor.~5.1 p.~72]{OB72}; this
useful implication is known as the ``property of exchange.''

In addition, (i)--(iv) are implied by the following condition:
 \begin{enumerate}\setcounter{enumi}{4}
  \item the first cohomology group of the fiber vanishes, $\uH^1(X_s,\mc
F_s)=0$.
 \end{enumerate} 
 Indeed, (v) implies that $\tu R^1f_*(\mc F\ox f^*\mc N)=0$ for all $\mc
N$ by \cite[7.5.3]{EGAIII2} or \cite[Cor.~2.1 p.~68]{OB72}; in turn, this
vanishing implies (ii) owing to the long exact sequence of higher direct
images.
 \end{sbs}

\begin{ex}\label{ex:gc&r}
  Assume $f\:X\to S$ is proper and flat, and its geometric fibers are
reduced and connected.  Show $\mc O_{\!S}\risom f_*\mc O_{\!X}$ holds
universally.
 \end{ex}

\begin{dfn}\label{df:LinSys}
Let $\mc L$ be an invertible sheaf on $X$.  Define a subfunctor
$\LinSys_{\mc L/X/S}$ of $\Div_{X/S}$ by the formula
  \begin{multline*}
   \LinSys_{\mc L/X/S}(T):=\{\,\text{relative effective divisors $D$
     on $X_T/T$ such that}\\
     \text{$\mc O_{\!X_T}(D)\simeq \mc L_T\ox f_T^*\mc N$
		for some invertible sheaf $\mc N$ on }T\,\}.
  \end{multline*}
 \end{dfn}

Notice the similarity of this definition with that,
Definition~\ref{df:Pfs}, of $\Pic_{X/S}$: both definitions work with
isomorphism classes of invertible sheaves on $X_T$ modulo $\Pic(T)$, in
the hope of producing a representable functor.  Here, this hope is
fulfilled under suitable hypotheses on $f$, according to the next
theorem.

\begin{thm}\label{th:LinSys} Assume $X/S$ is proper and flat, and its
geometric fibers are integral.  Let $\mc L$ be an invertible sheaf on
$X$, and let $\mc Q$ be the $\mc O_{\!S}$-module associated in
Subsection~\ref{sb:Q} to $\mc F:=\mc L$. Set $L:=\I P(\mc Q)$.  Then $L$
represents $\LinSys_{\mc L/X/S}$.
 \end{thm}

\begin{proof}
 Let $D\in\LinSys_{\mc L/X/S}(T)$.  Say $\mc O_{\!X_T}(D)\simeq \mc
L_T\ox f_T^*\mc N$.  Then $\mc N$ is determined up to isomorphism.
Indeed, let $\mc N'$ be a second choice.  Then $$\mc L_T\ox f_T^*\mc
N\simeq\mc L_T\ox f_T^*\mc N'.$$  So $f_T^*\mc N\simeq f_T^*\mc N'$
since $\mc L$ is invertible.  Hence $\mc N\simeq\mc N'$ by
Lemma~\ref{lm:fff}, which applies to $f_T\:X_T\to T$ since $\mc
O_{\!S}\risom f_*\mc O_{\!X}$ holds universally by Exercise~\ref{ex:gc&r}.

Say $D$ is defined by $\sigma\in\uH^0(X_T,\mc L_T\ox f_T^*\mc N)$.
Now, forming $\mc Q$ commutes with changing the base to $T$, and so
(\ref{eq:Q}) becomes
 \begin{equation}\label{eq:QT}
  \SHom(\mc Q_T, \mc N)\risom f_{\smash{T}*}(\mc L_T\ox f_T^*\mc N).
 \end{equation}
 Hence $\sigma$ corresponds to a map $u\:\mc Q_T\to\mc N$.
 
Let $t\in T$.  Since $D$ is a relative effective divisor on $X_T/T$, its
fiber $D_t$ is a divisor on $X_t$ by Lemma~\ref{lm:ctn}.  Since $D_t$ is
defined by $\sigma_t\in \uH^0(X_t,\mc L|X_t)$, necessarily $\sigma_t
\neq 0$.  But $\sigma_t$ corresponds to $u\ox k_t$, so $u\ox k_t\neq 0$.
Now, $\mc N$ is invertible, so $\mc N\ox k_t$ is a $k_t$-vector space of
dimension 1.  So $u\ox k_t$ is surjective.  Hence, by Nakayama's lemma,
$u$ is surjective at $t$.  But $t$ is arbitrary.  So $u$ is surjective
everywhere.

Therefore, $u\:\mc Q_T\to\mc N$ defines an $S$-map $p\:T\to L$ by
\cite[4.2.3]{EGAII}.  Since $(\mc N, u)$ is determined up to
isomorphism, a second choice yields the same $p$.

Plainly, this construction is functorial in $T$.  Thus we obtain a map
of functors,
\begin{equation*}
 	\Lambda\:\LinSys_{\mc L/X/S}(T)\to L(T).
 \end{equation*}
 Let us prove  $\Lambda$  is an isomorphism.

Let $p\in L(T)$, so  $p\:T\to L$ is an $S$-map.  Then $p$  arises from a
surjection $u\:\mc Q_T\onto\mc N$; namely, $u=p^*\alpha$ where
$\alpha\:\mc Q_L \onto\mc O(1)$ is the tautological map.  Moreover,
there is only one such pair $(\mc N, u)$ up to isomorphism.

Via the isomorphism in (\ref{eq:QT}), the surjection $u$
corresponds to a global section $\sigma\in\uH^0(X_T,\mc L_T\ox f_T^*\mc
N)$. Let $t\in T$.  Then $u\ox k_t$ is surjective, so $u\ox k_t\neq 0$.
But $u\ox k_t$ corresponds to $\sigma_t\in \uH^0(X_t,\mc L|X_t)$, so
$\sigma_t \neq 0$.  But $X_t$ is integral since the geometric fibers of
$X/S$ are integral by hypothesis.  Hence $\sigma_t$ is regular.

The section $\sigma$ defines a map $(\mc L_T\ox f_T^*\mc N)^{-1} \to \mc
O_{\!X_T}$.  Its image is the ideal of a closed subscheme $D\subset X$,
which is cut out locally by one element; moreover, on the fiber $X_t$,
this element corresponds to $\sigma_t$, so is regular.  Hence $D$ is a
relative effective divisor on $X_T/T$ by Lemma~\ref{lm:ctn}.  In fact,
$D\in \LinSys_{\mc L/X/S}(T)$.  Plainly, $D$ is the only such divisor
corresponding to $(\mc N, u)$, so mapping to $p$ under $\Lambda$.

Thus $\Lambda$ is an isomorphism.  In other words, $L$ represents
$\LinSys_{\mc L/X/S}$.
 \end{proof}

\begin{ex}\label{ex:LinSys}
 Under the conditions of Theorem~\ref{th:LinSys}, show that there exists a
natural relative effective divisor $W$ on $X_L\big/L$ such that
 $$\mc O_{X_L}(W)=\mc L_L\ox f_L^*\mc O_L(1).$$
 Furthermore, $W$ possesses the following universal property: given any
$S$-scheme $T$ and any relative effective divisor $D$ on $X_T/T$ such
that $\mc O_{\!X_T}(D)\simeq\mc L_T\ox f_T^*\mc N$ for some
invertible sheaf $\mc N$ on $T$, there exist a unique $S$-map $w\:T\to
L$ such that $(1\x w)^{-1}W=D$.
 \end{ex}

\section{The Picard scheme}\label{sc:Psch}

This section proves Grothendieck's main theorem about the Picard scheme,
which asserts its existence if $X/S$ is projective and flat and its
geometric fibers are integral; in fact, the functor $\Pic_{(X/S)\et}$ is
representable.  The proof involves Grothendieck's method of using
functors to prescribe patching.  The basic theory is developed in
\cite[Ch.~0, Sct.~4.5, pp.~102--107]{EGAG}, and is applied in
\cite[Ch.~1, Sct.~9, pp.~354--401]{EGAG} to the construction of
Grassmannians and related parameter schemes.  The present construction
involves the basic theory, but is more sophisticated and more
complicated because it works not simply with the Zariski topology, but
also with the \'etale topology.

\begin{dfn}\label{df:Psch}
 If any of the four relative Picard functors of Definition~\ref{df:Pfs}
is representable, then the representing scheme is called the {\it Picard
scheme\/} and denoted by $\IPic_{X/S}$.  Moreover, we say simply that the
Picard scheme $\IPic_{X/S}$ exists.
 \end{dfn}

Notice that, although there are four relative Picard functors, there is
at most one Picard scheme.  Of course, if any functor is representable,
then the representing scheme is uniquely determined, up to a unique
isomorphism that preserves the identification of the given functor
with the functor of points of the representing scheme.  But here, there
is more to the story.

Indeed, for example, say $\Pic_{(X/S)\et}$ is representable by
$\IPic_{X/S}$.  Then, by descent theory, $\Pic_{(X/S)\et}$ is already a
sheaf in the fppf topology; so it is equal to its associated sheaf
$\Pic_{(X/S)\fppf}$.  Hence, it too is representable by $\IPic_{X/S}$.
On the other hand, $\Pic_{X/S}$ may or may not be representable;
however, if it is, then it, as well, must be representable by $\IPic_{X/S}$.

\begin{ex}\label{ex:0sec}
 Assume $\IPic_{X/S}$ exists, and $\mc O_{\!S}\risom f_*\mc O_{\!X}$
holds universally.  Let $T$ be an $S$-scheme, and $\mc L$ an invertible
sheaf on $X_T$.  Show that there exist a subscheme $N\subset T$ and an
invertible sheaf $\mc N$ on $N$ with these three properties: first, $\mc
L_N\simeq f_N^*{\mc N}$; second, given any $S$-map $t\:T'\to T$ such
that $\mc L_{T'}\simeq f_{T'}^*{\mc N'}$ for some invertible sheaf $\mc
N'$ on $T'$, necessarily $t$ factors through $N$ and $\mc N'\simeq
t^*\mc N$; and third, $N$ is a closed subscheme if $\IPic_{X/S}$ is
separated.  Show also that the first two properties determine $N$ uniquely
and $\mc N$ up to isomorphism.
 \end{ex}

\begin{ex}\label{ex:univshf}
 Assume $\IPic_{X/S}$ exists.  An invertible sheaf $\mc P$ on
$X\x\IPic_{X/S}$ is called a {\it universal sheaf}, or {\it Poincar\'e
sheaf}, if $\mc P$ possesses the following property: given any
$S$-scheme $T$ and any invertible sheaf $\mc L$ on $X_T$, there exists a
unique $S$-map $h\:T\to\IPic_{X/S}$ such that, for some invertible sheaf
$\mc N$ on $T$,
	$$\mc L\simeq (1\x h)^*\mc P\ox f_T^*\mc N.$$
 Show that a universal sheaf $\mc P$ exists if and only if $\IPic_{X/S}$
represents $\Pic_{X/S}$.

Assume $\mc O_{\!S}\risom f_*\mc O_{\!X}$ holds universally.  Show that, if
$\mc P$ exists, then it is unique up to tensor product with the pullback
of a unique invertible sheaf on $\IPic_{X/S}$.

Show that, if also $f$ has a section, then a universal sheaf $\mc P$
exists.

Find an example where no universal sheaf $\mc P$ exists.
 \end{ex}

\begin{ex}\label{ex:bschg}
 Assume $\IPic_{X/S}$ exists, and let $S'$ be an $S$-scheme.
Show that $\IPic_{X_{S'}/S'}$ exists too, and in fact, that
	$$\IPic_{X_{S'}/S'}=\IPic_{X/S}\x_SS'.$$
 Thus forming the Picard scheme commutes with changing the base.

Find an example where  $\IPic_{X_{S'}/S'}$ represents $\Pic_{X_{S'}/S'}$,
but  $\IPic_{X/S}$ does not represent  $\Pic_{X/S}$.
 \end{ex}

\begin{ex}\label{ex:schpts}
 Assume $\IPic_{X/S}$ exists, and either it represents
$\Pic_{(X/S)\fppf}$ or $\mc O_{\!S}\risom f_*\mc O_{\!X}$ holds
universally.  Show the scheme points of $\IPic_{X/S}$ correspond, in a
natural bijective fashion, to the classes of invertible sheaves $\mc L$
on the fibers of $X/S$.  A class is, by definition, represented by an
$\mc L$ on an $X_k$ where $k$ is a field containing the residue field
$k_s$ of a (scheme) point $s\in S$; an $\mc L'$ on an $X_{k'}$
represents the same class if and only if there is a third field $k''$
containing the other two such that $\mc L|X_{k''}\simeq\mc L'|X_{k''}$.
 \end{ex}

\begin{dfn}\label{dfn:Abel}
 The {\it Abel map\/} is the natural map of functors
	$$A_{X/S}(T)\:\Div_{X/S}(T)\to \Pic_{X/S}(T)$$
 defined by sending a relative effective divisor $D$ on $X_T/T$ to the
sheaf $\mc O_{\!X_T}(D)$.  The target $\Pic_{X/S}$ may be replaced by any
of its associated sheaves.  
 If $\IPic_{X/S}$ exists,
then the term {\it Abel map\/} may refer to the corresponding map of
schemes
	$$\I A_{X/S}\:\IDiv_{X/S}\to\IPic_{X/S}.$$
 \end{dfn}

\begin{ex}\label{ex:PQ-Abel}
  Assume $X/S$ is proper and flat with integral geometric fibers. Assume
$\IPic_{X/S}$ exists, and denote it by $P$.  View $\IDiv_{X/S}$ as a
$P$-scheme via the Abel map.  Assume a universal sheaf $\mc P$ exists,
and let $\mc Q$ be the sheaf on $P$ associated to $\mc P$ as in
Subsection~\ref{sb:Q}.  Show $\I P(\mc Q)=\IDiv_{X/S}$ as $P$-schemes.
 \end{ex}

\begin{thm}[Main]\label{th:main}
 Assume $f\:X\to S$ is projective Zariski locally over $S$, and is flat
with integral geometric fibers.

\tu{(1)} Then $\IPic_{X/S}$ exists, is separated and locally of finite
type over $S$, and represents $\Pic_{(X/S)\et}$.

\tu{(2)} If also $S$ is Noetherian and $X/S$ is projective, then
$\IPic_{X/S}$ is a disjoint union of open
 subschemes, each an increasing union of open quasi-projective
$S$-schemes.
 \end{thm}

\begin{proof}
 By~\cite[(0, 4.5.5), p.~106]{EGAG}, it is a local matter on $S$ to
represent a Zariski sheaf on the category of $S$-schemes.  Moreover, it
is also a local matter on $S$ to prove that an $S$-scheme is separated
and locally of finite type.  Hence, in order to prove (1), we may assume
$S$ is Noetherian and $X/S$ is projective.

Plainly an $S$-scheme is separated if it is a disjoint union of
separated open subschemes, or if it is an increasing union of
separated open subschemes.  Hence (1) follows from (2).

To prove (2), owing to Yoneda's lemma, we may view the category of
schemes as a full subcategory of the category of functors by identifying
a scheme $T$ with its functor of points.  Denote this functor too by $T$
in order to lighten the notation. And say that the functor is a scheme,
as well as that it is representable.  Also, set
 $$P:=\Pic_{(X/S)\et}.$$
  Note $P(T)= \Hom(T,P)$. 

Given a polynomial $\phi\in\bb Q[n]$, let $P^\phi\subset P$ be the
\'etale subsheaf associated to the presheaf whose $T$-points are
represented by the invertible sheaves $\mc L$ on $X_T$ such that we have
 \begin{equation}\label{eq:phi}
	\chi(X_t,\mc L^{-1}_t(n)) = \phi(n) \text{ for all $t\in T$}.
 \end{equation}
 Notice that this presheaf is well defined, because (\ref{eq:phi})
remains valid after any base change $p\:T'\to T$; indeed, for any $t'\in
T'$, for any $i$, and for any $n$, we have
	$$\uH^i\bigl(X_{t'},\mc L^{-1}_{t'}(n)\bigr)
   =\uH^i\bigl(X_{p(t')},\mc L^{-1}_{p(t')}(n)\bigr)\ox_{k_t}k_{t'}$$
 because cohomology commutes with flat base change by \cite[Prp.~III,
9.3, p.~255]{Ha83}.  Hence $P^\phi$ is well defined too.

Fix a map $T\to P$, and represent it by means of an \'etale covering
$p\:T'\to T$ and an invertible sheaf $\mc L'$ on $X_{T'}$.  Consider the
subset ${T'}^{\phi}\subset T'$ defined as follows:
       $${T'}^{\phi}:=\{\,t'\in T'\mid
       \chi\bigl(X_{t'},\mc L^{-1}_{t'}(n)\bigr) = \phi(n)\,\}.$$
 Then ${T'}^{\phi}$  is  open by \cite[7.9.11]{EGAIII2}.

Set $T^\phi:=p({T'}^{\phi})$.  Then $T^\phi\subset T$ is
open as ${T'}^{\phi}\subset T'$  is  open and $p$ is  \'etale.

Moreover, ${T'}^{\phi}=p^{-1}(T^\phi)$.  Indeed, let $t'\in
p^{-1}(T^\phi)$.  Say $p(t')=p(t'_1)$ where $t'_1\in {T'}^\phi$.  Now,
there is an \'etale covering $T''\to T'\x_TT'$ such that the two
pullbacks of $\mc L'$ to $X_{T''}$ are isomorphic.  Let $t''\in T''$
have image $t'\in T'$ under the first map $T''\to T'$ and have the image
$t'_1\in T'$ under the second map.  Then
	$$\chi\bigl(X_{t'},\mc L^{-1}_{t'}(n)\bigr)
	 = \chi\bigl(X_{t''},\mc L^{-1}_{t''}(n)\bigr)
	 = \chi\bigl(X_{t'_1},\mc L^{-1}_{t'_1}(n)\bigr) = \phi(n).$$
 Hence $t'\in {T'}^\phi$.  Thus ${T'}^{\phi}\supset
p^{-1}(T^\phi)$. Therefore, ${T'}^{\phi}=p^{-1}(T^\phi)$.

Furthermore, $T^\phi$ is (represents) the fiber product of functors
$P^\phi\x_PT$.  Indeed, to see they have the same $R$-points, let
$r\:R\to T$ be a map; form $R':=R\x_TT'$ and $r'\:R'\to T'$.  Suppose
$r$ factors through $T^\phi$.  Then $r'$ factors through ${T'}^{\phi}$.
So $R'\to T'\to P$ factors through $P^\phi$ essentially by definition.
Now, $R'\to R$ is an \'etale covering.  Hence $R\to T\to P$ factors
through $P^\phi$ since $P^\phi$ is an \'etale sheaf.

Conversely, suppose $R\to T\to P$ factors through $P^\phi$.  Then $R\to
P$ is defined by means of an \'etale covering $R''\to R$ and an
invertible sheaf $\mc L''$ on $X_{R''}$ such that $\chi(X_u,\mc {\mc
L''}^{-1}_u(n)) = \phi(n)$ for all $u\in R''$.  Since both $\mc L''$ and
$\mc L'$ define $R\to P$, there is an \'etale covering $R'''\to
R''\x_RR'$ such that the pullbacks of $\mc L''$ and $\mc L'$ to
$X_{R'''}$ are isomorphic.  Hence the image of $r'\:R'\to T'$ lies in
${T'}^{\phi}$.  But the latter is open.  Hence $r'$ factors through it.
Therefore,   $r\:R\to T$ factors through $T^\phi$.  Thus  $T^\phi$ and
$P^\phi\x_PT$ have the same $R$-points.

Let $\phi$ vary.  Plainly the ${T'}^{\phi}$ are disjoint and cover $T'$.
So the ${T}^{\phi}$ are disjoint and cover $T$.  Hence, by a general
result~\cite[(0, 4.5.4), p.~103]{EGAG}, if the $P^\phi$ are
(representable by) schemes, then $P$ is their disjoint union.  Thus it
remains to represent each $P^\phi$ by an increasing union of open
quasi-projective $S$-schemes.

Fix $\phi$.  Given $m\in \bb Z$, let $P^\phi_m\subset P^\phi$ be the
\'etale subsheaf associated to the presheaf whose $T$-points are
represented by the $\mc L$ on $X_T$ such that, in addition to
(\ref{eq:phi}), we have
 \begin{equation}\label{eq:4b}
	\tu R^if_{\smash{T}*}\mc L(n)=0
    \text{ for all }i\ge1\text{ and }n\ge m.
  \end{equation}
 Notice that this presheaf is well defined, because (\ref{eq:4b})
remains valid after any base change $p\:T'\to T$, as is shown next.

First, let's see that (\ref{eq:4b}) is equivalent to the following
condition:
 \begin{equation}\label{eq:4c}
	\uH^i(X_t,\mc L_t(n))=0
    \text{ for all }i\ge1,\text{ all }n\ge m,\text{ and all }t\in T.
 \end{equation}
 Indeed, (\ref{eq:4c}) implies (\ref{eq:4b}); in fact,  for any given
$i$, $t$ and $n$, if $\uH^i(X_t,\mc L_t(n))=0$, then $\tu
R^if_{\smash{T}*}\bigl(\mc L(n)\ox f_T^*\mc N\bigr)_t=0$ for all
quasi-coherent $\mc N$ on $T$ by \cite[7.5.3]{EGAIII2} or \cite[Cor.~2.1
p.~68]{OB72}.

Conversely, assume (\ref{eq:4b}).  Fix $t$ and $n$.  Let's proceed by
descending induction on $i$ to prove $\uH^i(X_t,\mc L_t(n))$ vanishes.
It vanishes for $i\gg 1$ by Serre's Theorem \cite[2.2.2]{EGAIII1}.
Suppose it vanishes for some $i\ge2$.  Then $\tu
R^if_{\smash{T}*}\bigl(\mc L(n)\ox f_T^*\mc N\bigr)_t$ vanishes for all
quasi-coherent $\mc N$ on $T$, as just noted.  So $\tu
R^{i-1}f_{\smash{T}*}\bigl(\mc L(n)\ox f_T^*\mc N\bigr)_t$ is right
exact in $\mc N$ owing to the long exact sequence of higher direct
images.  Therefore, by general principles, there is a natural
isomorphism of functors
 	$$ \tu R^{i-1}f_{\smash{T}*}(\mc L(n))_t\ox \mc N_t\risom
 	 \tu R^{i-1}f_{\smash{T}*}\bigl(\mc L(n)\ox f_T^*\mc N\bigr)_t.$$
 Since (\ref{eq:4b}) holds, both source and target vanish.
Taking $\mc N:=k_t$ yields the vanishing of $\uH^{i-1}(X_t,\mc
L_t(n))$. Thus (\ref{eq:4b}) implies (\ref{eq:4c}).

Finally, for any $t'\in T'$, any $i$, and any $n$,  we have
	$$\uH^i\bigl(X_{t'},\mc L_{t'}(n)\bigr)
	=\uH^i\bigl(X_{p(t')},\mc L_{p(t')}(n)\bigr)\ox_{k_t}k_{t'}$$
 because cohomology commutes with flat base change.  So  (\ref{eq:4c})
remains valid after the base change $p\:T'\to T$; whence, (\ref{eq:4b})
does too.  Thus the presheaf is well defined, and so $P^\phi_m$ is too. 

Arguing much as we did for $P^\phi\x_PT$, we find, given a map $T\to
P^\phi$, that, as $m$ varies, the products $P^\phi_m\x_{P^\phi}T$ form a
nested sequence of open subschemes of $T$, whose union is $T$.  In the
argument, the key change is in proving openness.  In place of
\cite[7.9.11]{EGAIII2}, we use the following part of Serre's
Theorem~\cite[2.2.2]{EGAIII1}: given a coherent sheaf $\mc F$ on a
projective scheme over a Noetherian ring $A$, there are only finitely
many $i\ge1$ and $n\ge m$ such that $H^i(\mc F(n))$ is nonzero, and all
these nonzero $A$-modules are finitely generated.  Hence, if there is a
prime ${\bf p}$ of $A$ such that $H^i(\mc F(n))_{\bf p}=0$ for all
$i\ge1$ and $n\ge m$, then there is an $a\notin{\bf p}$ such that
$H^i(\mc F(n))_a=0$ for all $i\ge1$ and $n\ge m$.

Therefore, again by \cite[(0, 4.5.4), p.~103]{EGAG}, it suffices to
represent each $P^\phi_m$ by a quasi-projective $S$-scheme.

Fix $\phi$ and $m$.  Set $\phi_0(n):=\phi(m+n)$.  Then there is an
isomorphism of functors $P^\phi_m \risom P^{\phi_0}_0$, which is defined
as follows.  First, define an endomorphism $\varepsilon$ of $\Pic_{X/S}$
by sending an invertible sheaf $\mc L$ on an $X_T$ to $\mc L(m)$.
Plainly $\varepsilon$ is an automorphism.  So $\varepsilon$ induces an
automorphism $\varepsilon^+$ of the associated sheaf $P$.  Plainly
$\varepsilon^+$ carries $P^\phi_m$ onto $P^{\phi_0}_0$.  Thus it
suffices to represent $P^{\phi_0}_0$ by a quasi-projective $S$-scheme.

The function $s\mapsto\chi(X_s,\mc O_{\!X_s}(n))$ is locally constant on
$S$ by \cite[7.9.11]{EGAIII2}.  Hence we may assume it is constant by
replacing $S$ by an open and closed subset.  Set $\psi(n):=\chi(X_s,\mc
O_{\!X_s}(n))$.

Consider the Abel map $A_{X/S}\:\Div_{X/S}\to P$.  Note that
$\Div_{X/S}$ is a scheme, in fact, an open subscheme of the Hilbert
scheme $\IHilb_{X/S}$, by Theorem~\ref{th:repDiv}.  Form the product
$\smash{P_0^{\phi_0}}\x_P\Div_{X/S}$.  It is a scheme $Z$, in fact, an
open subscheme of $\Div_{X/S}$, by what was proved above.  Set
$\theta(n):=\psi(n)-\phi_0(n)$.  Plainly $Z$ lies in
$\IHilb_{X/S}^\theta(n)$, which is projective over $S$.  Hence $Z$ is
quasi-projective.

Let's now prove  the projection $\alpha\:Z\to\smash{P_0^{\phi_0}}$
is a surjection of \'etale sheaves.  In other words, given a $T$ and a
$\lambda\in \smash{P_0^{\phi_0}}(T)$, we have to find an \'etale
covering $T_1\to T$ and a $\lambda_1\in Z(T_1)$ such that
$\alpha(\lambda_1)\in \smash{P_0^{\phi_0}}(T_1)$ is equal to the image of
$\lambda$.

Represent $\lambda$ by means of an \'etale covering $p\:T'\to T$ and an
invertible sheaf $\mc L'$ on $X_{T'}$.  Virtually by definition, the
product $T'\x_{\smash{P_0^{\phi_0}}}Z$ is equal to $\LinSys_{\mc
  L'/X_{T'}/T'}$.  So, by Theorem~\ref{th:LinSys}, this product is equal
to $\I P(\mc Q)$ where $\mc Q$ is the $\mc O_{T'}$-module associated to
$\mc L'$ as in Subsection~\ref{sb:Q}.  Now, $m=0$, so $\uH^1(X_t,\mc
L_t)=0$ for all $t\in T'$ owing to (\ref{eq:4c}).  Since (v) implies
(i) in Subsection~\ref{sb:Q}, therefore $\mc Q$ is locally free.

Hence $\I P(\mc Q)$ is smooth over $T'$.  So there exist an \'etale
covering $T_1\to T'$ and a $T'$-map $T_1\to \I P(\mc Q)$ by
\cite[17.16.3 (ii)]{EGAIV4}. Then the composition $T_1\to \I P(\mc Q)\to
Z\to \smash{P_0^{\phi_0}}$ is equal to the composition $T_1\to T'\to
T\to \smash{P_0^{\phi_0}}$.  In other words, the map $T_1\to Z$ is a
$\lambda_1\in Z(T_1)$ such that $\alpha(\lambda_1)\in
\smash{P_0^{\phi_0}}(T_1)$ is equal to the image of $\lambda$.  Since
the composition $T_1\to T'\to T$ is an \'etale covering, $\alpha$ is
thus a surjection of \'etale sheaves.

Plainly, the map $\alpha\:Z\to\smash{P_0^{\phi_0}}$ is defined by the
invertible sheaf associated to the universal relative effective divisor
on $X_Z/Z$.  So taking $T:=Z$ and $T':=T$ above, we conclude that the
product $Z\x_{\smash{P_0^{\phi_0}}}Z$ is a scheme and that the first
projection is smooth and projective.  Therefore, the theorem now results
from the following general lemma.
 \end{proof}

\begin{lem}\label{lm:qt}
 Let $\alpha\:Z\to P$ be a map of \'etale sheaves, and set $R:=Z\x_PZ$.
Assume $\alpha$ is a surjection, $Z$ is representable by a
quasi-projective $S$-scheme, $R$ is representable by an $S$-scheme, and
the first projection $R\to Z$ is representable by a smooth and proper
map.  Then $P$ is representable by a quasi-projective $S$-scheme, and
$\alpha$ is representable by a smooth map.
 \end{lem}

\begin{proof} To lighten the notation, let $Z$ and $R$ denote the
corresponding schemes as well.  Since the structure map $Z\to S$ is
quasi-projective, it is separated; whence, the first projection
$Z\x_SZ\to Z$ is too.  But the first projection $R\to Z$ is proper, and
factors naturally through a map $h\:R\to Z\x_SZ$.  Hence $h$ is proper.
But $h$ is a monomorphism; that is, $h$ is injective on $T$-points.
Therefore, $h$ is a closed embedding by \cite[8.11.5]{EGAIV3}.

Plainly, for each $S$-scheme $T$, the subset $R(T)\subset
Z(T)\x_{S(T)}Z(T)$ is the graph of an equivalence relation.  Also, the
map of schemes $R\to Z$ is flat and proper, and $Z$ is a
quasi-projective $S$-scheme.  It follows that there exist a
quasi-projective $S$-scheme $Q$ and a faithfully flat and projective map
$Z\to Q$ such that $R=Z\x_QZ$.  (In other words, a flat and proper
equivalence relation on a quasi-projective scheme is effective.)  In
fact, $R$ defines a map from $Z$ to the Hilbert scheme $\IHilb_{Z/S}$, and
its graph lies in the universal scheme as a closed subscheme, which
descends to $Q\subset \IHilb_{Z/S}$; for more details, see \cite[Thm.~(2.9),
p.~70]{AK80}.

Since $Z\to Q$ is flat, it is smooth if (and only if) its fibers are
smooth.  But these fibers are, up to extension of the ground field, the
same as those of $R\to Z$.  And $R\to Z$ is smooth by hypothesis.  Thus
$Z\to Q$ is smooth.

It remains to see that $Q$ represents $P$.  First, observe that $Z\to Q$
is (represents) a surjection of \'etale sheaves.  Indeed, given an
element of $Q(T)$, set $A:=Z\x_QT$.  Then $A\to T$ is smooth.  So there
exist an \'etale covering $T'\to T$ and a $T$-map $T'\to A$ by
\cite[17.16.3 (ii)]{EGAIV4}.  Then $T'\to A\to Z$ is an element of
$Z(T')$, which maps to the element of $Q(T')$ induced by the given
element of $Q(T)$.

Since $Z\to Q$ is a surjection of \'etale sheaves, $Q$ is, in the
category of \'etale sheaves, the coequalizer of the pair of maps
$R\rightrightarrows Z$ by Exercise~\ref{ex:epi} below, which is an
elementary exercise in general Grothendieck topology.  By the same
exercise, $P$ too is the coequalizer of this pair of maps.  But, in any
category, the coequalizer is unique up to unique isomorphism.  Thus $Q$
represents $P$.
 \end{proof}
 
\begin{ex}\label{ex:epi}
 Given a map of \'etale sheaves $F\to G$,  show it is a surjection if
and only if $G$ is the coequalizer of the pair of maps
$F\x_GF\rightrightarrows F$.
 \end{ex}

\begin{ex}\label{ex:q-proj}
 Assume $X/S$ is projective and flat, its geometric fibers are integral,
and $S$ is Noetherian.  Let $Z\subset \IPic_{X/S}$ be a subscheme of
finite type.  Show $Z$ is quasi-projective.
 \end{ex}

\begin{ex}\label{ex:Abelpr}
 Assume $f\:X\to S$ is projective Zariski locally over $S$, and is flat
with integral geometric fibers.  First, show that, if a universal sheaf
$\mc P$ exists, then the Abel map $\I A_{X/S}\:\IDiv_{X/S} \to
\IPic_{X/S}$ is projective Zariski locally over $S$.

Second, show that, in general, $\I A_{X/S}\:\IDiv_{X/S} \to \IPic_{X/S}$
is proper.  Proceed by reducing this case to the first: just use
$f\:X\to S$ itself to change the base.
 \end{ex}

\begin{ex}\label{ex:Abel}
Assume $f\:X\to S$ is flat and projective Zariski locally over $S$.
Assume its geometric fibers $X_k$ are integral curves of arithmetic
genus $\dim\uH^1(\mc O_{\!X_k})$ at least 1.  Let $X_0\subset X$ be the
open subscheme where $X/S$ is smooth.  Show there is a natural closed
embedding $A\:X_0\into\IPic_{X/S}$.
 \end{ex}

\begin{eg}\label{eg:Mumford}
 In Theorem~\ref{th:main}, the geometric fibers of $f$ are assumed to be
integral.  This condition is needed not only for the proof to work, but
also for the statement to hold.  The matter is clarified by the following
example, which is attributed to Mumford and is described in
\cite[p.~236-1]{FGA} and in \cite[p.~210]{BLR}.

Let $\bb R[\![t]\!]$ be the ring of formal power series, $S$ its
spectrum.  Let $X\subset \I P^2_S$ be the subscheme with inhomogeneous
equation $x^2+y^2=t$, and $f\:X\to S$ the structure map.  The generic
fiber $X_\sigma$ is a nondegenerate conic.  The special fiber $X_0$ is a
pair of conjugate lines; $X_0$ is irreducible and geometrically
connected.  So $f$ is flat by the implication (iii)$\Rightarrow$(i) of
Lemma~\ref{lm:ctn} with $D:=X$.  And $\mc O_{\!S}\risom f_*\mc O_{\!X}$
holds universally by Exercise~\ref{ex:gc&r}.  However, as we'll see,
$\IPic_{X/S}$ does {\it not\/} exist!

On the other hand, set $S':=\bb C[\![t]\!]$ and $X':=X\x S'$.  Then
$f_{S'}$ has sections: for example, one section $g'$ is defined by
setting
 $$x:=\sqrt{-1}\text{ and } y:=\sqrt{1+t}=1+(1/2)t-(1/8)t^2+\dotsb.$$
Hence all four relative Picard functors of $X'/S'$ are equal by the
Comparison Theorem, Theorem~\ref{th:cmp}.  Furthermore, as we'll see,
$\IPic_{X'/S'}$ exists, but is  {\it not\/}  separated!

In fact, $\IPic_{X'/S'}$ is a disjoint union of isomorphic open
nonseparated subschemes $S'_n$ for $n\in\bb Z$.  Each $S'_n$ is obtained
from $S'$ by repeating the origin infinitely often; more precisely, to
form $S'_n$, take a copy $S'_{a,b}$ of $S'$ for each pair $a,b\in\bb Z$
with $a+b=n$, and glue the $S'_{a,b}$ together off their closed
points.  Each $S'_{a,b}$ parameterizes a different degeneration of the
invertible sheaf of degree $n$ on the generic fiber $X'_\sigma$; the
degenerate sheaf has degree $a$ on the first line, and degree $b$ on the
second.  Also, complex conjugation interchanges the two lines, so
interchanges $S'_{a,b}$ and $S'_{b,a}$.

Suppose $\IPic_{X/S}$ exists.  Then $\IPic_{X/S}\x_SS' = \IPic_{X'/S'}$
by Exercise~\ref{ex:bschg}.  Since the closed points of $S'_{a,b}$ and
$S'_{b,a}$ are conjugate, they map to the same point of $\IPic_{X/S}$.
This point lies in an affine open subscheme $U$.  So the two closed
points lie in the preimage $U'$ of $U$ in $\IPic_{X'/S'}$.  But $U'$ is
affine since $U$, $S$ and $S'$ are affine.  However, if $a\neq b$, then
the two closed points are distinct, and so lie in no common affine $U'$.
We have a contradiction.  Thus $\IPic_{X/S}$ cannot exist.

Finally, let's prove $\Pic_{X'/S'}$ is representable by $\amalg S'_n$.
First, note $X'$ is regular; in fact, in $\I P^2_\bb C\x\I A^1_\bb C$,
the equation $x^2+y^2=t$ defines a smooth surface.  So on $X'$ every
reduced curve is an effective divisor.  In particular, consider these
three: the line $L:x=\sqrt{-1}y,\,t=0$, the line
$M:x=-\sqrt{-1}y,\,t=0$, and the image $A$ of the section $g$ of
$f_{S'}$ defined above.

Set $\mc P'_{a,b}:=\mc O_{\!X'}(bL+nA)$ where $n=a+b$.  The restriction to
the generic fiber $\mc P'_{a,b}|X'_\sigma$ has degree $n$ since $L\cap
X'_\sigma$ is empty and $A$ is the image of a section.  And $\mc
P'_{a,b}|M$ has degree $b$ since $A\cap M$ is empty and $L\cap M$ is a
reduced $\bb C$-rational point. And $\mc P'_{a,b}|L$ has degree $a$
since, in addition, $\mc O_{\!X'}(L)\ox \mc O_{\!X'}(M)\simeq\mc O_{\!X'}$ as
the ideal of $L+M$ is generated by $t$.  Lastly, every invertible sheaf
on $S'$ is trivial; fix a $g$-rigidification $\mc O_{S'}\risom
g^*\mc P'$, and use it to identify the two sheaves.

On $X'\x_{S'}\amalg S'_n$, form an invertible sheaf $\mc P'$ by placing
$\mc P'_{a,b}$ on $S'_{a,b}$; plainly, the $\mc P'_{a,b}$ patch
together.  It now suffices to show this: given any $S'$-scheme $T$ and
any invertible sheaf $\mc L$ on on $X'_T$, there exist a unique $S'$-map
$q\:T\to\bb \amalg S'_n$ and some invertible sheaf $\mc N$ on $T$ such
that $(1\x q)^*\mc P'\simeq \mc L\ox f_T^*\mc N$.

Replace $\mc L$ by $\mc L\ox f_T^*g_T^*\mc L^{-1}$.  Then $g_T^*\mc
L=\mc O_{\!T}$ since $g_T^*f_T^*=1$.  Hence, if $q$ and $\mc N$ exist, then
necessarily $\mc N\simeq\mc O_{\!T}$ since $g_T^*(1\x q)^*\mc
P'=q^*g^*\mc P'$ and $g^*\mc P'=\mc O_{S'}$.

Plainly, we may assume $T$ is connected.  Then the function
$s\mapsto\chi(X'_t,\mc L_t)$ is constant on $T$ by
\cite[7.9.5]{EGAIII2}.  Set $n:=\chi(X'_t,\mc L_t)-1$.

Fix $a,b$ with $a+b=n$.  Set
 $$\mc M:=\mc L^{-1}\ox(1\x\tau)^*\mc P'_{a,b}
 \text{ and } \mc N:=f_{\smash{T}*}\mc M.$$
 Plainly $g_T^*\mc M=\mc O_{\!T}$.  Form the natural map $u\:f_T^*\mc N\to
\mc M$.

Let $T_\sigma$ be the generic fiber of the structure map $\tau\:T\to
S'$.  Then $T_\sigma\subset T$ is open.  Let $t\in T_\sigma$.  Then
$X'_t$ is a nondegenerate plane conic with a rational point $A_t$.  So
$X'_t\simeq \I P^1_{k_t}$.  Hence $\mc L_t\simeq \mc O_{\!X'_t}(nA_t)$.
So $\mc M_t\simeq\mc O_{\!X'_t}$.  Hence $\uH^1(X'_t,\mc M_t)=0$ and
$\uH^0(X'_t,\mc M_t)=k_t$ by Serre's explicit computation
\cite[2.1.12]{EGAIII1}.

Therefore, $\mc N$ is invertible at $t$, and forming $\mc N$ commutes
with passing to $X'_t$, owing to the theory in Subsection~\ref{sb:Q}.
So forming $u\:f_T^*\mc N\to \mc M$ commutes with passing to $X'_t$ too.
But $u$ is an isomorphism on $X'_t$.  Hence $u$ is surjective along
$X'_t$ by Nakayama's lemma. But both source and target of $u$ are
invertible along $X'_t$.  Hence $u$ is an isomorphism along $X'_t$, so
over $T_\sigma$ as $t\in T_\sigma$ is arbitrary.  Now, $g_T^*\mc
M=\mc O_{\!T}$ and $g_T^*f_T^*=1$.  Hence $\mc N|T_\sigma= \mc
O_{T_\sigma}$.  Therefore, $\mc L|X'_{T_\sigma}=(1\x\tau)^*\mc
P'_{a,b}|X'_{T_\sigma}$.
 
Let $T_{a,b}'$ be the set of $t\in T$ such that $\tau(t)\in S'$ is the
closed point and $\mc L|L_t$ has degree $a$ and $\mc L|M_t$ has degree
$b$.  Fix $t\in T_{a,b}'$.  Then $\mc M|L_t$ has degree 0, so $\mc
M|L_t\simeq \mc O_{L_t}$.  Similarly, $\mc M|M_t\simeq \mc O_{M_t}$.
Consider the natural short exact sequence
	$$0\to\mc O_{L_t}(-1)\to\mc O_{\!X'_t}\to\mc O_{M_t}\to 0.$$
 Twist  by $\mc M$ and take cohomology.  Thus $\uH^1(X'_t,\mc M_t)=0$ and
$\uH^0(X'_t,\mc M_t)=k_t$.  Hence, on $X'_t$, the map $u$ becomes a map $\mc
O_{\!X'_t}\to \mc M_t$.  This map is surjective as it is surjective after
restriction to $L_t$ and to $M_t$.  So this map is an isomorphism.

As above, we conclude $u$ is an isomorphism on an open neighborhood $V'$
of $X'_t$.  Set $W':=f_{T'}(X'_T-V')$.  Since $f$ is proper, $W'$ is
open.  But $f^{-1}W'\subset V'$.  So $u$ is an isomorphism over $W'$.
Hence $O_{\!X'_{t'}}\risom \mc M_{t'}$ for all $t'\in W'$.  So $W'\subset
T_{a,b}'\cup T_\sigma$.

Set $T_{a,b}:=T_{a,b}'\cup T_\sigma$.  Then $T_{a,b}\subset T$ is open
as it contains a neighborhood of each of its points.  Furthermore, $u$
is an isomorphism over $T_{a,b}$.  Hence $\mc N|T_\sigma= \mc
O_{T_\sigma}$, again since $g_T^*\mc M=\mc O_{\!T}$ and $g_T^*f_T^*=1$.
Therefore, $\mc L|X'_{T_{a,b}}=(1\x\tau)^*\mc P'_{a,b}|X'_{T_{a,b}}$.

Let $q_{a,b}\:T_{a,b}\to \amalg S'_n$ be the composition of the
structure map $\tau\:T\to S'$ and the inclusion of $S'$ in $S'_n$ as
$S'_{a,b}$.  Plainly $(1\x q_{a,b})^*\mc P'=\mc L|X'_{T_{a,b}}$.
Plainly, as $a$ and $b$ and $n$ vary, the $q_{a,b}$ patch to a map
$q\:T\to\bb \amalg S'_n$ such that $(1\x q)^*\mc P'= \mc L$.  Plainly
this map $q$ is the only $S'$-map $q$ such that $(1\x q)^*\mc P'\simeq
\mc L\ox f_T^*\mc N$ for some invertible sheaf $\mc N$ on $T$.  Thus
$\amalg S'_n$ represents $\Pic_{X'/S'}$, and $\mc P'$ is a universal
sheaf.
 \end{eg}

\begin{ex}\label{ex:PE} 
 Assume $X=\I P(\mc E)$ where $\mc E$ is a locally free sheaf on $S$ and is
everywhere of rank at least 2.  Show $\IPic_{X/S}$ exists, and
represents $\Pic_{X/S}$; in fact, $\IPic_{X/S}= \bb Z_S$ where $\bb
Z_S$ stands for the disjoint union of copies of $S$ indexed by $\bb Z$.
 \end{ex}

\begin{ex}\label{ex:PfsCtd}
 Consider the curve $X/\bb R$ of Exercise~\ref{ex:Pfs}.  Show
$\IPic_{X/S}= \bb Z_\bb R$.
 \end{ex}

\begin{prp}\label{prp:lft}
 If $\IPic_{X/S}$ exists and represents $\Pic_{(X/S)\fppf}$, then
$\IPic_{X/S}$ is locally of finite type.
 \end{prp}
 \begin{proof}
 Set $P:=\Pic_{(X/S)\fppf}$.  Owing to \cite[8.14.2]{EGAIV3}, we just
need to check the following condition.  For every filtered inverse system
of affine $S$-schemes $T_i$, the natural map is a bijection:
 \begin{equation*}\label{equation:lft}
 \varinjlim P(T_i) \risom P\bigl(\varprojlim T_i\bigr).
 \end{equation*}

To check it is injective, set $T:=\varprojlim T_i$.  Fix $i$ and let
$\lambda_i\in P(T_i)$.  Represent $\lambda_i$ by a sheaf $\mc L_i'$ on
$X_{T_i'}$ where $T_i'\to T_i$ is an fppf covering.  Set
$T':=T_i'\x_{T_i}T$.  Set $\mc L':=\mc L_i|X_{T'}$.  Let $\lambda$ be
the image of $\lambda_i$ in $P(T)$.  Then $\mc L'$ represents $\lambda$.
  
Suppose $\lambda=0\in P(T)$.  Then there exists an fppf covering $T''\to
T'$ such that $\mc L'|X_{T''}\simeq \mc O_{\!X_{T''}}$.  It follows from
\cite[8.8.2, 8.10.5(vi), 11.2.6, 8.5.2(i), and 8.5.2.4]{EGAIV3} that
there exist a $j\ge i$ and an fppf covering $T_j''\to T_j'$ with
$T_j':=T_i'\x_{T_j}T$ such that $\mc L'|X_{T_j''}\simeq \mc
O_{\!X_{T_j''}}$.  So  $\lambda_i$ maps to $0\in P(T_j)$.  Thus the map
is injective. 

Surjectivity can be proved similarly. \end{proof}

\begin{rmk}\label{rk:exist}
 There are three important existence theorems for $\IPic_{X/S}$, which
refine Theorem~\ref{th:main}.  They were proved soon after it, and each
involves new ideas.
 
First, Mumford proved the following generalization of
Theorem~\ref{th:main}, and it fits in nicely with
Example~\ref{eg:Mumford}.
\begin{sbsthm}\label{th:Mumford}
 Assume $X/S$ is projective and flat, and its geometric fibers are
reduced and connected; assume the irreducible components of its ordinary
fibers are  geometrically irreducible.  Then $\IPic_{X/S}$ exists.
 \end{sbsthm}

Mumford stated this theorem at the bottom of Page viii in \cite{Mm66}.
He said the proof is a generalization of that \cite[Lects.~19--21]{Mm66}
in the case where $S$ is the spectrum of an algebraically closed field
and $X$ is a smooth surface.  That proof is based on his theory of
independent 0-cycles.  This theory is further developed in
\cite[pp.~23--28]{AK79} and used to prove \cite[Thm.~(3.1)]{AK79}, which
asserts the existence of a natural compactification of $\IPic_{X/S}$
when $X/S$ is flat and locally projective with integral geometric
fibers.

On the other hand, Grothendieck \cite[p.~236-1]{FGA} attributed to
Mumford a slightly different theorem, in which neither the geometric
fibers nor the ordinary fibers are assumed connected (see
\cite[p.~210]{BLR} also).  Grothendieck said the proof is based on a
refinement of Mumford's construction of quotients, and referred to the
forthcoming notes of a Harvard seminar of Mumford and Tate's, held in
the spring of 1962.

Mumford was kind enough, in November 2003, to mail the present author
his personal folder of handwritten notes from the seminar; the folder is
labeled, ``Groth--Mumford--Tate,'' and contains notes from talks by each
of the three, and notes written by each of them.  Virtually all the content
has appeared elsewhere; Mumford's contributions appeared in Mumford's
books \cite{Mm65}, \cite{Mm66}, and \cite{Mm70}.

The notes contain, in Mumford's hand, a precise statement of the theorem
and a rough sketch of the proof.  This statement too is slightly
different from that of Theorem ~\ref{th:Mumford}: he crossed out the
hypothesis that the geometric fibers are connected; and he made the
weaker assumption that the ordinary fibers are connected.  The proof is
broadly like his proof in \cite[Lects.~19--21]{Mm66}.
 
\smallskip
Second, Grothendieck proved this theorem of ``generic
representability.''
  \begin{sbsthm}\label{th:genrep}
 Assume $X/S$ is proper, and $S$ is integral.  Then there exists a
nonempty open subset $V\subset S$ such that $\IPic_{X_V/V}$ exists,
represents $\Pic_{(X_V/V)\fppf}$, and is a disjoint union of open
quasi-projective subschemes.
 \end{sbsthm}

 A particularly important special case is covered by the
following corollary.
  \begin{sbscor}\label{cor:algsch}
  Assume $S$ is the spectrum of a field $k$, and $X$ is complete.  Then
$\IPic_{X/k}$ exists, represents $\Pic_{(X/k)\fppf}$, and is a disjoint
union of open quasi-projective subschemes.
 \end{sbscor}

Before Grothendieck discovered Theorem~\ref{th:genrep}, he
\cite[Cor.~6.6, p.~232-17]{FGA} proved Corollary~\ref{cor:algsch}
assuming $X/k$ is projective.  To do so, he developed a method of
``relative representability,'' by which Theorem~\ref{th:main} implies
the existence of the Picard scheme in other cases.  The method
incorporates a ``d\'evissage'' of Oort's \cite[\S 6]{Oo62}; the latter
yields the Picard scheme as an extension of a group subscheme of
$\IPic_{X_{\red}/k}$ by an affine group scheme.

In \cite[Rem.~6.6, p.~232-17]{FGA}, Grothendieck said it is ``extremely
plausible'' that Corollary~\ref{cor:algsch} holds in general, and can be
proved by extending the method of relative representability, so that it
covers the case of a surjective map $X'\to X$ with $X'$ projective, such
as the map provided by Chow's lemma \cite[Thm.~5.6.1]{EGAII}.
Furthermore, he conjectured, in \cite[Rem.~6.6, p.~232-17]{FGA}, that,
for any surjective map $X'\to X$ between proper schemes over a field,
the induced map on Picard schemes is affine.  He said he was led to the
conjecture by considerations of ``nonflat descent,'' a version of
descent theory where the maps are not required to be flat, but the
objects are.

Thanking Grothendieck for help, Murre \cite[(II.15)]{Mr63} gave the
first proof of the heart of Corollary~\ref{cor:algsch}: if $X/k$ is
proper, then $\IPic_{X/k}$ exists and
is locally of finite type.  Murre did not use the method of relative
representability.  Rather, he identified seven conditions
\cite[(I.2.1)]{Mr63} that are necessary and sufficient for the
representability of a functor from schemes over a field to Abelian
groups.  Then he showed the seven are satisfied by the Picard functor
localized in the fpqc topology.  To handle the last two conditions, he
used Chow's lemma and Theorem~\ref{th:main}.

In the meantime, Grothendieck had proved Theorem~\ref{th:genrep},
according to Murre \cite[p.~5]{Mr63}.  Later, the proof appeared in two
parts.  The first part established a key intermediate result, the
following theorem.
 \begin{sbsthm}\label{th:flatten}
  Assume $X/S$ is proper.  Let $\mc F$ be a coherent sheaf on $X$, and
$S_{\mc F}\subset S$ the subfunctor of all $S$-schemes $T$ such that
${\mc F}_T$ is $T$-flat.  Then $S_{\mc F}$ is representable by an
unramified $S$-scheme of finite type.
 \end{sbsthm}
 Murre sketched Grothendieck's proof of this theorem in \cite[Cor.~1,
   p.~294-11]{Mr65}.  The proof involves identifying and checking eight
 conditions that are necessary and sufficient for the representability
 of a functor by a separated and unramified $S$-scheme locally of finite
 type.  As in Murre's proof of Corollary~\ref{cor:algsch}, a key step is
 to show the functor is ``pro-representable'': there exist certain
 natural topological rings, and if the functor is representable by $Y$,
 then these rings are the completions of the local rings at the points
 of $Y$ that are closed in their fibers over $S$.

In the second part of Grothendieck's proof of Theorem~\ref{th:genrep},
the main result is the following theorem of relative representability.
 \begin{sbsthm}\label{th:rrep}
  Assume $S$ is integral. Let $X'\to X$ be a surjective map of
proper $S$-schemes.  Then there is a nonempty open subset $V\subset
S$ such that the map
	$$\Pic_{(X_V/V)\fppf}\to\Pic_{(X'_V/V)\fppf}$$
 is representable by quasi-affine maps of finite type.
 \end{sbsthm}
 \noindent In other words, for every $S$-scheme $T$ and map
$T\to\Pic_{(X'_V/V)\fppf}$, the fibered product
 $\Pic_{(X_V/V)\fppf}\x_{\Pic_{(X'_V/V)\fppf}} T$ is
representable and its projection to $T$ is a quasi-affine map of finite
type between schemes.

Raynaud \cite[Thm.~1.1]{Ra71} gave Grothendieck's proof of
Theorem~\ref{th:rrep}; the main ingredients are, indeed, Oort
d\'evissage and nonflat descent.  As a first consequence, Raynaud
\cite[Cor.~1.2]{Ra71} derived Theorem~\ref{th:genrep}.  In order to show
$\IPic_{X/S}$ is a disjoint union of open quasi-projective subschemes,
he used a finiteness theorem for $\IPic_{X/S}$, which Grothendieck had stated
under (v) in \cite[p.~C-08]{FGA}, and which is proved below as
Theorem~\ref{th:Ptaufin}.  As a second consequence of
Theorem~\ref{th:rrep}, Raynaud \cite[Cors.~1.5]{Ra71} established
Grothendieck's conjecture that, over a field, a surjective map of proper
schemes induces an affine map on Picard schemes.

\smallskip
The third important existence theorem for $\IPic_{X/S}$ is due to
M.~Artin, whose work greatly clarifies the situation.  Artin proved
$\IPic_{X/S}$ exists when it should, but not as a scheme.  Rather, it
exists as a more general object, called an ``algebraic space,'' which is
closer in nature to a (complex) analytic space.

   Grothendieck \cite[Rem.~5.2, p.~232-13]{FGA} had said: ``it is not
 ruled out that $\IPic_{X/S}$ exists whenever $X/S$ is proper and flat
 and such that $\mc O_{\!S}\risom f_*\mc O_{\!X}$ holds universally.  At least,
 this statement is proved in the context of analytic spaces when $X/S$
	   is, in addition, projective.''  Mumford's example,
  Example~\ref{eg:Mumford}, shows the statement is false for schemes;
    Artin's theorem shows the statement holds for algebraic spaces.

Algebraic spaces were introduced by Artin \cite[\S 1]{Ar68}, and the
theory developed by himself and his student Knutson \cite{Kn71}. The
spaces are constructed by gluing together schemes along open subsets
that are isomorphic locally in the \'etale topology.  Over $\bb C$, these
open sets are locally analytically isomorphic; so an algebraic space is
a kind of complex analytic space.  In general and more formally, an
algebraic space is the quotient in the category of \'etale sheaves of a
scheme divided by an \'etale equivalence relation.

Artin \cite[Thm.~3.4, p.~35]{Ar67}, proceeding in the spirit of
Grothendieck and Murre, identified five conditions on a functor that are
necessary and sufficient for it to be a locally separated algebraic
space that is locally of finite type over a field or over an excellent
Dedekind domain.  In the proof of necessity, a key new ingredient is
Artin's approximation theorem; it implies that the topological rings
given by pro-representability can be algebraized.  The resulting schemes
are then glued together via an \'etale equivalence relation.

By checking that the five conditions hold, Artin \cite[Thm.~7.3, p.~67]{Ar67}
proved the following theorem.
 \begin{sbsthm}\label{th:algsp}
  Let $f\:X\to S$ be a flat, proper, and finitely presented map of
algebraic spaces.  Assume forming $f_*\mc O_{\!X}$ commutes with changing
$S$ (or $f$ is ``cohomologically flat in dimension zero'').  Then
$\IPic_{X/S}$ exists as an algebraic space, which is locally of finite
presentation over $S$.
 \end{sbsthm}
 \noindent Since $f$ is finitely presented, the statement can be reduced
to the case where $S$ is locally of finite type over a field or over an
excellent Dedekind domain.  Artin's proof of Theorem~\ref{th:algsp} is
direct: it involves no reduction whatsoever to Theorem~\ref{th:main}.

On the other hand, it is not known, in general, whether
Theorem~\ref{th:algsp} implies Theorem~\ref{th:main}, which asserts
$\IPic_{X/S}$ is a scheme.  However, more is known over a {\it field}.
Indeed, given an algebraic space that is locally of finite type over a
field, Artin \cite[Lem.~4.2, p.~43]{Ar67} proved this: if the space is a
sheaf of groups, then it is a group scheme.  Thus Theorem~\ref{th:algsp}
implies the heart of Corollary ~\ref{cor:algsch}.

Thus, in Theorem~\ref{th:algsp}, the fibers over $S$ of the algebraic
space $\IPic_{X/S}$ are schemes; they are the Picard schemes of the
fibers $X_s$, though the $X_s$ need not be schemes.  In particular, if
$S$ and $X$ are schemes\emdash so the $X_s$ are too\emdash then the
$\IPic_{X_s/k_s}$ are schemes.  Furthermore, then they form a family;
its total space $\IPic_{X/S}$ is an algebraic space, but not necessarily
a scheme.
 \end{rmk}

\section{The connected component of the identity}\label{sc:Pic0}
 Having treated the existence of the Picard scheme $\IPic_{X/S}$, we now
turn to its structure.  In this section, we study the union
$\IPicz_{X/S}$ of the connected components of the identity element,
$\IPicz_{X_s/k_s}$, for $s\in S$.  We establish a number of basic
properties, especially when $S$ is the spectrum of a field.

It is remarkable how much we can prove about $\IPicz_{X/S}$ formally, or
nearly so, from general principles.  Notably, we can do without the
finiteness theorems proved in the next section.  In order to emphasize
the formal nature and corresponding generality of the arguments,
most of the results are stated with the general hypothesis that
$\IPic_{X/S}$ exists instead of with specific hypotheses that imply it
exists.

\begin{lem}\label{lem:agps}
 Let $k$ be a field, and $G$ a group scheme locally of finite
type.  Let $G^0$ denote the connected component of the identity element $e$.

\tu{(1)} Then $G$ is separated.

\tu{(2)} Then $G$ is smooth if it has a geometrically
reduced open subscheme.

 \tu{(3)} Then $G^0$ is an open and closed group subscheme of finite
type; it is geometrically irreducible; and forming it commutes with
extending $k$.
 \end{lem}
 \begin{proof}
 Since $e$ is a $k$-point, it is closed.  Define a map $\alpha\:G\x G\to
G$ on $T$-points by $\alpha(g,h):=gh^{-1}$.  Then $\alpha^{-1}e\subset
G\x G$ is a closed subscheme.  Its $T$-points are just the pairs $(g,g)$
for $g\in G(T)$; so $\alpha^{-1}e$ is the diagonal.  Thus (1) holds.

Suppose $G$ has a geometrically reduced open subscheme $V$.  To prove
$G$ is smooth, we may replace $k$ by its algebraic closure.  Then $V$
contains a nonempty smooth open subscheme $W$.  Furthermore, given any
two closed points $g,\,h\in G$, there is an automorphism of $G$ that
carries $g$ to $h$, namely, multiplication by $g^{-1}h$.  Taking $g\in
W$, we conclude that $G$ is smooth at $h$.  Thus (2) holds.

Consider (3).  By definition, $G^0$ is the largest connected subspace
containing $e$.  But the closure of a connected subspace is plainly
connected; so $G^0$ is closed.  By \cite[6.1.9]{EGAI}, in any locally
Noetherian topological space, the connected components are open.  Thus
$G^0$ is open too.

Since $G^0$ is connected and has a $k$-point, $G^0$ is geometrically
connected by \cite[4.5.14]{EGAIV2}.  Thus forming $G^0$ commutes with
extending $k$.

Furthermore, $G^0\x G^0$ is connected by \cite[4.5.8]{EGAIV2}.
So $\alpha(G^0\x G^0)\subset G$ is connected, and contains $e$,
so lies in $G^0$.  Thus $G^0$ is a subgroup.

To prove $G^0$ is geometrically irreducible and quasi-compact, we may
replace $k$ by its algebraic closure.  Then $G^0_{\red}\x G^0_{\red}$ is
reduced by \cite[4.6.1]{EGAIV2}.  Hence $\alpha$ induces a map from
$G^0_{\red}\x G^0_{\red}$ into $G^0_{\red}$.  So $G^0_{\red}$ is a
subgroup.  Thus we may replace $G^0$ by $G^0_{\red}$.

Since $G^0$ is reduced and $k$ is algebraically closed, $G^0$ contains a
nonempty smooth affine open subscheme $U$.  Take arbitrary $k$-points
$g\in U$ and $h\in G^0$.  Then $hg^{-1}U$ is smooth and open, and
contains $h$.  Hence $G^0$ is irreducible locally at $h$.  But $h$ is
arbitrary, and $G^0$ is connected.  So $G^0$ is irreducible by
\cite[6.1.10]{EGAI}.

Since $G^0$ is irreducible, its open subschemes $U$ and $hU$ meet.  So
their intersection contains a $k$-point $g_1$ since $k$ is algebraically
closed. Then $g_1=hh_1$ for some a $k$-point $h_1\in U$.  Then
$h=g_1h_1^{-1}$.  But $h\in G^0$ is arbitrary.  Hence $\alpha(U\x
U)=G^0$.  Now, $U\x U$ is affine, so quasi-compact.  Hence $G^0$ is
quasi-compact.  By hypothesis, $G$ is locally of finite type. Hence
$G^0$ is of finite type.  Thus (3) holds.
 \end{proof}

\begin{rmk}\label{rk:Gred}
 Let $G$ be a group scheme of finite type over a field, $G_{{\red}}$ its
reduction.  Then $G_{{\red}}$ need not be a group subscheme, because
$G_{{\red}}\x G_{{\red}}$ need not be reduced.  Waterhouse
\cite[p.~53]{Wa79} gives two conditions equivalent to reducedness when $G$
is finite in Exercise~9, and he gives a counterexample in Exercise~10.
 \end{rmk}

\begin{prp}\label{prp:pic0}
  Assume $S$ is the spectrum of a field $k$.  Assume $\IPic_{X/k}$
exists and represents $\Pic_{(X/k)\fppf}$.  Then $\IPic_{X/k}$ is
separated, and it is smooth if it has a geometrically reduced open
subscheme.  Furthermore, the connected component of the identity
$\IPicz_{X/k}$ is an open and closed group subscheme of finite type; it
is geometrically irreducible; and forming it commutes with extending
$k$.
 \end{prp}
\begin{proof}
 This result follows formally from Proposition~\ref{prp:lft}
and Lemma~\ref{lem:agps}.
 \end{proof}

\begin{thm}\label{th:qpp&p}
 Assume $S$ is the spectrum of a field $k$.  Assume $X/k$ is projective
and geometrically integral.  Then $\IPicz_{X/k}$ exists and is
quasi-projective.  If also $X/k$ is geometrically normal, then
$\IPicz_{X/k}$ is projective.
 \end{thm}
 \begin{proof} 
 Theorem~\ref{th:main} implies $\IPic_{X/k}$ exists and represents
$\Pic_{(X/k)\et}$, so $\Pic_{(X/k)\fppf}$.  Hence $\IPicz_{X/k}$ exists
and is of finite type by Proposition~\ref{prp:pic0} (in fact, here
Proposition~\ref{prp:lft} is logically unnecessary since $\IPic_{X/k}$
is locally of finite type by Theorem~\ref{th:main}).  So $\IPicz_{X/k}$
is quasi-projective by Exercise~\ref{ex:q-proj}.

Suppose $X$ is also geometrically normal.  Since $\IPicz_{X/k}$ is
quasi-projective, to prove it is projective, it suffices to prove it is
proper.  By Lemma~\ref{lem:agps}, forming $\IPicz_{X/k}$ commutes
with extending $k$.  And by \cite[2.7.1(vii)]{EGAIV2}, a $k$-scheme is
complete if (and only if) it is after extending $k$.  So we may, and do,
assume $k$ is algebraically closed.

Recall the structure theorem of Chevalley and Rosenlicht for algebraic
groups, or reduced connected group schemes of finite type over $k$; see
\cite[Thm.~1.1, p.~3]{Co02}.  The theorem says that every algebraic
group is an extension of an Abelian variety (or complete algebraic
group) by a linear (or affine) algebraic group.  Recall also that every
solvable linear algebraic $k$-group is triangularizable (the
Lie--Kolchin theorem); so, if it's nontrivial, then it contains a copy
of the multiplicative group or of the additive group; see \cite[(10.5)
and (10.2)]{Bo69}.  Now,  $\IPic_{X/k}$ is commutative, so solvable.
Hence it suffices to show that, if $T$ denotes the
affine line minus the origin, then every $k$-map $t\:T\to
(\IPicz_{X/k})_{\red}$ is constant.

Since $k$ is algebraically closed, $X/k$ has a section.  So $t$
arises from an invertible sheaf $\mc L$ on $X\x T$ by the Comparison
Theorem, Theorem~\ref{th:cmp}.  Since $X\x T$ is integral, there is a
divisor $D$ such that $\mc O(D)=\mc L$ by \cite[Ex.~II, 6.15,
p.~145]{Ha83}.

Form the projection $p\:X\x T \to X$.  Restrict $\mc L$ to its generic
fiber.  This restriction is trivial as $T$ is an open subset of the
line.  So there is a rational function $\phi$ on $X\x T$ such that
$(\phi)+D$ restricts to the trivial divisor.  Let $s\:X\to X\x T$ be a
section.  Set $E:=s^*((\phi)+D)$; then $E$ is a well-defined divisor on
$X$.  Plainly, $p^*E$ and $(\phi)+D$ coincide as cycles; whence, they
coincide as divisors by \cite[Prp.~(3.10), p.~139]{AK70}, since $X\x T$
is normal.  Therefore, $\mc L=p^*\mc O(E)$.  Thus $t\:T\to \IPic_{X/k}$
is constant.
 \end{proof}

\begin{cor}\label{cor:Poincare}
 Assume $S$ is the spectrum of an algebraically closed field $k$.
Assume $X$ is projective and integral.  Set $P:=\IPicz_{X/k}$, and let
$\mc P$ be the restriction to $X_P$ of a Poincar\'e sheaf.  Then a\/
{\rm Poincar\'e family} $W$ exists; by definition, $W$ is a relative
effective divisor on $X_P/P$ such that
	$$\mc O_{X_P}\bigl(W-(W_0\x P)\bigr) \simeq \mc P\ox f_P^*\mc N$$
 where $W_0$ is the fiber over $0\in P$ and where $\mc N$ is an
invertible sheaf on $P$.
 \end{cor}
\begin{proof}
 Note that $P$ exists and is quasi-projective by Theorem~\ref{th:qpp&p}
and that $\mc P$ exists by Exercise \ref{ex:gc&r} and Exercise
\ref{ex:univshf}.  Since $P$ is Noetherian, Serre's Theorem
\cite[2.2.1]{EGAIII1} implies there is an $N$ such that
$R^if_{\smash{P}*}\mc P(n)=0$ for all $i>0$ and $n\ge N$.  Recall that
(\ref{eq:4b}) implies (\ref{eq:4c}); similarly, $\uH^i(\mc P_t(n))=0$
for all $t \in P$.  Fix an $n\ge N$ such that $\dim \uH^0(\mc
O_X(n))>\dim P$.

Say $\lambda\in \IPic_{X/k}$ represents $\mc O_X(n)$.  Form the
automorphism of $\IPic_{X/k}$ of multiplication by $\lambda$.  Plainly,
$P$ is carried onto the connected component, $P'$ say, of $\lambda$.
Let $q\:P\risom P'$ be the induced isomorphism.  Let $\mc P'$ be the
restriction to $X_{P'}$ of a Poincar\'e sheaf.  Plainly, $(1\x q)^*\mc
P'\simeq\mc P(n)\ox f_P^*\mc N$ for some invertible sheaf $\mc N$.

By Exercise~\ref{ex:PQ-Abel}, there is a coherent sheaf
$\mc Q$ on $P$ such that $\I P(\mc Q)=\IDiv_{X/S}$.  Moreover, $\mc
Q|P'$ is locally free of rank $\dim \uH^0(\mc O_X(n))$ owing to
Subsection~\ref{sb:Q}.  So $\mc Q|P'$ is of rank at least $1+\dim P$.
Now, there is an $m$ such that the sheaf $\SHom(\mc Q|P', \,\mc O_P)(m)$
is generated by finitely many global sections; so a general linear
combination of them vanishes nowhere by a well-known lemma
\cite[p.~148]{Mm66} attributed to Serre.  Hence
 there is a surjection $\mc Q|P'\onto\mc O_P(m)$.
  Correspondingly, there is a
$P'$-map $h'\:P'\to \I P(\mc Q|P')$; in other words, $h'$ is a section
of the restriction over $P'$ of the Abel map $\I A_{X/S}\:\IDiv_{X/S}\to
\IPic_{X/S}$.

Let $W'\subset X_{P'}$ be the pullback under $1\x h'$ of the universal
relative effective divisor.  Then $\mc O_{X_{P'}}(W')=\mc P'$ since $h'$
is a section of $\I A_{X/S}|P'$.  So, in particular, $\mc
O_X(W'_\lambda)= \mc O_X(n)$.  Set $W:=(1\x q)^{-1}W'$.  Plainly, $W$ is
a Poincar\'e family.
 \end{proof}

\begin{rmk}\label{rk:agqp}
  More generally, Theorem~\ref{th:qpp&p} holds whenever $X/k$ is proper,
whether $X$ is integral or not.  The proof is essential the same, but
uses Corollary~\ref{cor:algsch} in place of Theorem~\ref{th:main}.  In
fact, the proof of the quasi-projectivity assertion is easier, and does
not require Proposition~\ref{prp:pic0}, since $\IPicz_{X/k}$ is given as
contained in a quasi-projective scheme.  On the other hand, the proof of
the projectivity assertion requires an additional step: the reduction,
when $k$ is algebraically closed, to the case where $X$ is irreducible.
Here is the idea: since $X$ is normal, $X$ is the disjoint union of its
irreducible components $X_i$, and plainly
$\IPic_{X/k}=\prod\IPic_{X_i/k}$; hence, if the $\IPicz_{X_i/k}$ are
complete, so is $\IPicz_{X/k}$.

Chow \cite[Thm.~p.~128]{Ch57} proved every algebraic group \emdash indeed,
every homogeneous variety \emdash is quasi-projective.  Hence, in
Lemma~\ref{lem:agps} , if $G$ is reduced, then $G^0$ is
quasi-projective.  In characteristic 0, remarkably $G$ is smooth, so
reduced; this result is generally attributed to Cartier, and is proved
in \cite[p.~167]{Mm66}.

It follows that Theorem~\ref{th:qpp&p} holds in characteristic 0
whenever $X$ is a proper algebraic $k$-space.  Indeed, the
quasi-projectivity assertion holds in view of the preceding discussion
and of the discussion at the end of Remark~\ref{rk:exist}.  The proof of
the projectivity assertion has one more complication: it is necessary to
work with an \'etale covering $U\to X$ where $U$ is a scheme and with an
invertible sheaf $\mc L$ on $U\x T$.  However, the proof shows that, in
arbitrary characteristic, if $X$ is a proper and normal algebraic
$k$-space, then $\IPicz_{X/k}$ is proper.
 \end{rmk}

\begin{ex}\label{ex:sm=>pr}
 Assume $X/S$ is projective and smooth, its geometric fibers are
irreducible, and $S$ is Noetherian. Using the Valuative Criterion
\cite[Thm.~4.7, p.~101]{Ha83} rather than the Chevalley--Rosenlicht
structure theorem, prove that a closed subscheme $Z\subset\IPic_{X/S}$ is
projective over $S$ if it is of finite type.
 \end{ex}

\begin{rmk}\label{rmk:RamSam}
 There are three interesting alternative proofs of the second assertion
of Theorem~\ref{th:qpp&p}.  The first alternative uses
Exercise~\ref{ex:sm=>pr}.  It was sketched by Grothendieck
\cite[p.~236-12]{FGA}, and runs basically as follows.  Proceed by
induction on the dimension $r$ of $X/k$.  If $r=0$, then $\IPicz_{X/k}$
is $S$, so trivially projective.  If $r=1$, then $\IPicz_{X/k}$ is
projective by Exercise~\ref{ex:sm=>pr}.

Suppose $r\ge2$.  As in the proof of the theorem, reduce to the case
where $k$ is algebraically closed.  Let $Y$ be a general hyperplane
section of $X$.  Then $Y$ too is normal \cite[Thm.~$7'$, p.~376]{Se50}.
Plainly the inclusion $\vf\:Y\into X$ induces a map
	$$\vf^*:\IPicz_{X/k}\to\IPicz_{Y/k}.$$
 By induction, $\IPicz_{Y/k}$ is projective.  So $\IPicz_{X/k}$ is
projective too if $\vf^*$ is finite.

In order to handle $\vf^*$, Grothendieck suggested using a version of
the ``known equivalence criteria.''  In this connection, he
\cite[p.~236-02]{FGA} announced that \cite{SGA2} contains some key
preliminary results, which must be combined with the existence theorems
for the Picard scheme.  In fact, \cite[Cor.~3.6, p.~153]{SGA2} does
directly imply $\vf^*$ is injective for $r\ge3$.  Hence $\vf^*$ is
generically finite.  So  $\vf^*$ is finite since it is homogeneous.

For any $r\ge2$, \cite[Lem.~3.11, p.~639, and Rem.~3.12, p.~640]{Kl73}
assert $\ker\vf^*$ is finite and unipotent, so trivial in characteristic
0; however, the proofs in \cite{Kl73} require the compactness of
$\IPicz_{X/k}$.  It would be good to have a direct proof of this
assertion also when $r=2$, a proof in the spirit of \cite{SGA2}.

In characteristic 0, Mumford \cite[p.~99]{Mm67} proved as follows
 $\ker\vf^*$ vanishes.  Suppose not. First of all, $\ker\vf^*$ is
reduced by Cartier's theorem.  So $\ker\vf^*$ contains a point $\lambda$
of order $n>1$.  And $\lambda$ defines an unramified Galois cover $X'/X$
with group ${\bb Z}/n$.  Set $Y':=Y\x_XX'$.  Then $Y'$ is a disjoint
union of $n$ copies of $Y$ because $\lambda\in \ker\vf^*$.  On the other
hand, $Y'$ is ample since $Y$ is; hence, $Y'$ is connected (by
Corollary~\ref{cor:conn} for example).  We have a contradiction.  Thus
$\ker\vf^*$ vanishes.

In characteristic 0, the injectivity of $\vf^*$ also follows from the
Kodaira Vanishing Theorem.  Indeed, as just noted, $\ker\vf^*$ is
reduced.  So $\vf^*$ is injective if and only if its differential is
zero.  Now, this differential is, owing to Theorem~\ref{thm:tgtsp}
below, equal to the natural map
	$$\tu H^1(\mc O_{\!X})\to\tu H^1(\mc O_{\!Y}).$$
 Its kernel is equal to $\tu H^1(\mc O_{\!X}(-1))$ owing to the long exact
sequence of cohomology.

In characteristic 0, if $X$ is smooth, then $\tu H^1(\mc O_{\!X}(-1))$
vanishes by the Kodaira Vanishing Theorem.  If also $r=2$, then the dual
group $\tu H^1(\Omega^2_{\!X}(1))$ vanishes by the theorem on the
regularity of the adjoint system, which was proved by Picard in 1906
using Abelian integrals and by Severi in 1908 using algebro-geometric
methods.  For more information, see \cite[pp.~181, 204--206]{Za35} and
\cite[pp.~94--97]{Mm67}.

For any $r\ge2$, there are, as Grothendieck \cite[p.~236-12]{FGA}
suggested, finitely many smooth irreducible curves $Y_i\subset X$
such that the induced map is injective:
	 $$\IPic_{X/k}^0 \to \textstyle\prod \IPic_{Y_i/k}^0.$$
 So, again, since the $\IPic_{Y_i/k}^0$ are projective, $\IPicz_{X/k}$
is projective too.

To find the $Y_i$, use the final version of the ``equivalence criteria''
proved by Weil \cite[Cor.~2, p.~159]{We54}; it says in other words that,
if $W/T$ is the family of smooth 1-dimensional linear-space sections of
$X$ (or even a nonempty open subfamily), then the induced map is
injective:
	$$\IPic_{X/k}^0(k)\to\IPic_{W/T}^0(T).$$

 For each finite set $F$ of $k$-points of $T$, let $K_F$ be the kernel
of the map
     $$\IPic_{X/k}^0 \to\textstyle \prod_{t\in F} \IPic_{W_t/k}.$$
 Since $K_F$ is closed, we may assume, by Noetherian induction, that
$K_F$ contains no strictly smaller $K_G$.  Suppose $K_F$ has a nonzero
$k$-point.  It yields a nonzero $T$-point of $\Pic_{W/T}$, so a nonzero
$k$-point of $\IPic_{W_t/k}$ for some $k$-point $t$ of $T$.  Let $G$ be the
union of $F$ and $\{t\}$.  Then $K_G$ is strictly smaller than $K_F$.  So
$K_F = 0$.  Take the $Y_i$ to be the $W_t$ for $t\in F$. 

In characteristic 0 or if $r=2$, another way to finish is to take a
desingularization $\vt\:X'\to X$.  Since $X$ is normal, a divisor $D$ on
$X$ is principal if $\vt^*D$ is principal; hence, the induced map
	$$\vt^*:\IPicz_{X/k}\to\IPicz_{X'/k}$$
  is injective.  But $\IPicz_{X'/k}$ is projective by
Exercise~\ref{ex:sm=>pr}.  So $\IPicz_{X/k}$ is projective too.

The second alternative proof of Theorem~\ref{th:qpp&p} is similar to the
proof in Answer~\ref{ex:sm=>pr}.  (It too may be due to
Grothendieck\emdash see \cite[21.14.4, iv]{EGAIV4}\emdash but was
indicated to the author by Mumford in a private communication in 1974).
Again we reduce to the case where $k$ is algebraically closed.  Then,
using a more refined form of the Valuative Criterion (obtained modifying
\cite[Ex.~4.11, p.~107]{Ha83} slightly), we need only check this
statement: given a $k$-scheme $T$ of the form $T=\Spec(C)$ where $C$ is
a complete discrete valuation ring with algebraically closed residue
field $k_0$, and with fraction field $K$ say, and given a divisor $D$ on
$X_K$, its closure $D'$ is a divisor on $X_T$.

This statement follows from the Ramanujam--Samuel Theorem
\cite[21.14.1]{EGAIV4}, a result in commutative algebra.  Apply it
taking $B$ to be the local ring of a closed point of $X_T$, and $A$ to be the
local ring of the image point in $X_{k_0}$.  The hypotheses hold because
$A$ and $B$ share the residue field $k_0$.   The completion $\wh A$ is a
normal domain by a theorem of Zariski's  \cite[Thm.~32, p.~320]{ZS60}.
And $A\to B$ is formally smooth because $C$ is a formal power series
ring over $k_0$ by a theorem of Cohen's \cite[Thm.~32, p.~320]{ZS60}.

The third alternative proof is somewhat like the second, but involves
some geometry instead of the Ramanujam--Samuel Theorem; see
\cite[Thm.~19, p.~138]{AK74}.  Moreover, $C$ need not be complete, just
discrete, and $k_0$ need not be algebraically closed.  Here is the idea.
Let $E\subset X_{k_0}$ be the closed fiber of $D'/T$.  In
$\IHilb_{X_L/L}$, form the sets $U$ and $V$ parameterizing the
divisorial cycles linearly equivalent to those of the form $D_L+H$ and
$E_{L}+H$ as $H$ ranges over the divisors whose associated sheaves are
algebraically equivalent to $\mc O_{\!X_L}(n)$ for a suitably large $n$.

It can be shown that $U$ and $V$ are dense open subsets of the same
irreducible component of $\IHilb_{X_L/L}$.  Hence they have a common
point.  Let $I$ and $J$ be the ideals of $D$ and $E$.  Then there are
invertible sheaves $\mc L$ and $\mc M$ on $X_L$ such that $I_L\ox\mc L$
and $J_L\ox\mc M$ are isomorphic.  Since $I_L$ is invertible, so is
$J_L$. Hence so is $J$.  Thus $E$ is a divisor, as desired.

More generally, if $X$ is not necessarily projective, but is simply
complete and normal, then $\IPicz_{X/k}$ is still complete, whether $X$
is a scheme or algebraic space.  Indeed, the original proof and its
second alternative work without essential change.  The first and third
alternatives require $X$ to be projective.  However, it is easy to see
as follows
that this case implies the general case.

Namely, we may assume $k$ is algebraically closed.  By Chow's lemma,
there is a projective variety $Y$ and a birational map $\gamma\:Y\to X$.
Since $X$ and $Y$ are normal, a divisor $D$ on $X$ is the divisor of a
function $h$ if and only if $\gamma^*D$ is the divisor of $\gamma^*h$.
Hence the induced map $\IPicz_{X/k}\to\IPicz_{Y/k}$ is injective.  It
follows, as above, that $\IPicz_{X/k}$ is complete since $\IPicz_{Y/k}$
is.
 \end{rmk}

\begin{dfn}\label{dfn:algeq}
 Assume $S$ is the spectrum of a field $k$.  Let $\mc L$ and $\mc N$ be 
invertible sheaves on $X$. Then $\mc L$ is said to be {\it algebraically
equivalent to}  $\mc N$  if, for some $n$ and all $i$ with $1\le i\le n$,
there exist a connected $k$-scheme of finite type $T_i$, geometric
points $s_i,\,t_i$ of $T_i$ with the same field, and an invertible sheaf
$\mc M_i$ on $X_{T_i}$ such that
  $$\mc L_{s_1}\simeq\mc M_{1,s_1},\
    \mc M_{1,t_1}\simeq\mc M_{2,s_2},\ \dotsc,\
    \mc M_{n-1,t_{n-1}}\simeq\mc M_{n,s_n},\
    \mc M_{n,t_n}\simeq\mc N_{t_n}.$$
  \end{dfn}

\begin{prp}\label{prp:algeq}
 Assume $S$ is the spectrum of a field $k$.  Assume $\IPic_{X/k}$ exists
and represents $\Pic_{(X/k)\fppf}$.  Let $\mc L$ be an invertible sheaf
on $X$, and $\lambda\in\IPic_{X/k}$ the corresponding point.  Then $\mc
L$ is algebraically equivalent to $\mc O_{\!X}$ if and only if
$\lambda\in\IPicz_{X/k}$.
 \end{prp}
 \begin{proof}
 Suppose $\mc L$ is algebraically equivalent to $\mc O_{\!X}$, and use the
notation of Definition~\ref{dfn:algeq}.  Then $\mc M_i$ defines a map
$\tau_i\:T_i\to \IPic_{X/k}$.  Now, $\mc M_{n,t_n} \simeq\mc
O_{\!X_{t_n}}$.  So $\tau_n(t_n)\in\IPicz_{X/k}$.  Suppose
$\tau_i(t_i)\in\IPicz_{X/k}$.  Then $\tau_i(T_i)\subset \IPicz_{X/k}$
since $T_i$ is connected.  So $\tau_i(s_i)\in \IPicz_{X/k}$.  But $\mc
M_{i,s_i}\simeq\mc M_{i-1,t_{i-1}}$.  So
$\tau_{i-1}(t_{i-1})\in\IPicz_{X/k}$.  Descending induction yields
$\tau_1(s_1)\in \IPicz_{X/k}$.  But $\mc M_{1,s_1}\simeq\mc L_{s_1}$.
Thus $\lambda\in\IPicz_{X/k}$.

Conversely, suppose $\lambda\in\IPicz_{X/k}$.  The inclusion
$\IPicz_{X/k}\into \IPic_{X/k}$ is defined by an invertible sheaf $\mc
M$ on $X_T$ for some fppf covering $T\to \IPicz_{X/k}$.  Let $t_1,\,t_2$
be geometric points of $T$ lying over $\lambda,\,0\in\IPicz_{X/k}$.  Let
$T_1,\,T_2\subset T$ be irreducible components containing $t_1,\,t_2$,
and $T_1',\,T_2'\subset \IPicz_{X/k}$ their images.  The latter contain
open subsets because an fppf map is open by \cite[2.4.6]{EGAIV2}.  Since
$\IPicz_{X/k}$ is irreducible by Lemma~\ref{lem:agps} , these open
subsets contain a common point.  Say it is the image of geometric points
$t_1,\,s_2$ of $T_1,\,T_2$.  Then $\mc M_{1,t_1}\simeq\mc M_{2,s_2}$ by
Exercise~\ref{ex:gpts}.  Set $\mc M_i:=\mc M|T_i$.  Thus $\mc L$ is
algebraically equivalent to $\mc O_X$.
 \end{proof}

\begin{thm}\label{thm:tgtsp}
 Assume $S$ is the spectrum of a field $k$.  Assume $\IPic_{X/k}$ exists
and represents $\Pic_{(X/k)\et}$.  Let $\tu T_0\IPic_{X/k}$ denote the
tangent space at $0$.  Then
	$$\tu T_0\IPic_{X/k}=\tu H^1(\mc O_{\!X}).$$
 \end{thm}
 \begin{proof} Let $P$ be any $k$-scheme locally of finite type, $e\in
P$ a rational point, $A$ its local ring, and $\I m$ its maximal ideal.
Usually, by the ``tangent space'' $\tu T_eP$ is meant the Zariski
tangent space $\Hom(\I m/\I m^2, k)$.  However, as in differential
geometry, $\tu T_eP$ may be viewed as the vector space of
$k$-derivations $\delta\:A\to k$.  Indeed, $\delta(\I m^2)=0$; so
$\delta$ induces a linear map $\I m/\I m^2\to k$.  Conversely, every
such linear map arises from a $\delta$, and $\delta$ is unique because
$\delta(1)=0$.

Let $k_\ve$ be the ring of ``dual numbers,'' the ring obtained from $k$
by adjoining an element $\ve$ with $\ve^2=0$.  Then any (fixed)
derivation $\delta$ induces a local homomorphism of $k$-algebras
$u\:A\to k_\ve$ by $u(a):=\overline a+\delta(a)\ve$ where $\overline
a\in k$ is the residue of $a$.  Conversely, every such $u$ arises from a
unique $\delta$.

On the other hand, to give a $u$ is the same as to give a $k$-map
$t_\ve$ from the ``free tangent vector'' $\Spec(k_\ve)$ to $P$ such that
the image of $t_\ve$ has support at $e$.  Denote the set of $t_\ve$ by
$P(k_\ve)_e$.  Thus, as sets,
 \begin{equation}\label{equation:5.7a}
 T_eP=P(k_\ve)_e.
 \end{equation}
 
The vector space structure on $T_eP$ transfers as follows.  Given $a\in
k$, define a $k$-algebra homomorphism $\mu_a\:k_\ve\to k_\ve$ by
$\mu_a\ve:= a\ve$.  Now, let $\delta$ be a derivation, and $u$ the
corresponding homomorphism.  Then $a\delta$ corresponds to $\mu_au$.
Thus multiplication by $a$ transfers as $P(\mu_a)\:P(k_\ve)_e\to
P(k_\ve)_e$.

  Let $k_{\ve,\ve'}$ denote the ring obtained from $k_\ve$ by adjoining
an element $\ve'$ with $\ve\ve'=0$ and $(\ve')^2=0$.  Define a
homomorphism $\sigma_1\:k_{\ve,\ve'}\to k_\ve$ by $\ve\mapsto\ve$ and
$\ve'\mapsto0$.  Define another $\sigma_2\:k_{\ve,\ve'}\to k_{\ve}$ by
$\ve\mapsto0$ and $\ve'\mapsto\ve$.  The $\sigma_i$ induce a map of sets
	$$\pi\:P(k_{\ve,\ve'})_e\to P(k_\ve)_e\x P(k_\ve)_e$$
 where $P(k_{\ve,\ve'})_e$ is the set of maps with image supported at
$e$. Plainly $\pi$ is bijective.

Define a third homomorphism $\sigma\:k_{\ve,\ve'}\to k_\ve$ by
$\ve\mapsto\ve$ and $\ve'\mapsto\ve$.  Given two derivations $\delta,\,
\delta'$,  define a homomorphism $v\:A\to k_{\ve,\ve'}$ by $v(a):=\overline
a+\delta(a)\ve+\delta'(a)\ve'$.  Then  $\delta+\delta'$ corresponds to
$\sigma v$.  Therefore, addition on $T_eP$   transfers as
	$$\alpha\:P(k_\ve)_e\x P(k_\ve)_e\to P(k_\ve)_e
	\text{ where }\alpha:=P(\sigma)\pi^{-1}.$$

Suppose $P$ is a group scheme, $e\in P$ the identity.  The natural ring
homomorphism $\rho\:k_\ve\to k$ induces a group homomorphism
$P(\rho)\:P(k_\ve)\to P(k)$.  Plainly
 \begin{equation}\label{equation:5.7b}
 P(k_\ve)_e=\ker P(\rho).
 \end{equation}
 The left side $\tu T_0P$ is a vector space; the right side is a group.
Does addition on the left match multiplication on the right?  Yes,
indeed!  We know $\alpha$ is the addition map.  We must show $\alpha$ is
also the multiplication map.  Let us do so.

Since $\pi$ and $P(\sigma)$ arise from ring homomorphisms, both are
group homomorphisms.  Now $\alpha:=P(\sigma)\pi^{-1}$.  Hence
$\alpha$ is a group homomorphism too.  So
	$$\alpha(m,n)=	 \alpha(m,e)\cdot \alpha(e,n)$$
 where $e\in P(k_\ve)_e$ is the identity.  So we have to show
    $\alpha(m,e)=m\text{ and }\alpha(e,n)=n$.

Consider the inclusion $\iota\:k_\ve\to k_{\ve,\ve'}$.  Plainly
$\sigma_2\iota\:k_\ve\to k_{\ve'}$ factors through $\rho\:k_\ve\to k$.
Hence $P(\sigma_2)P(\iota)(m)=e$ for any $m\in P(k_\ve)_e$ owing to
Formula~(\ref{equation:5.7b}).  On the other hand, $\sigma_1\iota$ is the
identity of $k_\ve$.  Thus $\pi P(\iota)(m)=(m,e)$

Plainly $\sigma\iota\:k_\ve\to k_\ve$ is also the identity of $k_\ve$.
So $P(\sigma)P(\iota)(m)=m$.  Hence
	$$\alpha(m,e)=P(\sigma)\pi^{-1}\pi P(\iota)(m)=m.$$
 Similarly $\alpha(e,n)=n$.  Thus $\alpha$ is the multiplication map.

Take $P:=\IPic_{X/k}$, so $e=0$.  Then Formulas~(\ref{equation:5.7a})
and (\ref{equation:5.7b}) yield
 \begin{equation}\label{equation:5.7d}
  \tu T_0\IPic_{X/k}=\ker\bigl(\IPic_{X/k}(k_\ve)
  \to \IPic_{X/k}(k)\bigr).
 \end{equation}

To compute this kernel, set $X_\ve:=X\ox_kk_\ve$, and form the truncated
exponential sequence of sheaves of Abelian groups:
	$$0\to \mc O_{\!X} \to \mc O_{\!X_\ve}^* \to \mc O_{\!X}^* \to 1,$$
 where the first map takes a local section $b$ to $1+b\ve$.  This
sequence is split by the map $\mc O_{\!X}^*\to \mc O_{\!X_\ve}^*$ defined by
$a\mapsto a+0\cdot\ve$.  Hence taking cohomology yields this split exact
sequence of Abelian groups:
	$$0\to \tu H^1(\mc O_{\!X}) \to \tu H^1(\mc O_{\!X_\ve}^*)
	 \to \tu H^1(\mc O_{\!X}^*) \to 1.$$

However,  $\IPic_{X/k}$ represents  $\Pic_{(X/k)\et}$, which is the sheaf
associated to  the presheaf $T\mapsto \tu H^1(\mc O_{\!X_T}^*)$.  So 
there is a natural commutative square of groups
 \begin{equation}\label{CD:tgtsp}\begin{CD}
 \tu H^1(\mc O_{\!X_\ve}^*) @>>> \tu H^1(\mc O_{\!X}^*)\\
               @VVV	               @VVV\\
 \IPic_{X/k}(k_\ve)       @>>> \IPic_{X/k}(k).
 \end{CD} \end{equation}\vskip-\smallskipamount
 Hence, there is an induced homomorphism between the horizontal kernels.
Owing to Formula~(\ref{equation:5.7d}), this homomorphism is an additive
map
	$$v\:\tu H^1(\mc O_{\!X}) \to \tu T_0\IPic_{X/k}.$$
 
Let $a\in k$.  On $T_0\IPic_{X/k}$, scalar multiplication by $a$ is,
owing to the discussion after Formula~(\ref{equation:5.7a}), the map
induced by $\mu_a\:k_\ve\to k_\ve$.  Now, $\mu_a$ induces an
endomorphism of the above square.  At the top, it arises from the map of
sheaves of groups $\mc O_{\!X}^* \to \mc O_{\!X}^*$ defined by $\ve\mapsto
a\ve$.  So the induced endomorphism of $\tu H^1(\mc O_{\!X})$ is scalar
multiplication by $a$.  Thus $v$ is a map of $k$-vector spaces.

 Square (\ref{CD:tgtsp}) maps to the corresponding square obtained by
making a field extension $K/k$.  Since the kernels are vector spaces,
there is an induced square
 $$\begin{CD}
 \tu H^1(\mc O_{\!X})\ox_k K    @>>> \tu H^1(\mc O_{\!X_K})\\
   @V v\ox_k K VV                             @VVV\\
 \tu T_0\IPic_{X/k}\ox_k K @>>> \tu T_0\IPic_{X_K/K}.
  \end{CD}$$\vskip-\smallskipamount
  The two horizontal maps are isomorphisms.  Hence, if the right-hand
map is an isomorphism, so is $v$.  Thus we may assume $k$ is
algebraically closed.

In Square (\ref{CD:tgtsp}), the two vertical maps are isomorphisms by
Exercise~\ref{ex:Alr} since $\IPic_{X/k}$ represents $\Pic_{(X/k)\et}$ and
$k$ is algebraically closed.  Therefore, $v$
too is an isomorphism, as desired.
 \end{proof}

\begin{rmk}\label{rmk:tgtsp}
  There is a relative version of Theorem~\ref{thm:tgtsp}.  Namely,
assume $\IPic_{X/S}$ exists, represents $\Pic_{(X/S)\et}$, and is
locally of finite type, but let $S$ be arbitrary.  Then $\tu R^1f_*{\mc
O}_X$ is equal to the normal sheaf of $\IPic_{X/S}$ along the identity
section, its ``Lie algebra''; the latter is simply the dual of the
restriction to this section of the sheaf of relative differentials.  For
more information, see the recent treatment \cite[\S\,1]{LLR} and the
references it cites.
 \end{rmk}

\begin{cor}\label{cor:sm} Assume $S$ is the spectrum of a field
$k$. Assume $\IPic_{X/k}$ exists and represents $\Pic_{(X/k)\et}$.  Then
	$$\dim \IPic_{X/k}\le \dim\tu H^1(\mc O_{\!X}).$$
 Equality holds if and only if $\IPic_{X/k}$ is smooth at $0;$ if so, then
$\IPic_{X/k}$ is smooth of dimension $\dim\tu H^1(\mc O_{\!X})$
everywhere.
 \end{cor}
 \begin{proof} Plainly we may assume $k$ is algebraically closed. Then,
given any closed point $\lambda\in\IPic_{X/k}$, there is an automorphism
of $\IPic_{X/k}$ that carries 0 to $\lambda$, namely, ``multiplication''
by $\lambda$.  So $\IPic_{X/k}$ has the same dimension at $\lambda$ as
at 0, and $\IPic_{X/k}$ is smooth at $\lambda$ if and only if it is
smooth at 0.

By general principles, $\dim_0 \IPic_{X/k}\le \dim\tu T_0\IPic_{X/k}$,
and equality holds if and only if $\IPic_{X/k}$ is regular at $0$.
Moreover, $\IPic_{X/k}$ is regular at $0$ if and only if it is smooth
at $0$ since $k$ is algebraically closed.  Therefore, the corollary
results from Theorem~\ref{thm:tgtsp}. 
 \end{proof}

\begin{cor}\label{cor:ch0}
 Assume $S$ is the spectrum of a field $k$. Assume $\IPic_{X/k}$ exists
and represents $\Pic_{(X/k)\et}$. If $k$ is of characteristic $0$, then
$\IPic_{X/k}$ is smooth of dimension $\dim\tu H^1(\mc O_{\!X})$ everywhere.
 \end{cor}
 \begin{proof} Since $k$ is of characteristic $0$, any group scheme
locally of finite type over $k$ is smooth by Cartier's theorem
\cite[Thm.~1, p.~167]{Mm66}.  So the assertion follows from
Corollary~\ref{cor:sm}.
 \end{proof}

\begin{rmk}\label{rmk:Igusa}
  Over a field $k$ of positive characteristic, $\IPic_{X/k}$ need not be
smooth, even when $X$ is a connected smooth projective surface.
Examples were constructed by Igusa \cite{Ig55} and Serre
\cite[n$^\circ$~20]{Sr56}.

On the other hand, Mumford \cite[Lect.~27, pp.~193--198]{Mm66} proved that
$\IPic_{X/k}$ is smooth if and only if all of Serre's Bockstein
operations $\beta_i$ vanish; here
	$$\beta_1:\tu H^1(\mc O_{\!X})\to \tu H^2(\mc O_{\!X})
 \text{ and }\beta_i\:\ker\beta_{i-1}\to\cok\beta_{i-1}\text{ for }i\ge2.$$
 In fact, the tangent space to $\IPic_{X_{\red}/k}$ is the subspace of $\tu
H^1(\mc O_{\!X})$ given by
	$$\tu T_0\IPic_{X_{\red}/k}=\bigcap\ker\beta_i.$$
 Moreover, here $X$ need not be smooth or 2-dimensional.

The examples illustrate further pathologies.  Set
	$$g:=\dim \IPic_{X/k},\ h^{0,1}:=\dim \tu H^1(\mc O_{\!X}),
	\text{ and }h^{1,0}:=\dim \tu H^0(\Omega^1_{\!X}).$$
 In Igusa's example, $g=1$, $h^{0,1}=2$, and $h^{1,0}=2$; in Serre's,
$g=0$, $h^{0,1}=1$, and $h^{1,0}=0$.  Moreover, Igusa had just proved
that, in any event, $g\le h^{1,0}$.

By contrast, in characteristic 0, Serre's Comparison Theorem
\cite[Thm.~2.1, p.~440]{Ha83} says that $h^{0,1}$ and $h^{1,0}$ can be
computed by viewing $X$ as a complex analytic manifold.  Hence Hodge
Theory yields
	$$h^{0,1}=h^{1,0}\text{ and }h^{0,1}+h^{1,0}=b$$
 where $b$ is the first Betti number; see \cite[p.~200]{Za35}.
Therefore, the following exercise now yields the Fundamental  Theorem of
Irregular Surfaces (\ref{equation:FTIS}).
 \end{rmk}

\begin{ex}\label{ex:q=0}
 Assume $S$ is the spectrum of a field $k$.  Assume $X$ is a projective,
smooth, and geometrically irreducible surface.  According to the
original definitions as stated in modern terms, the ``geometric genus''
of $X$ is the number $p_g:=\dim\tu H^0(\Omega^2_{\!X})$; its
``arithmetic genus'' is the number $p_a:=\phi(0)-1$ where $\phi(n)$ is
the polynomial such that $\phi(n)=\dim\tu H^0(\Omega^2_{\!X}(n))$ for
$n\gg0$; and its ``irregularity'' $q$ is the difference between the two
genera, $q:=p_g-p_a$.

Show $\dim \IPic_{X/k}\le q$, with equality in characteristic 0.
 \end{ex}

\begin{ex}\label{ex:Enriques}
 Assume $S$ is the spectrum of an algebraically closed field $k$.
Assume $X$ is projective and integral.  Set $q:=\dim\tu H^1(\mc O_X)$.

Show $q=0$ if and only if every algebraic system of curves is
``contained completely in a linear system.''  The latter condition means
just that, given any relative effective divisor $D$ on $X_T/T$ where $T$
is a connected $k$-scheme, there exist invertible sheaves $\mc L$ on $X$
and $\mc N$ on $T$ such that $\mc O_{\!X_T}(D)\simeq \mc L_T\ox f_T^*\mc
N$.  The condition may be put more geometrically: in the notation of
Exercise~\ref{ex:LinSys}, it means there is a map, necessarily unique,
$w\:T\to L$ such that $(1\x w)^{-1}W=D$.

In characteristic 0, show  $q=0$ if the condition
holds for all smooth such $T$.
 \end{ex}

\begin{rmk}\label{rmk:charsys}
 Assume $S$ is the spectrum of a field $k$.  Assume $X/k$ is projective
and geometrically integral. Let $D\subset X$ be an effective
divisor, and $\mc N_D$ its normal sheaf.  Let $\delta\in\IDiv_{X/k}$ be
the point representing $D$, and $\lambda\in \IPic_{X/k}$ the point
representing $\mc O_{\!X}(D)$.

Then the tangent space at $\delta$ is given by the formula
\begin{equation}\label{equation:ccs}
 \tu T_\delta\IDiv_{X/k}=\tu H^0(\mc N_D),
 \end{equation}
 which respects the vector space structure of each side.  This formula
can be proved with a simple elementary computation; see \cite[Cor.,
p.~154]{Mm66}.

Form the fundamental exact sequence of sheaves
	$$0\to \mc O_{\!X} \to  \mc O_{\!X}(D) \to \mc N_D\to 0,$$
 and consider its associated long exact sequence of cohomology groups
$$\begin{CD}
 0 @>>>\tu H^0(\mc O_{\!X})@>>>\tu H^0(\mc O_{\!X}(D))@>>>\tu H^0(\mc N_D)
\\@()\\
\hphantom0@>\partial^0>>\tu H^1(\mc O_{\!X})@>>>
        \tu H^1(\mc O_{\!X}(D))@>u>>\tu H^1(\mc N_D)
@>\partial^1>>\tu H^2(\mc O_{\!X})
\end{CD}
$$
 Another elementary computation shows that the boundary map $\partial^0$
is equal to the tangent map of the Abel map, $\tu T_\delta\IDiv_{X/k}\to
T_\lambda\IPic_{X/k}$; see \cite[Prp., p.~165]{Mm66}.

By definition, $D$ is said to be ``semiregular'' if the boundary map
$\partial^1$ is injective.  Plainly, it is equivalent that $u=0$.  So it
is equivalent that $\dim\tu H^0(\mc N_D)=R$ where
 $$R:=\dim\tu H^1(\mc O_{\!X})+\dim\tu H^0(\mc O_{\!X}(D))
	-1-\dim\tu H^1(\mc O_{\!X}(D)).$$
 Semiregularity was recognized in 1944 by Severi as precisely the right
positivity condition for the old (1904) theorem of the completeness of
the characteristic system, although Severi formulated the condition in
an equivalent dual manner.

In its modern formulation, the theorem of the completeness of the
characteristic system asserts that $\IDiv_{X/k}$ is smooth of dimension
$R$ at $\delta$ if and only if $D$ is semiregular, provided the
characteristic is 0, or more generally, $\IPic_{X/k}$ is smooth.
Indeed, Formula~(\ref{equation:ccs}) says that the characteristic
system of $\IDiv_{X/k}$ on $D$ is always equal to the complete linear
system of the invertible sheaf $\mc N_D$.  But $\IDiv_{X/k}$ can have a
nilpotent or a singularity at $\delta$.  So, in effect, Enriques and
Severi had simply sought conditions guaranteeing $\IDiv_{X/k}$ is smooth
of dimension $R$ at $\delta$.

The first purely algebraic discussion of the theorem was made by
Grothendieck \cite[Sects.~221-5.4 to 5.6]{FGA}, and he proved it in the
two most important cases.  Specifically, he noted that, if $\tu H^1(\mc
O_{\!X}(D))=0$, then the Abel map is smooth; see the end of the proof of
Theorem~\ref{th:main}.  Hence $\IDiv_{X/k}$ is smooth at $\delta$
if and only if $\IPic_{X/k}$ is smooth; furthermore, $\IPic_{X/k}$ is
smooth in characteristic zero by Cartier's theorem.  Grothendieck
pointed out that this case had been treated with transcendental means by
Kodaira in 1956.

Grothendieck also observed $\tu H^1(\mc N_D)$ houses the obstruction to
deforming $D$ in $X$.  Hence, if this group vanishes, then $\IDiv_{X/k}$
is smooth at $\delta$ in any characteristic.  Mumford
\cite[pp.~157--159]{Mm66} explicitly worked out the obstruction and its
image under $\partial^1$.  If $\partial^1=0$, so if $D$ is semiregular,
then $\IDiv_{X/k}$ is smooth if this image vanishes.  Inspired by work
of Kodaira and Spencer in 1959, Mumford used the exponential in
characteristic 0 and proved the image vanishes.  Mumford did not use
Cartier's theorem; so the latter results from taking a $D$ with $\tu
H^1(\mc O_{\!X}(D))=0$.

A purely algebraic proof of the full completeness theorem is given in
\cite[Thm., p.~307]{Kl73}.  This proof was inspired by Kempf's
(unpublished) thesis.  The proof does not use obstruction theory, but
only simple formal properties of a scheme of the form $\I P(\mc Q)$
where $\mc Q$ arises from an invertible sheaf $\mc F$ as in
Subsection~\ref{sb:Q}.  This proof works, more generally, if $S$ is
arbitrary and if $X/S$ is projective and flat and has integral geometric
fibers.  Here $D$ is the divisor on the geometric fiber through
$\delta$.  Then provided $\IPic_{X/S}$ is smooth, $\IDiv_{X/S}$ is
smooth of relative dimension $R$ at $\delta$ if and only if $D$ is
semiregular.

There is a celebrated example, valid over an algebraically closed field $k$ 
of any characteristic, where $\IDiv_{X/k}$ is nonreduced at $\delta$.
The example was discovered by Severi and Zappa in the 1940s, and is
explained in \cite[pp.~155--156]{Mm66}.  Here is the idea.  Let $C$ be
an elliptic curve, and $0\to \mc O_{\!C}\to\mc E\to \mc O_{\!C}\to0$ the
nontrivial extension; set $X:=\I P(\mc E)$.

Let $D$ be the section of $X/C$ defined by $\mc E\to \mc O_{\!C}$.  Then
$\mc N_D= \mc O_{\!C}$ by \cite[Prp.~2.8, p.~372]{Ha83}; so $\dim
T_\delta\IDiv_{X/k}=1$.  However, $\delta$ is an isolated point.
Otherwise, the connected component of $\delta$ contains a second closed
point.  And it represents a curve $D'$ algebraically equivalent to $D$.
So $\deg\mc O_{\!D}(D')=\deg\mc N_D=0$.  Hence $D$ and $D'$ are
disjoint.  Let $F$ be a fiber.  Then $\deg\mc O_{\!F}(D)=\deg\mc
O_{\!F}(D')=1$.  Hence $D'$ is a second section.  Therefore, $\mc E$ is
decomposable, a contradiction.
 \end{rmk}

\begin{prp}\label{prp:H2}
 Assume $\IPic_{X/S}$ exists and represents $\Pic_{(X/S)\et}$.  Let
$s\in S$ be a point such that $\tu H^2(\mc O_{\!X_s})=0$.  Then there
exists an open neighborhood of $s$ over which $\IPic_{X/S}$ is smooth.
 \end{prp}
 \begin{proof}
 By the Semicontinuity Theorem \cite[7.7.5-I)]{EGAIII2}, there exists an
open neighborhood $U$ of $s$ such that $\tu H^2(\mc
O_{\!X_t})=0$ for all $t\in U$.  Replace $S$ by $U$.

By \cite[{\bf0}-10.3.1, p.~20]{EGAIII1}, there is a flat local homomorphism
from $\mc O_s$ into a Noetherian local ring $B$ whose residue field is
algebraically closed.  By \cite[6.8.3]{EGAIV2}, if $f_B$ is smooth, then
$f\:X\to S$ is smooth along $X_s$.  Replace $S$ by $\Spec(B)$.

By the Infinitesimal Criterion for Smoothness \cite[p.~67]{SGA1}, it
suffices to show this: given any $S$-scheme $T$ of the form $T=\Spec(A)$
where $A$ is an Artin local ring that is a finite $\mc O_s$-algebra and
given any closed subscheme $R\subset T$ whose ideal $\mc I$ has square
0, every $R$-point of $\IPic_{X/S}$ lifts to a $T$-point.

The residue field of $A$ is a finite extension of $k_s$, which is
algebraically closed; so the two fields are equal. Hence the $R$-point
is defined by an invertible sheaf on $X_{R}$ by Exercise~\ref{ex:Alr}.  So
we want to show invertible sheaves on $X_{R}$ lift to $X_T$.

Since $\mc I^2=0$, we can form the truncated exponential sequence
  $$0\to f_T^*{\mc I} \to \mc O_{\!X_T}^* \to \mc O_{\!X_R}^* \to1.$$
  It yields the following exact sequence:
	$$\tu H^1(\mc O_{\!X_T}^*)\to \tu H^1(\mc O_{\!X_R}^*)
	\to \tu H^2(f_T^*{\mc I}).$$
   Hence it suffices to show  $\tu H^2(f_T^*{\mc I})=0$.

Since $T$ is affine, $\tu H^2(f_T^*{\mc I})=H^0(\tu R^2f_{T*}f_T^*{\mc
I})$ owing to \cite[Prp.~8.5, p.~251]{Ha83}.  But $\tu H^2(\mc
O_{\!X_t})=0$ for all $t\in S$.  Hence $\tu R^2(f_T^*{\mc I})=0$ by the
Property of Exchange \cite[7.7.5 II and 7.7.10]{EGAIII2}.  Thus $\tu
H^2(f_T^*{\mc I})=0$, as desired.
 \end{proof}

\begin{prp}\label{prp:P0}
 Assume $\IPic_{X/S}$ exists, and represents $\Pic_{(X/S)\fppf}$.  For $s\in
S$, assume all the $\IPicz_{X_s/k_s}$ are smooth of the same
dimension.  Then $\IPic_{X/S}$ has an open group subscheme
$\IPicz_{X/S}$ of finite type whose fibers are the $\IPicz_{X_s/k_s}$.
Furthermore, if $S$ is reduced, then $\IPicz_{X/S}$ is smooth over $S$.
Moreover, if all the $\IPicz_{X_s/k_s}$ are complete and if
$\IPic_{X/S}$ is separated over $S$, then $\IPicz_{X/S}$ is closed in
$\IPic_{X/S}$ and proper over $S$.
 \end{prp}
 \begin{proof} First off,  $\IPic_{X/S}$ is locally of finite type by
 Proposition~\ref{prp:lft}.  Now, 
 for every $s\in S$, the schemes $\IPic_{X/S}\ox k_s$ and
$\IPic_{X_s/k_s}$ coincide by Exercise~\ref{ex:bschg}.  And the
$\IPicz_{X_s/k_s}$ are smooth of the same dimension by hypothesis.  So
the $\IPicz_{X_s/k_s}$ form an open subscheme $\IPicz_{X/S}$ of
$\IPic_{X/S}$ and the structure map $\sigma\:\IPicz_{X/S}\to S$ is
universally open by \cite[15.6.3 and 15.6.4]{EGAIV3}.  Furthermore, if
$S$ is reduced, then $\sigma$ is flat by \cite[15.6.7]{EGAIV3}, so
smooth by \cite[17.5.1]{EGAIV4}.

Define a map $\alpha\:\IPicz_{X/S}\x_S \IPicz_{X/S}\to\IPic_{X/S}$ by
$\alpha(g,h):=gh^{-1}$.  Then $\alpha$ factors through the open
subscheme $\IPicz_{X/S}$ because forming $\alpha$ commutes with passing
to the fibers.  Hence $\IPicz_{X/S}$ is a subgroup

To prove $\IPicz_{X/S}$ is of finite type, we may work locally on $S$,
and so assume $S$ is Noetherian.  Since $\IPic_{X/S}$ is
locally of finite type, we need only prove  $\IPicz_{X/S}$ is
quasi-compact.

Let $V\subset S$ be a nonempty affine open subscheme, and
$U\subset\sigma^{-1}V$ another.  Then $\sigma U\subset S$ is open since
$\sigma$ is an open map.  Set $U':=\sigma^{-1}\sigma U$.  Then $U'$ is
open, and $\alpha$ restricts to a map $\alpha'\:U\x_VU\to U'$.  In fact,
$\alpha'$ is surjective because its geometric fibers are surjective by
an argument at the end of the proof of Lemma~\ref{lem:agps}.  Now,
$U\x_VU$ so quasi-compact.  Hence $U'$ is quasi-compact.

Set $T:=S-\sigma U$.  By Noetherian induction, we may assume
$\sigma^{-1}T$ is quasi-compact.  Therefore, $U'\cup\sigma^{-1}T$ is
quasi-compact.  But it is equal to $\IPicz_{X/S}$.  Thus $\IPicz_{X/S}$
is quasi-compact, as desired.

Moreover, if all the $\IPicz_{X_s/k_s}$ are complete and if
$\IPic_{X/S}$ is separated over $S$, then $\IPicz_{X/S}$ is proper over
$S$ by \cite[15.7.11]{EGAIV3}.  Finally, consider the inclusion map
$\IPicz_{X/S} \into \IPic_{X/S}$.  It is proper as $\IPicz_{X/S}$ is
proper and $\IPic_{X/S}$ is separated; hence, $\IPicz_{X/S}$ is closed.
 \end{proof}

\begin{rmk}\label{rmk:Vistoli}
 Assume the characteristic is 0 and $f\:X\to S$ is smooth and proper.
Then, in the last two propositions, more can be said, as Vistoli
explained to the author in early May 2004.  These additions result from
Part~(i) of Theorem (5.5) on p.~123 in Deligne's article \cite{De68}
(which uses Hodge theory on p.~121).  Deligne's theorem asserts that,
under the present conditions, all the sheaves $\tu R^qf_*\Omega^p_{X/S}$
are locally free of finite rank, and forming them commutes with changing
the base.

In Proposition~\ref{prp:H2}, {\it $\IPic_{X/S}$ is smooth over the
connected component $S_0$ of $s$.}  Indeed, Deligne's theorem implies
$\tu R^2f_*\mc O_{X/S}\big|S_0$ vanishes since it is locally free, so of
constant rank, and its formation commutes with passage to every fiber,
in particular, to that over $s\in S_0$.  Hence, $\tu H^2(\mc
O_{\!X_t})=0$ for every $t\in S_0$.  Therefore,
Proposition~\ref{prp:H2}, as it stands, implies $\IPic_{X/S}$ is smooth
over $S_0$.

In Proposition~\ref{prp:P0}, {\it there is no need to assume all the
$\IPicz_{X_s/k_s}$ are smooth of the same dimension.}  Indeed, all the
$\IPicz_{X_s/k_s}$ are smooth of dimension $\smash{\dim\tu H^1(\mc
O_{\!X_s})}$ by Corollary~\ref{cor:sm}.  But $\smash{\dim\tu H^1(\mc
O_{\!X_s})}$ is constant on each connected component of $S$ owing to
Deligne's theorem.

{\it Furthermore, if $f\:X\to S$ is smooth and projective Zariski
locally over $S$, then $\IPicz_{X/S}$ is smooth whether or not $S$ is
reduced.}  Indeed, we may assume $S$ is of finite type over $\bb C$; the
reduction is standard, and sketched by Deligne at the beginning of his
proof of his Theorem (5.5).  By the Infinitesimal Criterion for
Smoothness \cite[p.~67]{SGA1}, it suffices to show this: given any
$S$-scheme $T$ of the form $T=\Spec(A)$ where $A$ is an Artin local ring
that is a finite $\bb C$-algebra and given any closed subscheme
$R\subset T$, every $R$-point of $\IPicz_{X/S}$ lifts to a $T$-point.
But, $\IPicz_{X/S}$ is open in $\IPic_{X/S}$ by the above argument.
Hence it suffices to show every $R$-point of $\IPicz_{X/S}$ lifts to a
$T$-point of $\IPic_{X/S}$.

The residue field of $A$ is a finite extension of $\bb C$; so the two
fields are equal.  Hence the $R$-point is defined by an invertible sheaf
on $X_{R}$ by Exercise~\ref{ex:Alr}.  So we want to show every
invertible sheaf $\mc L$ on $X_{R}$ lifts to $X_T$.  By Serre's
Comparison Theorem \cite[Thm.~2.1, p.~440]{Ha83}, it suffices to lift
$\mc L$ to an analytic invertible sheaf since $f_T\:X_T\to T$ is
projective.  So pass now to the analytic category.

Let $X_0$ denote the closed fiber of $X_T$, and form the exponential
sequence:
	$$0\to \bb Z_{X_0}\to \mc O_{X_0}\to \mc O_{X_0}^* \to 0.$$
Consider the class of $\mc L_{X_0}$ in $\tu H^1(\mc O_{X_0}^*)$.  It
maps to 0 in $\tu H^2(\bb Z_{X_0})$ because this group is discrete and
$\mc L$ defines an $R$-point of $\IPicz_{X/S}$.

Form the exponential map $\mc O_{X_R}\to \mc O_{X_R}^*$, form its kernel
$\mc Z$, and form the natural map $\kappa\:\mc Z\to Z_{X_0}$.  Then
$\kappa$ is bijective.  Indeed, let $a$ be a local section of $\mc Z$.
Then $1+a+a^2/2+\dotsb =1$.  Set $u:=1+a/2+\dotsb$.  Then $au=0$.
Suppose $a$ maps to the local section 0 of $\mc O_{X_0}$.  Then $a$ is
nilpotent.  So $u$ is invertible.  Hence $a=0$.  Thus $\kappa$ is
injective.  But $\kappa$ is obviously surjective.  Thus $\kappa$ is
bijective.

Consider the class $\lambda$ of $\mc L$ in $\tu H^1(\mc O_{X_R}^*)$.  It
follows that $\lambda$ maps to 0 in $\tu H^2(\mc Z)$.  Hence $\lambda$
comes from a class $\gamma$ in $\tu H^1(\mc O_{X_R})$.  By Deligne's
theorem, $\gamma$ lifts to a class $\gamma'$ in $\tu H^1(\mc O_{X_T})$.
The image of $\gamma'$ in $\tu H^1(\mc O_{X_T}^*)$ gives the desired
lifting of $\mc L$ to $X_T$.


 \end{rmk}

\begin{eg}\label{eg:Vistoli}
 The following example complements Proposition~\ref{prp:P0}, and was
provided, in early May 2004, by Vistoli.  The example shows
$\IPicz_{X/S}$ {\it can be smooth and proper over $S$ and open and
closed in $\IPic_{X/S}$, although $\IPic_{X/S}$ isn't smooth.}  In fact,
$f\:X\to S$ is smooth and projective, its geometric fibers are integral,
and $S$ is a smooth curve over an algebraically closed field $k$ of
arbitrary characteristic.  Furthermore, $\IPic_{X/S}$ has a component
that is a reduced $k$-point; so it is not smooth over $S$, nor even
flat.

To construct $f\:X\to S$, set $P:=\bb P^3_k$ and fix $d\ge4$.  Set $\mc
Q:=\tu H^0(\mc O_P(d))^*$ where the `$*$' means dual. Set $H:=\bb P(\mc
Q)$.  Then $H$ represents the functor (on $k$-schemes $T$) whose
$T$-points are the $T$-flat closed subschemes of $P_T$ whose fibers are
surfaces of degree $d$ by Exercise~\ref{ex:Enriques} and
Theorem~\ref{th:LinSys}.  Let $W\subset P\x H$ be the universal
subscheme.  Its ideal $\mc O_{P\x H}(-W)$ is equal to the tensor product
of the pullbacks of $\mc O_P(-d)$ and $\mc O_H(-1)$ by
Exercise~\ref{ex:LinSys}.  Set $N:=\dim H$.

Let $G$ be the Grassmannian of lines in $P$, and $L\subset P\x G$  the
universal line.  Let $\pi \:L\to G$ be the projection.  Form the exact
sequence of locally free sheaves
	$$0\to\mc K\to \pi _*\mc O_{P\x G}(d)\to \pi _*\mc O_L(d)\to 0,$$
 which defines $\mc K$; the right-hand map is surjective, since forming
it commutes with passing to the fibers, and on the fibers, it is plainly
surjective.  Set $I:=\I P(\mc K^*)$.  Then $I$ is smooth, irreducible,
and of dimension $4+N -(d+1)$, or $N-d+3$.

Note that $\pi _*\mc O_{P\x G}(d)=\mc Q^*_G$.  So there is a surjection
$\mc Q_G\to \mc K^*$, and it induces a closed embedding $I\subset G\x
H$.  Furthermore, given a $T$-point of $G\x H$, it lies in $I$ if and
only if $L_T\subset W_T$.  Indeed, the latter means the ideal of
$W_T\subset P_T$ maps to 0 in $\mc O_{L_T}$; in other words, the
composition
    $$\mc O_{P_T}(-d)\ox\mc O_H(-1)_{P_T}\to O_{P_T}\to\mc O_{L_T}$$
 vanishes.  Equivalently, the composition
     $$\mc O_H(-1)_T\to g_{T*}O_{P_T}(d)\to g_{T*}\mc O_{L_T}(d)$$
 vanishes.  But the first map is equal to the natural map $\mc
O_H(-1)_T\to \mc Q^*_T$.  So it is equivalent that this map factors
through $\mc K_T$, or that $\mc Q_T\to \mc O_H(1)_T$ factors through
$\mc K^*_T$.  Equivalently, the $T$-point of $G\x H$ lies in $I$.  Thus
$I$ is the graph of the incidence correspondence.

Let $J\subset H$ be the image of $I$.  Note $J\ne H$ since $\dim I=\dim
H-d+3$ and $d\ge4$.  Consider the open set $U\subset H$ over which $W\to
H$ is smooth. Then $U\cap J$ is nonempty.  Indeed, choose coordinates
$w,x,y,z$ for $P$.  Then $U\cap J$ contains the point representing the
surface $\{w^d-x^d=y^d-z^d\}$ if $p\nmid d$ or the surface
$\{wx^{d-1}=y^{d-1}\}$ if $p\mid d$, because, in either case, the
surface is smooth and  contains a line, either $\{w=x,\ y=z\}$ or 
$\{w=0,\ y=0\}$.

Let $s\in U\cap J$ be a simple $k$-point.  Take a line $S\subset H$
through $s$ and transverse to $J$.  Replace $S\/$ by $S\cap U$.  Let
$X\subset W$ be the preimage of $S$, and $f\:X\to S$ the induced map.
Then $f$ is smooth and projective, and its geometric fibers are
integral.  Hence $\IPic_{X/S}$ exists, is separated, and represents
$\Pic_{(X/S)\et}$ by Theorem~\ref{th:main}.  Moreover, $\tu H^1(\mc
O_{\!X_t})=0$ for each $t\in S$; hence, Corollary~\ref{cor:sm} implies
$\IPic_{X_t/k_t}$ is reduced and discrete.  In particular,
$\IPicz_{X_t/k_t}$ is smooth, of dimension 0, and complete.  Therefore,
Proposition~\ref{prp:P0} implies $\IPicz_{X/S}$ is smooth and proper
over $S$ and open and closed in $\IPic_{X/S}$

Since $s\in J$, the surface $X_s\subset P$ contains a line $M$.  And
$X_s$ is smooth since $s\in U$; hence $M$ is a divisor.  Say
$\mc O_{\!X_s}(M)$ defines the $k$-point $\mu\in\IPic_{X/S}$.  View $\mu$ as
a reduced closed subscheme.  Then $\mu$ is a connected component of its
fiber $\IPic_{X_s/k}$ because, as just noted, this fiber is reduced and
discrete.  It remains to prove $\mu$ is a connected component of
$\IPic_{X/S}$.

Suppose not, and let's find a contradiction.  Let $Q$ be the connected
component of $\mu\in\IPic_{X/S}$.  Let $k_\ve$ be the ring of ``dual
numbers,'' and set $T:=\Spec(k_\ve)$.  Then there is a closed embedding
$T\into \IPic_{X/S}$ supported at $\mu$.  However, $T$ does not embed
into the fiber $\IPic_{X_s/k}$ because the latter is reduced and
discrete.  Hence the structure map $\IPic_{X/S}\to S$ embeds $T$ into
$S$.

The embedding $T\into \IPic_{X/S}$ corresponds to an invertible sheaf
$\mc M$ on $X_T$ by Exercise~\ref{ex:Alr} since $k$ is algebraically
closed.  Moreover, $\mc M|X_s\simeq \mc O_{\!X_s}(M)$ since the embedding
is supported at $\mu$.

Note $\tu H^1(\mc O_{\!X_s}(M))=0$.  Indeed, by Serre duality
\cite[Cor.~7.7, p.~244, and Cor.~7.12, p.~246]{Ha83}, it suffices to
show $\tu H^1(\bigomega_{\!X_s}(-M))=0$.  Form the sequence
	$$0\to\bigomega_{\!X_s}(-M)\to\mc \bigomega_{\!X_s}
	 \to\bigomega_{\!X_s}|M \to0.$$
 Now, $\bigomega_{X_s}\simeq\mc O_{X_s}(d-4)$ since
$\bigomega_{X_s}\simeq\bigomega_P\ox \mc O_{X_s}(X_s)$ by
\cite[Prp.~8.20, p.~182]{Ha83} and $\bigomega_P\simeq\mc O_P(-4)$ by
\cite[Eg.~8.20.1, p.~182]{Ha83}.  But $\tu H^1(\mc O_{\!X_s}(d-4))=0$
because $X_s\subset P$ is a hypersurface, and $\tu H^0(\mc
O_{\!X_s}(d-4))\to \tu H^0(\mc O_M(d-4))$ is surjective because $\tu
H^0(\mc O_P(d-4))\to \tu H^0(\mc O_M(d-4))$ is.  Thus $\tu H^1(\mc
O_{\!X_s}(M))=0$.

Therefore, $\uH^0(\mc M)\ox k\to\uH^0(\mc O_{\!X_s}(M))$ is surjective
by the implication (v)$\Rightarrow$(iv) of Subsection~\ref{sb:Q}.  So
the section of $\mc O_{\!X_s}(M)$ defining $M$ extends to a section of
$\mc M$.  The extension defines a relative effective divisor on $X_T$,
which restricts to $M$, owing to the implication (iii)$\Rightarrow$(i)
of Lemma~\ref{lm:ctn}.  It follows that the embedding $T\into H$ factors
through $I$, so through $J$.  However, $J$ and $S$ meet transversally at
$s$; whence, $T$ cannot embed into $J\cap S$.  Here is the desired
contradiction.  The discussion is now complete.
 \end{eg}

\begin{ex}\label{ex:jac}
 Assume $X/S$ is projective and flat, its geometric fibers are integral
curves of arithmetic genus $p_a$, and $S$ is Noetherian.  Show the
``generalized Jacobians'' $\IPicz_{X_s/k_s}$ form a smooth
quasi-projective family of relative dimension $p_a$.  And show this
family is projective if and only if $X/S$ is smooth.
 \end{ex}

{\def\wh#1{#1^*}
\begin{rmk}\label{rmk:Ablsch}
 Assume $X$ is an Abelian $S$-scheme of relative dimension $g$; that is,
$X$ is a smooth and proper $S$-group scheme with geometrically connected
fibers of dimension $g$.  Then $X$ needn't be projective Zariski locally
over $S$.

  Indeed, according to Raynaud \cite[Rem.~8.c, p.~1315]{Ra66},
Grothendieck found two such examples: one where $S$ is reduced and
1-dimensional, and another where $S$ is the spectrum of the ring of dual
numbers of a field of characteristic 0.  In \cite[Ch.~XII]{Ra70},
Raynaud gave detailed constructions of similar examples.

Assume $X/S$ is projective, and $S$ is Noetherian.  Then $\IPicz_{X/S}$
exists, and is also a projective Abelian $S$-scheme of relative
dimension $g$. Set $\wh X:=\IPicz_{X/S}$.

Indeed, $\IPic_{X/S}$ exists, represents $\Pic_{(X/S)\et}$, and is
locally of finite type by the main theorem, Theorem~\ref{th:main}.  For
$s\in S$, the $\IPicz_{X_s/k_s}$ are smooth and proper of dimension $g$
by \cite[\S\,13]{Mm70}.  Hence $\wh X$ exists by
Proposition~\ref{prp:P0}, and is projective by Exercise~\ref{ex:sm=>pr}.
Finally, $\wh X$ is smooth owing to a more sophisticated version of the
proof of Proposition~\ref{prp:H2}; see \cite[pp.~117--118]{Mm65}.

There exists a universal sheaf $\mc P$ on $X\x\IPic_{X/S}$ by
Exercise~\ref{ex:univshf}
 since $f$ has a section $g$, namely, the identity section.  Normalize
$\mc P$ by tensoring it with $f_{\wh X}^*g_{\wh X}^* \mc P$.  Then its
restriction to $X\x\wh X$ defines a map, which is a ``duality''
isomorphism
 \def\ddX{X^{**}}  
  $$\pi\:X\risom\ddX.$$
 Indeed, forming $\pi$ commutes with changing $S$, and $\pi$'s
geometric fibers are isomorphisms by \cite[Cor., p.~132]{Mm70}.  But $X$
and  $\ddX$ are proper over $S$, and $X$  is flat.    Therefore, $\pi$
is an isomorphism by \cite[4.6.7]{EGAIII1}.
 \end{rmk}

\begin{rmk}\label{rmk:Alb}
Assume $S$ is the spectrum of an algebraically closed field $k$, and $X$
is normal, integral, and projective.  Then $\IPicz_{X/k}$ is irreducible
and projective by Proposition~\ref{prp:pic0} and Theorem~\ref{th:qpp&p}.
Set $P:=(\IPicz_{X/k})_{\red}$ and $A:=\IPicz_{P/k}$.  Then $P$ is
plainly an Abelian variety; whence, $A$ is an Abelian variety too by
Remark~\ref{rmk:Ablsch}.

Fix a point $x\in X(k)$.  Let $B$ be an Abelian variety, and set
$B^*:=\IPicz_{B/k}$.  Then $B^*$ is an Abelian variety too, and there is
a canonical isomorphism $B\risom B^{**}$ by Remark~\ref{rmk:Ablsch}.
Let $\xi\:X\to X\x B^*$ be the map defined by $0\in B^*(k)$, and
let $\beta\:B^*\to X\x B^*$ be the map defined by $x$. 

Consider a map $a\:B^*\to P$ such that $a(0)=0$.  By the Comparison
Theorem, Theorem~\ref{th:cmp}, $a$ corresponds to an invertible sheaf
$\mc L$ on $X\x B^*$ such that $\xi^*\mc L\simeq \mc O_X$.  Normalize
$\mc L$ by tensoring it with $(\beta^*\mc L)_X$.  Then $\mc L$ defines a
map $b\:X\to B$ such that $b(x)=0$.

Reversing the preceding argument, we see that every such $b$ arises from
a unique map $a\:B^*\to \IPic_{X/k}$ such that $a(0)=0$.  Since $B^*$ is
integral, $a$ factors through $P$.  Thus the maps $a\:B^*\to P$ and
$b\:X\to B$ are in bijective correspondence.  Plainly, this
correspondence is compatible with maps $b'\:B\to B'$ such that
$b'(0)=0$.  In particular, $1_P$ corresponds to a natural map $u\:X\to
A$ such that $u(x)=0$, and every map $b\:X\to B$ factors uniquely
through $u$.
 \end{rmk}

\begin{rmk}\label{rmk:Jac}
 Assume $X/S$ is projective and smooth, its geometric fibers are
connected curves of genus $g>0$, and $S$ is Noetherian.  Set
$J:=\IPicz_{X/S}$; it exists and is a projective Abelian $S$-scheme by
Exercise~\ref{ex:jac}.  Set $\wh J:=\IPicz_{J/S}$; it exists, is a
projective Abelian scheme, and is ``dual'' to $J$ by
Remark~\ref{rmk:Ablsch}.

Suppose $X$ has an invertible sheaf  $\mc L$ whose fibers $\mc L_s$ are
of degree 1.  Define an associated ``Abel'' map
	$$A_{\mc L}\:X\to J$$
 directly on $T$-points as follows.  Given $t\:T\to X$, its graph
subscheme $\Gamma_t\subset X\x T$ is a relative effective divisor; see
Answer~\ref{ex:Abel}.  Use $\mc L_T\ox\mc
O_{\!X_T}(-\Gamma)$ to define $A_{\mc L}(t)$.

Then  $A_{\mc L}$  induces, via pullback, an ``auto-duality'' isomorphism
	$$A_{\mc L}^*\:\wh J\risom J.$$
 This isomorphism is independent of the choice of $\mc L$; in fact, it
exists even if no $\mc L$ does.  These facts are proved in
\cite[Thm.~2.1, p.~595]{EGK}.  In fact, a more general autoduality
result is proved: it applies to the natural compactification of $J$,
which parameterizes torsion-free sheaves, when the geometric fibers of
$X$ are not necessarily smooth, but integral with double points at
worst.  And the proof starts from scratch, recovering the original
case of a single smooth curve over a field.
 \end{rmk}
}

\begin{rmk}\label{rmk:Jacsp}
 Assume $X/S$ is proper and flat.  Assume its geometric fibers are
curves, but not necessarily integral.  Then there are two remarkable
theorems asserting the existence of $\IPic_{X/S}$ as an algebraic space
and of $\IPicz_{X/S}$ as a separated $S$-scheme.  These theorems are
important in the theory of N\'eron models; so in \cite[Sect.~9.4]{BLR},
their proofs are sketched, and the original papers, cited.

One theorem is due to Raynaud.  He assumes, in addition, that $S$ is the
spectrum of a discrete valuation ring, that $X$ is normal, and that $\mc
O_{\!S}\risom f_*\mc O_{\!X}$ holds.  Furthermore, given any geometric
fiber of $X/S$, he measures the lengths of the local rings at the
generic points of its irreducible components, and he assumes their
greatest common divisor is 1.  Then he proves the above existence
assertions.

The other theorem is due to Deligne.  Instead, he assumes, in addition,
$X/S$ is semi-stable; that is, its geometric fibers are reduced and
connected, and have, at worst, ordinary double points.  Then he proves,
in addition, $\IPicz_{X/S}$ is smooth and quasi-projective; in fact, it
carries a canonical $S$-ample invertible sheaf.
 \end{rmk}

\section{The torsion component of the identity}\label{sc:Pictau}
 This section establishes the two main finiteness theorem for
$\IPic_{X/S}$, when $X/S$ is projective and its geometric fibers are
integral.  The first theorem asserts the finiteness of the torsion
component $\IPict_{X/S}$, an open and closed group subscheme.  By
definition, it consists of the points with a multiple in the connected
component $\IPicz_{X/S}$, which was studied in the previous section.

The second theorem asserts the finiteness of a larger sort of subset
$\IPicp_{X/S}$.  Its points represent the invertible sheaves with a
given Hilbert polynomial $\phi$.  The section starts by developing
numerical characterizations of $\IPict_{X/S}$, or rather of the
corresponding invertible sheaves, when $S$ is the spectrum of an
algebraically closed field.  This development assumes some familiarity
with basic intersection theory, which is developed in Appendix B.

\begin{dfn}\label{dfn:numeq}
 Assume $S$ is the spectrum of a field.  Let $\mc L$ and $\mc N$ be
invertible sheaves on $X$. Then $\mc L$ is said to be {\it
$\tau$-equivalent to} $\mc N$ if, for some nonzero $m$ depending on $\mc
L$ and $\mc N$, the $m$th power $\mc L^{\ox m}$ is algebraically
equivalent to $\mc N^{\ox m}$.

In addition, $\mc L$ is said to be {\it numerically equivalent to} $\mc
N$ if, for every complete curve $Y\subset X$, the corresponding
intersection numbers are equal:
     $$\textstyle\int c_1\mc L \cdot [Y]=\int c_1\mc N \cdot [Y].$$
 It is sufficient, by additivity, to take $Y$ to be complete and
integral.  It is then equivalent that $\deg\mc L_Y=\deg\mc N_Y$ or that
$\deg\mc L_{Y'}=\deg\mc N_{Y'}$ where $Y'$ is the normalization of $Y$,
because, in any event,
   $$\textstyle\int c_1\mc L \cdot [Y]=\deg\mc L_Y=\deg\mc L_{Y'}.$$
 \end{dfn}

\begin{dfn}\label{dfn:bddset}
 Assume $S$ is the spectrum of a field.  Let $\Lambda$ be a family of
invertible sheaves on $X$.  Then $\Lambda$ is said to be {\it bounded\/}
if there exist an $S$-scheme $T$ of finite type and an invertible sheaf
$\mc M$ on $X_T$ such that, given $\mc L\in \Lambda$, there exists a
geometric point $t$ of  $T$ such that $\mc L_t\simeq\mc M_t$.
 \end{dfn}

\begin{thm}\label{thm:numeq}
  Assume $S$ is the spectrum of an algebraically closed field, and $X$
is projective.  Let $\mc L$ be an invertible sheaf on $X$.  Then the
following conditions are equivalent:\smallskip

\tu{(a)} The sheaf $\mc L$ is $\tau$-equivalent to $\mc O_{\!X}$.

\tu{(b)} The sheaf $\mc L$ is numerically equivalent to $\mc O_{\!X}$.

\tu{(c)} The family $\{\mc L^{\ox p}\mid p\in\bb Z\}$ is bounded.

\tu{(d)} For every coherent sheaf $\mc F$ on $X$, we have
 $\chi(\mc F\ox\mc L)=\chi(\mc F)$.
 
\tu{(e)} For every closed integral curve $Y\subset X$, we have 
 $\chi(\mc L_Y)=\chi(\mc O_{\!Y})$.
 
\tu{(f)} For every integer $p$, the sheaf $\mc L^{\ox p}(1)$ is ample.
\smallskip

\noindent If  $X$ is irreducible, then all the conditions above are
equivalent to the following one:\smallskip
 
\tu{(g)} For every pair of integers $p,\,n$, we have
 $\chi(\mc L^{\ox p}(n))=\chi(\mc O_{\!X}(n))$.\smallskip

\noindent If $X$ is irreducible of dimension $r\ge2$, then all the
conditions above are equivalent to the following one:\smallskip
 
\tu{(h)} Setting $\ell:=c_1\mc L$ and $h:=c_1\mc O_{\!X}(1)$, we have
 $\int \ell h^{r-1}=0$ and $\int\ell^2 h^{r-2}=0$.
 \end{thm}
\begin{proof}
 Let us proceed by establishing the following implications:
\begin{gather*}
 \tu{(c)} \To \tu{(a)} \To \tu{(d)} \To \tu{(e)} \Longleftrightarrow
\tu{(b)} \To\tu{(c)}
 \To \tu{(f)} \To \tu{(b)};\\
 \tu{(d)} \To \tu{(g)} \To \tu{(h)}
 \To \tu{(b)}; \text{ and } \tu{(g)} \To \tu{(b)} \text{ if } \dim X=1. 
 \end{gather*}

Assume (c).  Then, by definition, there exist an $S$-scheme $T$ of
finite type and an invertible sheaf $\mc M$ on $X_T$ such that, given $
p\in\bb Z$, there exists a geometric point $t$ of $T$ such that $\mc
L_t^{\ox p}=\mc M_t$.  Apply the Pigeonhole Principle: say $\mc L^{\ox
p_1}$ and $\mc L^{\ox p_2}$ belong to the same connected component of
$T$, but $p_1\neq p_2$.  Set $m:=p_1-p_2$.  Then $\mc L^{\ox m}$ is
algebraically equivalent to $\mc O_X$.  Thus (a) holds.

Furthermore, for each $t\in T$, there exists an $n$ such that $\mc M_t(n)$
is ample by \cite[4.5.8]{EGAII}.  So $t$ has a neighborhood $U$ such
that $\mc M_u(n)$ is ample for every $u\in U$ by \cite[4.7.1]{EGAIII1}.
Since $T$ is quasi-compact, $T$ is covered by finitely many of the $U$.
Let $N$ be the product of the corresponding $n$.  Then $\mc M_t(N)$ is
ample for every $t\in T$.  In particular, $\mc L^{\ox p}(N)$ is ample
for every $ p\in\bb Z$.  So $\bigl(\mc L^{\ox q}(1)\bigr)^{\ox N}$ is
ample for every $ q\in\bb Z$.  Thus (f) holds.

Assume (a).  The function $n\mapsto\chi(\mc F\ox\mc L^{\ox n})$ is a
polynomial.  To prove it is constant, we may replace $\mc L$ by $\mc
L^{\ox m}$ for any nonzero $m$.  Thus we may assume $\mc L$ is
algebraically equivalent to $\mc O_X$.  So let $T$ be a connected
$S$-scheme, and $\mc M$ an invertible sheaf on $X_T$.  Then for a fixed
$n$, as $t\in T$ varies, the function $t\mapsto \chi(\mc F\ox\mc
L_t^{\ox n})$ is constant by \cite[7.9.5]{EGAIII2}.  It follows that (d)
holds.

Assume (d).  Taking $\mc F:=\mc O_{\!Y}$, we get (e).  Taking $\mc F:=\mc
L^{\ox p}(n)$, we get $\chi(\mc L^{\ox (p+1)}(n))=\chi(\mc L^{\ox
p}(n))$.  Thus whether or not $X$ is irreducible, (g) holds,.

Assume (g).  Then $X$ is irreducible, say of dimension $r$.  Set
$\ell:=c_1\mc L$ and $h:=c_1\mc O_{\!X}(1)$.
 Write
	$$\chi(\mc L^{\ox p}(n))=
 \sum_{0\le i,\,j\le r}	a_{ij}\binom{p+i}i\binom{n+j}j$$
 where $a_{ij}=\int\ell^ih^j$ if $i+j=r$.  Then (g) implies $a_{ij}=0$
if $i\ge1$.  If $r\ge2$, then (h) follows.  Suppose $r=1$.  Then
$\int\ell=0$.  Now, $X_{\red}$ is the only closed integral curve
contained in $X$, and $\int\ell$ is a multiple of
$\int\ell\cdot[X_{\red}]$.  Thus (b) holds.

Conditions (e) and (b) are equivalent since $\chi(\mc L_Y)=\deg(\mc
L_Y)+\chi(\mc O_{\!Y})$ by Riemann's Theorem.

Assume (f).  Then for every closed integral curve $Y\subset X$, we have 
	$$0\le\deg(\mc L_Y^{\ox p}(1))=p\deg(\mc L_Y)+\deg(\mc O_{\!Y}(1))$$
  for every integer $p$.  So $\deg(\mc L_Y)=0$.  Thus (b) holds.

Assume (h).  Then $X$ is irreducible of dimension $r\ge2$.  To prove
(b), plainly we may replace $X$ by its reduction.  We proceed by
induction on $r$.  If $r=2$, then (b) holds by the Hodge Index Theorem.

Suppose $r\ge3$.  Given a complete integral curve $Y\subset X$, take $n$
so that the twisted ideal $\mc I_{Y,X}(n-1)$ is generated by its global
sections.  View these sections as sections of $\mc O_{\!X}(n-1)$.  Then they
define a linear system that is free of base points on $X-Y$.  So the
global sections of $\mc I_Y(n)$ define a linear system that is very
ample on $X-Y$.  In particular, this system maps $X-Y$ onto a variety of
dimension at least 2; in other words, the system is not ``composite with
a pencil.''

Hence the generic member $H_\eta$ is geometrically irreducible by
\cite[Thm.~I.6.3, p.~30]{Za58}.  In the first instance, we must apply
the cited theorem to the induced system on the normalization of $X$, and
we conclude that the preimage of $H_\eta$ is geometrically irreducible.
But then $H_\eta$ is too.  Therefore, by \cite[9.7.7]{EGAIV3}, a general
member $H$ is irreducible.

Set $\ell_1:=c_1\mc L_H$ and $h_1:=c_1\mc O_{\!H}(1)$.  Then, by the
Projection Formula,
  $$\textstyle\int \ell_1 h_1^{r-2}=n\int \ell h^{r-1}=0 \text{ and }
	\int\ell_1^2 h_1^{r-3}=n\int\ell^2 h^{r-2}=0.$$
 So by induction, $\mc L_H$ is numerically equivalent to $\mc O_{\!X}$
on $H$.  But $Y\subset H$ since $H$ arises from a section of $\mc
I_Y(n)$.  Hence, by the Projection Formula,
	  $$\textstyle\int\ell\cdot[Y]=\int\ell_1\cdot[Y]=0.$$
 Thus (b) holds.

Finally, assume (b).  By Lemma~\ref{lem:hp} below, there is an $m$ such
that, if $\mc N$ is an invertible sheaf on $X$ numerically equivalent to
$\mc O_{\!X}$, then $\mc N$ is $m$-regular.  So $\mc N(m)$ is generated
by its global sections, and its higher cohomology groups vanish.

Set $\phi(n):=\chi(\mc O_{\!X}(n))$ and $M:=\phi(m)$.  Then $\dim\tu
H^0(\mc N(m))=M$, since $\chi(\mc N(n))=\phi(n)$ also by
Lemma~\ref{lem:hp} below.  Set $\mc F:=\mc O_{\!X}(-m)^{\oplus M}$.
Then $\mc N$ is a quotient of $\mc F$.

Set $T:=\IQuot^\phi_{\mc F/X/k}$.  Then $T$ is of finite type.  Let $\mc
M$ be the universal quotient.  Then there exists a $k$-point $t\in T$
such that $\mc N=\mc M_t$.  Let $U\subset X\x T$ be the open set on
which $\mc M$ is invertible.  Let $R\subset T$ be the image of the
complement of $U$.  Then $R$ is closed.  Replace $T$ by $T-R$, and $\mc
M$  by its restriction.  Then $t\in T$ still.  Thus the invertible
sheaves on $X$ numerically equivalent to $\mc O_{\!X}$ form a bounded
family.  In particular, (c) holds.
 \end{proof}

\begin{ex}\label{ex:Altpf}
 Consider the preceding paragraph, the last one in the proof of
Theorem~\ref{thm:numeq}.  Using $\IDiv_{X/k}$ instead of
$\IQuot^\phi_{\mc F/X/k}$, give another proof that the invertible
sheaves $\mc N$ on $X$ numerically equivalent to $\mc O_{\!X}$ form a
bounded family.
 \end{ex}

\begin{lem}\label{lem:bd}
  Assume $S$ is the spectrum of an algebraically closed field, and $X$
is projective of dimension $r$.  Let $\mc F$ be a coherent sheaf on $X$.
Then there is a number $B_{\mc F}$ such that, if $\mc L$ is any invertible
sheaf on $X$ numerically equivalent to $\mc O_{\!X}$, then
 $$\dim \tu H^0(\mc L\ox\mc F(n)) \le B_{\mc F}\tbinom{n+r}r
 \text{ for all } n\ge0.$$
 \end{lem}
\begin{proof}
 Suppose $r=0$.  Then $\mc L\ox\mc F(n)=\mc F$.  So we may take $B_{\mc
F}=\dim \tu H^0(\mc F)$.

Suppose $r\ge1$.  Given a short exact sequence $0\to\mc F'\to\mc F\to\mc
F''\to0$, we have
	$$\dim \tu H^0(\mc L\ox\mc F(n))
 \le \dim \tu H^0(\mc L\ox\mc F'(n))+\dim \tu H^0(\mc L\ox\mc F''(n)).$$
 So given $B_{\mc F'}$ and  $B_{\mc F''}$, we can take $B_{\mc F}:=
B_{\mc F'}+B_{\mc F''}$.  

 Say $X=\Proj(A)$ and $\mc F=\wt M$ where $M$ is a finitely generated
graded $A$-module.  Then there is a filtration by graded submodules
 $$M=:M_q\supset M_{q-1} \supset\cdots\supset M_1\supset M_0:=0$$
 such that $M_{i+1}/M_i\simeq (A/P_i)[p_i]$ where $P_i$ is a homogeneous
prime for each $i$.  It follows that we may assume $X$ is integral and
$\mc F=\mc O_{\!X}(p)$.

Let $\mc L$ be numerically equivalent to $\mc O_{\!X}$.  Set $\ell:=c_1\mc
L$ and $h:=c_1\mc O_{\!X}(1)$.  Suppose $\mc L(p)$ has a nonzero section.
It defines a divisor $D$, possibly 0.  Hence
  $$\textstyle 0\le\int h^{r-1}[D]=\int h^{r-1}\ell+p\int h^r.$$
 But $\int h^{r-1}\ell=0$ because $h^{r-1}$ is represented by a curve
since $r\ge1$.  And $\int h^r>0$.  Hence $p\ge0$. Thus $\tu H^0(\mc
L(-1))=0$.

Let $H$ be a hyperplane section of $X$.  Then there is an exact sequence
	$$0\to\mc L(n-1)\to\mc L(n)\to \mc L_H(n)\to0.$$
 By induction on $r$, we may assume there is a number $B$ such that
  $$\dim \tu H^0(\mc L_H(n)) \le B\tbinom{n+r-1}{r-1} \text{ for all }
  n\ge0;$$
 moreover, $B$ works for every $\mc L$.  Hence
	$$\dim \tu H^0(\mc L(n)) -\dim \tu H^0(\mc L(n-1)) \le
	B\tbinom{n+r-1}{r-1} \text{ for all } n\ge0.$$
 But $\tu H^0(\mc L(-1))=0$.  Since $\binom{n+r-1}r+\binom{n+r-1}{r-1} =
\binom{n+r}r$, induction on $n$ yields
  $$\dim \tu H^0(\mc L(n)) \le B\tbinom{n+r}r\text{ for all }n\ge0.$$

Recall $\mc F=\mc O_{\!X}(p)$.  If $p\le0$, then $\mc F\subset\mc O_{\!X}$;
so we may take $B_{\mc F}:=B$.  But if $p\ge0$, then $\binom{p+n+r}r
\le \binom{p+r}r\binom{n+r}r$ since every monomial of degree $p+n$ is
the product of one of degree $p$ and one of degree $n$; so we may
take $B_{\mc F}:=B\binom{p+r}r$.
 \end{proof}

\begin{lem}\label{lem:hp}
  Assume $S$ is the spectrum of an algebraically closed field, and $X$
is projective.  Then there is an integer $m$ such that, if $\mc L$ is any
invertible sheaf on $X$ numerically equivalent to $\mc O_{\!X}$, then $\mc L$ is
$m$-regular, and 
	$$\chi(\mc L(n))=\chi(\mc O_{\!X}(n)) \text{ for all } n.$$
 \end{lem}
\begin{proof}
 Set $r:=\dim(X)$, and proceed by induction on $r$.  If $r=0$, then both
assertions are trivial.  So assume $r\ge1$. 

First, let us establish the asserted equation.  Given an $\mc L$, fix
$q\ge1$ such that $\mc L(q)$ is very ample.  Take effective divisors $F$
and $G$ such that
 $$\mc O_{\!X}(F)=\mc O_{\!X}(q) \text{ and } \mc O_{\!X}(G)=\mc L(q).$$
  For every $p$, plainly $\mc L^{\ox p}_F$ and $\mc L^{\ox p}_G$ are
numerically equivalent to $\mc O_{\!F}$ and $\mc O_{\!G}$.

Form $0\to\mc O_{\!X}(-F)\to\mc O_{\!X}\to\mc O_{\!F}\to0$ and $0\to\mc O_{\!X}(-G) \to
\mc O_{\!X}\to\mc O_{\!G}\to0$.  Tensor them with $\mc L^{\ox p}(n+q)$.  We
get
 \begin{align}
 0\to\mc L^{\ox p}&(n)\to\mc L^{\ox p}(n+q)\to\mc L^{\ox p}_F(n+q)\to0,
 \label{seq:F}\\
 0\to\mc L^{\ox p-1}&(n)\to\mc L^{\ox p}(n+q)\to\mc L^{\ox p}_G(n+q)\to0.
  \notag\end{align}
Apply $\chi(\bullet)$ and subtract.  We get
	$$\chi(\mc L^{\ox p}(n))-\chi(\mc L^{\ox p-1})(n))
	=\chi(\mc L_G^{\ox p}(n+q))-\chi(\mc L_F^{\ox p}(n+q)).$$
 By induction, the right hand side varies as a polynomial in $n$, which
is independent of $p$.  Hence there are polynomials $\phi_1(n)$ and
$\phi_0(n)$ such that
 \begin{equation}\label{equation:phi}
 \chi(\mc L^{\ox p}(n))=\phi_1(n)p+\phi_0(n).
 \end{equation}
Suppose $\phi_1\neq0$.  Say $\phi_1(n)\neq0$ for all $n\ge n_1$.

By induction, there is an integer $n_2$ such that $\mc L^{\ox p}_F$ is
$n_2$-regular for every $p$.  So $\tu H^i(\mc L^{\ox p}_F(n))=0$ for
$i\ge1$, for $n\ge n_2-i$, and for every $p$.  Hence, owing to
Sequence~(\ref{seq:F}), there is an isomorphism
	$$\tu H^i(\mc L^{\ox p}(n))\risom\tu H^i(\mc L^{\ox p}(n+q))
	 \text{ for $i\ge2$, for $n\ge n_2$, and for every $p$}.$$
 Now, for each $i\ge2$, each $p$, and each $n$, there is a $j\ge0$ such that
$$\tu H^i(\mc L^{\ox p}(n+jq))=0$$ by Serre's Theorem.  Therefore,
	$$\tu H^i(\mc L^{\ox p}(n))=0
	 \text{ for $i\ge2$, for $n\ge n_2$ and for every $p$}.$$

Hence $\tu H^0(\mc L^{\ox p}(n))\ge\chi(\mc L^{\ox p}(n))$ for $n\ge n_2$
and every $p$.  Take $n:=\max(n_2,n_1)$.  Owing to
Equation~(\ref{equation:phi}), then $\tu H^0(\mc L^{\ox p}(n))\to\infty$
as $p\to\infty$ if $\phi_1(n)>0$ or as $p\to-\infty$ if $\phi_1(n)<0$.
However, by Lemma~(\ref{lem:bd}), there is a number $B$ such that $\tu
H^0(\mc L^{\ox p}(n))\le B$ for any $p$.  This contradiction means
$\phi_1=0$.  Hence Equation~(\ref{equation:phi}) yields
$\chi(\mc L(n))=\chi(\mc O_{\!X}(n))$ for all $n$, as desired.

Finally, in order to prove there is an $m$ such that every $\mc L$
numerically equivalent to $\mc O_{\!X}$ is $m$-regular, we must modify Mumford's
original work \cite[Lect.~14]{Mm66} because these $\mc L$ are not ideals.
However, as Mumford himself points out \cite[pp.~102--103]{Mm66}, the
hypothesis that his sheaf $\mc I$ is an ideal enters only through the
bound $\dim\tu H^0(\mc I(n))\le \binom{n+r}r$.  Plainly, this bound can
be replaced by the bound $B_{\mc F}$ with $\mc F:=\mc O_{\!X}$ of
Lemma~\ref{lem:bd}.  Of course, in addition, we must use the fact we
just proved, that all the $\mc L$ have the same Hilbert polynomial.
 \end{proof}

\begin{ex}\label{ex:h0L}
 Assume $S$ is the spectrum of an algebraically closed field, and $X$ is
projective and integral of dimension $r\ge1$.  Set $h:=c_1\mc O_{\!X}(1)$.
Let $\mc L$ be an invertible sheaf on $X$, and set $\ell:=c_1\mc L$.  Say
	$$\chi(\mc L(n))=\textstyle\sum_{0\le i\le r}a_i\binom{n+i}i
	\text{ and } a:=\int\ell h^{r-1}.$$
 
Suppose $a<a_r$.  Show $\dim \tu H^0(\mc L(n)) \le a_r\binom{n+r}r$ for
all $n\ge0$.  Furthermore, modifying Mumford's work
\cite[pp.~102--103]{Mm66} slightly, show there is a polynomial $\Phi_r$
depending only on $r$ such that $\mc L$ is $m$-regular with
$m:=\Phi_r(a_0,\dotsc,a_{r-1})$.

In general, show there is a polynomial $\Psi_r$ depending only on $r$
such that $\mc L$ is $m$-regular with $m:=\Psi_r(a_0,\dotsc,a_r;a)$.
 \end{ex}

\begin{dfn}\label{dfn:Gtau}
 Let $G/S$ be a group scheme.  For $n>0$, let $\vf_n\:G\to G$ denote the
$n$th power map.  Then $G^\tau$ is the set defined by the formula
	$$G^\tau:=\textstyle\bigcup_{n>0}\vf_n^{-1}G^0$$
 where $G^0$ is the union of the connected components of the
identity $G_s^0$ for $s\in S$.
 \end{dfn}

\begin{lem}\label{lem:Gtauk}
 Let $k$ be a field, and $G$ a commutative group scheme locally of
finite type.  Then $G^\tau$ is an open group subscheme, and forming it
commutes with extending $k$.  Moreover, if $G^\tau$ is quasi-compact,
then it is closed and of finite type.
 \end{lem}
\begin{proof}
 By Lemma~\ref{lem:agps}, $G^0$ is an open and closed group subscheme of
finite type, and forming it commutes with extending $k$.  Now, the $n$th
power map $\vf_n$ is continuous, and forming it commutes with extending
$k$; also, $\vf_n$ is a homomorphism since $G$ is commutative.  So
$G^\tau$ is the filtered union of the open and closed group subschemes
$\vf_n^{-1}G^0$, and forming them commutes with extending $k$.  Hence
$G^\tau$ is an open group subscheme, and forming it commutes with
extending $k$.  Moreover, if $G^\tau$ is quasi-compact, then 
  $G^\tau$ is the union of finitely many of the $\vf_n^{-1}G^0$, and so
$G^\tau$ is closed; also, then $G^\tau$ is of finite type since $G$ is
locally of finite type.
 \end{proof}

\begin{ex}\label{ex:Gtauk}
 Let $k$ be a field, and $G$ a commutative group scheme locally of
finite type.  Let $H\subset G$ be a group subscheme of finite type.
Show  $H\subset G^\tau$.  (Thus, if $G^\tau$ is of finite type, then it
is the largest group subscheme of finite type.)
 \end{ex}

\begin{ex}\label{ex:teq}
 Assume $S$ is the spectrum of a field $k$.  Assume $\IPic_{X/k}$ exists
and represents $\Pic_{(X/k)\fppf}$.  Let $\mc L$ be an invertible sheaf
on $X$, and $\lambda\in\IPic_{X/k}$ the corresponding point.  Show $\mc
L$ is  $\tau$-equivalent to $\mc O_{\!X}$ if and only if
$\lambda\in\IPict_{X/k}$.
 \end{ex}

\begin{prp}\label{prp:pictk}
 Assume $S$ is the spectrum of a field $k$.  Assume $\IPic_{X/k}$ exists
and represents $\Pic_{(X/k)\fppf}$.  Then $\IPict_{X/k}$ is an open
group subscheme, and forming it commutes with extending $k$.  Moreover,
if $X$ is projective, then $\IPict_{X/k}$ is closed and of finite type.
 \end{prp}
 \begin{proof}
 By Proposition~\ref{prp:lft}, $\IPic_{X/S}$ is locally of finite type.
So owing to Lemma~\ref{lem:Gtauk}, we need only prove $\IPict_{X/k}$ is
quasi-compact when $X$ is projective.  Since forming $\IPict_{X/k}$
commutes with extending $k$, we may also assume $k$ is algebraically
closed.

At the very end of the proof of Theorem~\ref{thm:numeq}, we proved that
the invertible sheaves $\mc L$ numerically equivalent to $\mc O_{\!X}$ form
a bounded family.  In other words, there is a $k$-scheme $T$ of finite
type and an invertible sheaf $\mc M$ on $X_T$ such that the $\mc L$
appear among the fibers $\mc M_t$.

  Then $\mc M$ defines a map $\theta\:T\to \IPic_{X/k}$.  Owing to
Theorem~\ref{thm:numeq} and to Exercise~\ref{ex:teq}, we have $\theta(T)
\supset \IPict_{X/k}$.  Since $T$ is Noetherian, so is $\theta(T)$;
whence, so is any subspace of $\theta(T)$.  Thus $\IPict_{X/k}$ is
quasi-compact, as needed.
 \end{proof}

\begin{ex}\label{ex:ptauproj}
  Assume $S$ is the spectrum of a field $k$.  Assume $X$ is projective
and geometrically integral.  Show $\IPict_{X/k}$ is quasi-projective.
If also $X$ is geometrically normal, show $\IPict_{X/k}$ is projective.
 \end{ex}

\begin{rmk}\label{rmk:numeq}
 In Proposition~\ref{prp:pictk}, if $X$ is projective, then
$\IPic_{X/k}$ does exist and represent $\Pic_{(X/k)\fppf}$ according to
Corollary~\ref{cor:algsch}.  In fact, this corollary asserts
$\IPic_{X/k}$ exists and represents $\Pic_{(X/k)\fppf}$ whenever $X$ is
complete; $X$ need not be projective.  Furthermore, although we used
projective methods to prove Proposition~\ref{prp:pictk}, we can infer it
whenever $X$ is complete, as follows.

Assume $X$ is complete.  By Chow's lemma, there is a projective variety
$X'$ and a surjective map $\gamma\:X'\to X$.  By Theorem~\ref{th:rrep},
the induced map
	$$\gamma^*\:\IPic_{X/k}\to\IPic_{X'/k}$$
  is of finite type.  Set $H:=(\gamma^*)^{-1}\IPict_{X'/k}$.  Then $H$
is of finite type since $\IPict_{X'/k}$ is by
Proposition~\ref{prp:pictk}.  Now, plainly $\gamma^*\IPicz_{X/k}\subset
\IPicz_{X'/k}$.  So since $\gamma^*$ is a homomorphism,
$\gamma^*\IPict_{X/k}\subset\IPict_{X'/k}$; whence, $\IPict_{X/k}\subset
H$ (in fact, the two are equal by Exercise~\ref{ex:Gtauk}).  Since
$\IPict_{X/k}$ is open, it is therefore a subscheme of finite type.

 Similarly, in Theorem~\ref{thm:numeq}, Conditions (a)--(e) continue to
make sense and to remain equivalent whenever $X$ is complete.  Indeed,
our proofs of the implications
 $$\tu{(c)} \To \tu{(a)} \To \tu{(d)} \To \tu{(e)} \Longleftrightarrow
	\tu{(b)} \text{ and } \tu{(c)}\To \tu{(f)} \To \tu{(b)}$$
 work without change.  However, we used projective methods to prove
$\tu{(b)} \To\tu{(c)}$.  Nevertheless, we can infer this implication
whenever $X$ is complete, as follows.

Let $\mc L$ be numerically equivalent to $\mc O_{\!X}$.  Let $\gamma\:X'\to
X$ be as above.  Then $\gamma^*\mc L$ is numerically equivalent to $\mc
O_{\!X'}$.  Indeed, Let $Y'\subset X'$ be a closed curve.  Then $\int
c_1\gamma^*\mc L \cdot [Y']$ is a multiple of $\int c_1\mc L \cdot
[\gamma Y']$ by the Projection Formula; hence, the former number
vanishes as the latter does.

Let $\lambda\in\IPic_{X/k}$ represent $\mc L$.  Then
$\gamma^*\lambda\in\IPic_{X'/k}$ represents $\gamma^*\mc L$.  Now,
$\gamma^*\mc L$ is numerically equivalent to $\mc O_{\!X'}$.  Hence
$\gamma^*\lambda\in\IPict_{X'/k}$ owing to Theorem~\ref{thm:numeq} and
to Exercise~\ref{ex:teq}.  So $\lambda\in H:=(\gamma^*)^{-1}
\IPict_{X'/k}$.

The inclusion $H\into \IPic_{X/k}$ is defined by an invertible sheaf
$\mc M$ on $X_T$ for some fppf covering $T\to H$.  (Although $k$ is
algebraically closed, possibly $\mc O_{\!S}\risom f_*\mc O_{\!X}$ does
not hold universally, so we cannot simply take $T:=H$.)  Replace $T$ by
an open subscheme so that $T\to H$ is of finite type and surjective.
Since $H$ is of finite type, so is $T$.

Let $t\in T$ be a $k$-point that maps to $\lambda\in H$.  Then $\mc
M_t\simeq \mc L$.  Now, for every $p\in \bb Z$, plainly $\mc L^{\ox p}$
is numerically equivalent to $\mc O_{\!X}$.  So similarly $\mc L^{\ox p}
\simeq \mc M_{t_p}$ for some $k$-point $t_p\in T$.  Thus (c) holds.
 \end{rmk}

\begin{ex}\label{ex:bddsbset}
 Assume $\IPic_{X/S}$ exists and represents $\Pic_{(X/S)\fppf}$.  Let
$\Lambda$ be an arbitrary subset of $\IPic_{X/S}$, and $L$ the
corresponding family of classes of invertible sheaves on the fibers of
$X/S$ in the sense of Exercise~\ref{ex:schpts}.  Show $\Lambda$ is
quasi-compact (with the induced topology) if and only if $L$ is {\it
bounded} in the following sense: there exist an $S$-scheme $T$ of finite
type and an invertible sheaf $\mc M$ on $X_T$ such that every class in
$L$ is represented by a fiber $\mc M_t$ for some $t\in T$.
 \end{ex}

\begin{thm}\label{th:Ptaufin}
 Assume $f\:X\to S$ is projective Zariski locally over $S$, and flat
with integral geometric fibers.  Then $\IPict_{X/S}$ is an open and
closed group subscheme of finite type, and forming it commutes with
changing $S$.  If also $X/S$ is projective, and $S$ is Noetherian, then
$\IPict_{X/S}$ is quasi-projective.
 \end{thm}
\begin{proof}
 The second assertion follows from the first and
Exercise~\ref{ex:q-proj}.  The first assertion is local on $S$; so, to
prove it, we may assume $X/S$ is projective.

 Theorem~\ref{th:main} asserts $\IPic_{X/S}$ exists and is locally of
finite type.  So the $\IPict_{X_s/k_s}$ are subgroups of the
$\IPic_{X_s/k_s}$, and forming $\IPict_{X_s/k_s}$ commutes with
extending $k_s$ by Exercise~\ref{ex:bschg} and Lemma~\ref{lem:Gtauk}.
However, plainly $\IPict_{X/S}=\bigcup_{s\in S}\IPict_{X_s/k_s}$ as
sets.  Hence, forming $\IPict_{X/S}$ commutes with changing $S$;
moreover, in order to infer $\IPict_{X/S}$ is a group subscheme, we need
only prove it is open.

Plainly, a subset $A$ of a topological space $B$ is open or closed if
(and only if), for every member $B_i$ of an open covering of $B$, the
intersection $A\cap B_i$ is so in $B_i$.  Hence, in order to infer
$\IPict_{X/S}$ is an open and closed group subscheme, we need only prove
that, for any affine open subscheme $U$ of $\IPic_{X/k}$, the
intersection $U\cap \IPict_{X/S}$ is open and closed in $U$.

Theorem~\ref{th:main} also asserts $\IPic_{X/k}$ represents
$\Pic_{(X/S)\et}$.  Thus the inclusion $U\into\IPic_{X/S}$ is defined by
an invertible sheaf $\mc M$ on $X_T$ for some \'etale covering $T\to U$.
Replace $T$ by an open subscheme so that $T\to U$ is of finite type and
surjective.  Let $T^\tau$ be the set of $t\in T$ where $M_t$ is
$\tau$-equivalent to $\mc O_{\!X_t}$.  Then $T^\tau$ is the preimage of
$\IPict_{X/S}$ in $T$ owing to Exercise~\ref{ex:teq}.  So we have to
prove $T^\tau$ is open and closed.

Since $U$ is affine, it is quasi-compact, so of finite type.  Hence $T$
is of finite type.  So it has only finitely many connected components.
But, for every pair $p,\,n$, the function $t\mapsto\chi(\mc M_t^{\ox
p}(n))$  is constant on each connected component of $T$.  Therefore,
Theorem~\ref{thm:numeq} implies $T^\tau$ is open and closed, as desired.

It remains to prove $\IPict_{X/S}$ is of finite type.  Let $L$ be the
corresponding family of classes of invertible sheaves on the fibers of
$X/S$.  By Exercise~\ref{ex:bddsbset}, we need only prove $L$ is
bounded.  At the very end of the proof of Theorem~\ref{thm:numeq}, we
proved essentially this statement when $S$ is the spectrum of an
algebraically closed field, and $X$ is projective, but not necessarily
integral.  We can argue similarly here, but must make two important
modifications.

First, if a class in $L$ is represented by an invertible sheaf $\mc L$
on a fiber $X_k$ where $k$ is a field containing the field $k_s$ of a
point $s\in S$, then $\chi(\mc L(n))=\chi(\mc O_{\!X_s}(n))$ for all $n$,
owing to Lemma~\ref{lem:hp}.  But $\chi(\mc O_{\!X_s}(n))$ can vary with
$s$.  Nevertheless, it must remain the same on each connected component
of $S$.  And the matter in question is local on $S$.  So we may and must
assume $S$ is connected.

Second, at the end of the proof of Lemma~\ref{lem:hp}, when we modified
Mumford's work, we used the bound $B_{\mc F}$ with $\mc F:=\mc O_{\!X}$ of
Lemma~\ref{lem:bd}; in fact, in the induction step, we implicitly used
the corresponding bounds for various subschemes of $X$.  Unfortunately,
it is not clear, in general, how these bounds vary with $X$.  But
in place of Lemma~\ref{lem:hp}, we can use Exercise~\ref{ex:h0L}, which
provides a uniform $m$ such that $\mc L$ is $m$-regular for every $\mc
L$ representing a class in $L$.
 \end{proof}

\begin{cor}\label{cor:torgp}
 Assume $S$ is Noetherian.  Assume $f\:X\to S$ is projective Zariski
locally over $S$, and is flat with geometrically integral fibers.  For
each $s\in S$, let $k_s'$ be the algebraic closure of the residue field
$k_s$.  Then the torsion group \begin{equation}\label{equation:tg}
 \IPict_{X_{k_s'}/k_s'}(k_s')\bigm/\IPicz_{X_{k_s'}/k_s'}(k_s')
 \end{equation}
 is finite, and its order is bounded.
 \end{cor}
 \begin{proof}
 Since $\IPicz_{X_{k_s'}/k_s'}$ is open in $\IPict_{X_{k_s'}/k_s'}$ by
Proposition~\ref{prp:pic0}, the order of their quotient is equal to the
number of connected components of $\IPict_{X_{k_s'}/k_s'}$.  This number
is finite because $\IPict_{X_{k_s'}/k_s'}$ is of finite type by
Proposition~\ref{prp:pictk}.

Moreover, $\IPict_{X_{k_s'}/k_s'}$ is equal to $\IPict_{X_s/k_s}\ox
k_s'$ again by Proposition~\ref{prp:pictk}, so equal to $\IPict_{X/S}\ox
k_s'$ essentially by Definition~\ref{dfn:Gtau}.  But, since
$\IPict_{X/S}$ is of finite type by Theorem~\ref{th:Ptaufin}, the number
of connected components of $\IPict_{X/S}\ox k_s'$ is constant for $s$ in
a nonempty open subset of $S$ by \cite[9.7.9]{EGAIV3}.  Hence the number
is bounded by Noetherian induction.
 \end{proof}

\begin{ex}\label{ex:sm=>Pictpr}
 Assume $X/S$ is projective and smooth, its geometric fibers are
irreducible, and $S$ is Noetherian. Show that $\IPict_{X/S}$ is
projective.
 \end{ex}

\begin{rmk}\label{rmk:Ptaufin}
 Assume $X/S$ is proper.  If $\IPic_{X/S}$ exists and it represents
$\Pic_{(X/k)\fppf}$, then $\IPict_{X/S}$ is an open group subscheme of
finite type.  This fact can be derived from Theorem~\ref{th:Ptaufin}
through a series of reduction steps; see \cite[Thm.~4.7, p.~647]{Kl71}.

Assume $S$ is Noetherian in addition.  Then, whether or not
$\IPic_{X/k}$ exists, the torsion group~(\ref{equation:tg}) is finite
and its order is bounded.  This fact follows from the preceding one via
the proof of Corollary~\ref{cor:torgp}, since there is
an nonempty open subscheme $V$ of $S_{\red}$ such that $\IPic_{X_V/V}$
exists and represents $\Pic_{(X_V/V)\fppf}$ by Grothendieck's
Theorem~\ref{th:genrep}.

Furthermore, the rank of the corresponding ``N\'eron--Severi'' group
	$$\IPic_{X_{k_s'}/k_s'}(k_s')\bigm/\IPict_{X_{k_s'}/k_s'}(k_s')$$
 is finite and bounded.  This fact is far deeper; see \cite[Thm.~5.1,
p.~650, and Rem.~5.3, p.~652]{Kl71}, and see \cite[pp.~121--124]{Za35}.
Moreover, the rank is arithmetic in nature: its value need not be a
constructible function of $s\in S$; a standard example is discussed in
\cite[p.~235]{BLR}.
 \end{rmk}

\begin{thm}\label{thm:Pphifin}
 Assume $X/S$ is projective and flat with integral geometric fibers.
Given a polynomial $\phi\in\bb Q[n]$, let $\IPicp_{X/S}\subset
\IPic_{X/S}$ be the set of points representing invertible sheaves $\mc
L$ such that $\chi(\mc L(n))=\phi(n)$ for all $n$.  Then the
$\IPicp_{X/S}$ are open and closed subschemes of finite type; they are
disjoint and cover; and forming them commutes with changing $S$.  If
also $S$ is Noetherian, then $\IPicp_{X/S}$ is quasi-projective.
 \end{thm}
\begin{proof}
 Plainly, the $\IPicp_{X/S}$ are disjoint and cover, and forming them
commutes with changing $S$.  The rest of the proof is similar to the
proof of Theorem~\ref{th:Ptaufin}.  In fact, the present case is simpler
because the sheaves in question have the same Hilbert polynomials by
hypothesis; there is no need to appeal to Theorem~\ref{thm:numeq} nor to
Lemma~\ref{lem:hp}.
 \end{proof}

\begin{ex}\label{ex:curves}
 Assume $X/S$ is locally projective over $S$ and flat, and its geometric
fibers are integral curves.  Given an integer $m$, let
$\IPicm_{X/S}\subset \IPic_{X/S}$ be the set of points representing
invertible sheaves $\mc L$ of degree $m$.

Show the $\IPicm_{X/S}$ are open and closed subschemes of finite type;
show they are disjoint and cover; and show that forming them commutes
with changing $S$.

Show there is no abuse of notation: the fiber of $\IPicz_{X/S}$ over
$s\in S$ is the connected component of $0\in\IPic_{X_s/k_s}$.  Show
there is no torsion: $\IPicz_{X/S}=\IPict_{X/S}$.  Show each
$\IPicm_{X/S}$ is an fppf-torsor under $\IPicz_{X/S}$; that is, the
latter acts naturally on the former, and the two become isomorphic after
base change by an fppf-covering.

Show the $\IPicm_{X/S}$ are quasi-projective if $X/S$ is projective and
$S$ is Noetherian.
 \end{ex}

\begin{rmk}\label{rmk:curves}
 There is another important case where $\IPicz_{X/S}=\IPict_{X/S}$,
namely, when $X$ is an Abelian $S$-scheme.  Indeed, the equation holds
if it does on each geometric fiber of $X/S$; so we may assume that $S$
is the spectrum of an algebraically closed field.  In this case,
  a modern proof was given by Mumford \cite[Cor.2, p.~178]{Mm70}.
 \end{rmk}

\begin{eg}\label{eg:Pphifin}
 Theorem~\ref{thm:Pphifin} can fail if a geometric fiber of $X/S$ is
reducible.  For example, let $S$ be the spectrum of a field $k$, and let
$X$ be the union of two disjoint lines.  For each pair $a,b\in\bb Z$,
let $\mc L_{a,b}$ be the invertible sheaf that restricts to $\mc O(a)$
on the first line and to $\mc O(b)$ on the second.

 By Riemann's Theorem, $\chi(\mc L_{a,b}(n)) = (a+b)+2n + 1$.  And it is
easy to see that $\Pic_{X/k}$ is the disjoint union of copies of
$\Spec(k)$ indexed by $\bb Z\x \bb Z$; compare with
Exercise~\ref{ex:PE}.  Moreover, $\mc L_{a,b}$ is represented by the
point with index $(a,b)$.  Hence, each set $\IPicp_{X/k}$ is infinite,
and so not of finite type.
 \end{eg}

\begin{rmk}\label{rmk:Pphifin}
 Theorem~\ref{thm:Pphifin} can be modified as follows.  Assume $S$ is
Noetherian.  Assume $X/S$ is projective and flat with integral geometric
fibers of dimension $r$.  Then a subset $\Lambda\subset\IPic_{X/S}$ is
of finite type if, in the Hilbert polynomials $\sum_{i=0}^r
a_i\binom{n+i}i$ of the corresponding invertible sheaves, $a_{r-1}$
remains bounded from above and below, and $a_{r-2}$ remains bounded from
below alone.  Moreover, $\Lambda$ is of finite type if, instead, $\int
\ell h^{r-1}$ remains bounded from above and below, and $\int\ell^2
h^{r-2}$ remains bounded from below alone, where $\ell$ and $h$ are as
usual.

These facts can be derived from Theorem~\ref{thm:Pphifin} by reducing to
the case where the fibers are normal and by showing, in this case using
simple elementary means, that the given bounds imply bounds on all the
$a_i$; see  \cite[Thm.~3.13, p.~641]{Kl71}.

The first fact is essentially equivalent, given
Theorems~\ref{thm:Pphifin} and \ref{th:rrep}, to the following fact.
Assume $S$ is Noetherian, and $X/S$ is projective and flat with 
fibers whose irreducible components have dimension at least 3.  Let $Y$
be a relative effective divisor whose associated sheaf is $\mc O_{\!X}(1)$.
Then the induced map of Picard functors is representable by maps of
finite type.  See  \cite[Thm.~3.8, p.~636]{Kl71}.

The latter fact was first proved directly using the ``equivalence
criterion'' mentioned in Remark~\ref{rmk:RamSam} by Grothendieck
\cite[p.~C-10]{FGA}.
 \end{rmk}

\begin{cor}\label{cor:cc}
  Assume $X/S$ is projective and flat with integral geometric fibers.
Then the connected components of $\IPic_{X/S}$ are open and closed
subschemes of finite type.
 \end{cor}
\begin{proof}
 By construction, $\IPic_{X/S}$ is locally Noetherian.  Hence, its
connected components are open and closed; see the proof of
Lemma~\ref{lem:agps}.  Now, a connected component is always contained in
any open and closed set it meets.  Hence the connected components of
$\IPic_{X/S}$ are of finite type owing to Theorem~\ref{thm:Pphifin}.
 \end{proof}

\begin{rmk}\label{rmk:cc}
 In Corollary~\ref{cor:cc}, $X/S$ must be projective, not simply proper,
nor even projective Zariski locally over $S$.  Indeed, Grothendieck
\cite[Rem.~3.3, p.~232-07]{FGA} gave an example where $\IPic_{X/S}$ has
a connected component that is not of finite type: here $S$ is a curve
with two components that meet in two points, such as the union of a
smooth conic and a line in the plane over an algebraically closed field,
and $X$ is projective over a neighborhood of each component of $S$.
 \end{rmk}

\begin{cor}\label{cor:phin}
 Assume $X/S$ is projective and flat with integral geometric fibers.
For $n\neq0$,  the \tu nth power map 
	$\vf_n\:\IPic_{X/S}\to\IPic_{X/S}$
is of finite type. 
 \end{cor}
\begin{proof}
 Owing to Corollary~\ref{cor:cc}, we need only prove that, given any
connected component $U$ of $\IPic_{X/S}$, the preimage $\vf_n^{-1}U$ is
of finite type too.  Since $\IPic_{X/k}$ represents $\Pic_{(X/S)\et}$ by
Theorem~\ref{th:main}, the inclusion $U\into\IPic_{X/S}$ is defined by
an invertible sheaf $\mc M$ on $X_T$ for some \'etale covering $T\to U$.

Fix $t\in T$, and set $\psi(p,q):=\chi(\mc M_t^p(q))$.  Fix $p,\,q$, and
form the set $T'$ of points $t'$ of $T$ such that $\chi(\mc M_{t'}^p(q))
= \psi(p,q)$.  By  \cite[7.9.4]{EGAIII2}, the set $T'$ is open, and so
is its complement.  Hence their images are open in $U$, and plainly
these images are disjoint.  But $U$ is connected.  Hence $T'=T$. 

Let $\lambda\in\vf_n^{-1}U$.  Represent $\lambda$ by an invertible sheaf
$\mc L$.  Set $\theta(m,q):=\chi(\mc L^m(q))$.  Say $\theta(m,q)=\sum
a_i(m)\binom{q+i}i$ where $a_i(m)$ is a polynomial.  Now,
$\vf_n(\lambda)\in U$.  So $\theta(mn,q)=\psi(m,q)$.  So $a_i(mn)$ is
independent of the choice of $\lambda$ for all $m$.  Hence $a_i(m)$ is
too.  Set $\phi(n):=\theta(1,q)$.  Then $\vf_n^{-1}U\subset
\IPicp_{X/S}$.  Hence Theorem~\ref{thm:Pphifin} implies $\vf_n^{-1}U$ is
of finite type.
 \end{proof}

\begin{rmk}\label{rmk:phin}
 Corollary~\ref{cor:phin} holds in greater generality.  Assume $X/S$ is
proper, and  assume $\IPic_{X/S}$ exists and represents $\Pic_{(X/k)\fppf}$.
Then, for $n\neq0$,  the \tu nth power map 
	$\vf_n\:\IPic_{X/S}\to\IPic_{X/S}$
 is of finite type.  This fact can be derived from the corollary through
a series of reduction steps similar to those used to generalize
Theorem~\ref{th:Ptaufin}; see \cite[Thm.~3.6, p.~635]{Kl71}.
 \end{rmk}

\begin{ex}\label{ex:conv}
 Assume $S$ is Noetherian, and $X/S$ is projective and flat.  Assume
$\IPic_{X/S}$ exists and represents $\Pic_{(X/S)\fppf}$.  Let $\Lambda$
be an arbitrary subset of $\IPic_{X/S}$, and $\Pi$ the corresponding set
of Hilbert polynomials.  Show $\Pi$ is finite if $\Lambda$ is
quasi-compact.  Show $\Pi$ has only one element if $\Lambda$ is
connected.
 \end{ex}

\appendix
 \section{Answers to all the exercises}\label{sc:A}
 The exercises are not meant to be tricky, but are designed to help you
check, solidify, and expand your understanding of the ideas and methods.
So to promote your own mathematical health, try seriously to do each
exercise before you read its answer here.  Note that the answer key is
the same number as the exercise key.

\begin{ans}{Alr}
 Since $A$ is local, $\Pic(T)$ is trivial; so  Definitions~
\ref{df:aPf} and \ref{df:Pfs} yield the first isomorphism,
$\Pic_X(A)\risom\Pic_{X/S}(A)$.

Let $\mc L$ be an invertible sheaf on $X_A$.  Suppose the isomorphism
class of $\mc L$ maps to $0$ in $\Pic_{(X/S)\zar}(A)$.  Then there is a
Zariski covering $T'\to T$ and an isomorphism $v'\:\mc L_{T'}\risom\mc
O_{\!X_{T'}}$.  Now, $T'$ is a disjoint union of Zariski open subschemes
of $T$.  One of them contains the closed point, so is equal to $T$.
Restricting $v'$ yields an isomorphism $\mc L\risom\mc O_{\!X_A}$.  Thus
$\Pic_X(A)\to\Pic_{(X/S)\zar}(A)$ is injective.

Given $\lambda\in\Pic_{(X/S)\zar}(A)$, represent $\lambda$ by an
invertible sheaf $\mc L'$ on $X_{T'}$ for a suitable Zariski covering
$T'\to T$.  Again, $T'$ contains a copy of $T$ as an open and closed
subscheme.  Restricting $\mc L'$ yields an invertible sheaf $\mc L$ on
$X_A$.  Since $\mc L'$ represent $\lambda$, there is an isomorphism
$v''$ between the two pullbacks of $\mc L'$ to $T'\x_T T'$.  Restricting
$v''$ to $T'\x_T T$, or $T'$, yields an isomorphism between $\mc L'$ and
$\mc L_{T'}$.  Hence $\mc L$ too represents $\lambda$.  Thus
$\Pic_X(A)\to\Pic_{(X/S)\zar}(A)$ is surjective too.

Assume $A$ is Artin local with algebraically closed residue field.  Then
every \'etale $A$-algebra $B$ of finite type is a direct product of
Artin local algebras, each isomorphic to $A$, owing to \cite[17.6.2 and
17.6.3]{EGAIV4}.  So if $T'\to T$ is any \'etale covering, then $T'$ is
a disjoint union of open subschemes, each a copy of $T$.  Hence,
reasoning as above, we conclude $\Pic_X(A)\risom\Pic_{(X/S)\et}(A)$.

Assume $A=k$ where $k$ is an algebraically closed field.  Then any fppf
covering $T'\to T$ has a section; indeed, at any closed point of $T'$,
the local ring is essentially of finite type over $k$, and so the residue
field is equal to $k$ by the Hilbert Nullstellensatz.  This point is not
necessarily isolated in $T'$; nevertheless, reasoning essentially as
above, we conclude $\Pic_X(k)\risom\Pic_{(X/S)\fppf}(k)$.
 \end{ans}

\begin{ans}{Pfs}
 The extension $\bb C/\bb R$ is \'etale.  So if the two pullbacks of
$\vf^*\mc O(1)$ to $X_{\bb C\ox_\bb R\bb C}$ are isomorphic, then
$\vf^*\mc{O}(1)$ defines an element $\lambda$ in $\Pic _{(X/\bb
R)\et}(\bb R)$.

Take an indeterminate $z$, and identify $\bb C$ with $\bb
R[z]/(z^2+1)$. Then, by extension of scalars, $\bb C\ox_\bb R\bb C$
becomes identified with $\bb C[z]/(z-1)(z+1)$.  So, by the Chinese
Remainder Theorem, $\bb C\ox_\bb R\bb C$ is isomorphic to the product
$\bb C\x\bb C$.  (Correspondingly, $c\ox1$ and $1\ox c$ are identified
with $(c,c)$ and $(c,\overline c)$ where $\overline c$ is the conjugate
of $c$, but this fact is not needed here.)

Thus $X_{\bb C\ox_\bb R\bb C}$ is isomorphic, over $\bb R$, to the
disjoint union of two copies of $X_\bb C$.  Now, for any field $k$, an
invertible sheaf $\mc L$ on $\I P^1_k$ is determined, up to isomorphism,
by a single integer its Euler characteristic $\chi(\mc L)$ by
\cite[Cor.~6.17, p.~145]{Ha83}.  Hence, the two pullbacks of $\vf^*\mc
O(1)$ to $X_{\bb C\ox_\bb R\bb C}$ are isomorphic.  Thus
$\vf^*\mc{O}(1)$ defines a $\lambda$ in $\Pic _{(X/\bb R)\et}(\bb R)$.

(Similarly, the isomorphism class of $\vf^*\mc{O}(1)$ on $X_{\bb C}$ is
independent of the choice of the isomorphism $\vf\:X_{\bb C}\risom\I
P^1_{\bb C}$.  So $\lambda$ is independent too.  But this fact too is
not needed here.)

Finally, we must show  $\lambda$ is not in the image of
$\Pic_{(X/\bb R)\zar}(\bb R)$.  By way of contradiction, suppose
$\lambda$ is.  Then $\lambda$ arises from an invertible sheaf $\mc L$ on
$X$.  A priori, the pullback $\mc L|X_{\bb C}$ need not be isomorphic to
$\vf^*\mc{O}(1)$.  Rather, these two invertible sheaves need only become
isomorphic after they are pulled back to $X_A$ where $A$ is some \'etale
$\bb C$-algebra.

However, cohomology commutes with flat base change.  So
	$$\dim_\bb R\uH^0(\mc L)=\rank_A\uH^0(\mc L|X_A)
	=\dim_\bb C\uH^0(\vf^*\mc O(1))=2.$$
 Hence, $\mc L$ has a nonzero section.  It defines an exact sequence
	$$0\to\mc O_{\!X}\to\mc L\to\mc O_{\!D}\to0.$$
 Similarly $\uH^1(\mc O_{\!X})=0$.  Hence $\dim_\bb R\uH^0(\mc O_{\!D})=1$.
Therefore, $D$ is an $\bb R$-point of $X$.  But $X$ has no $\bb
R$-point.  Thus $\lambda$ is not in the image of
$\Pic_{(X/\bb R)\zar}(\bb R)$.
 \end{ans}

\begin{ans}{gpts}
 First of all, we have $\Pic_{(X/S)\fppf}(k)=\Pic_{(X_k/k)\fppf}(k)$
essentially by definition, because a map $T'\to T$ of $k$-schemes is an
fppf-covering if and only if it is an fppf-covering when viewed as a map
of $S$-schemes.  And a similar analysis applies to the other three
functors.  Now, $f_k\:X_k\to k$ has a section; indeed, $f_k$ is of
finite type and $k$ is algebraically closed, and so any closed point of
$X$ has residue field $k$ by the Hilbert Nullstellensatz.  Hence, by
Part 2 of Theorem~\ref{th:cmp}, the $k$-points of all four functors are
the same.  Finally, $\Pic_{X/S}(k)=\Pic(X_k)$ because $\Pic(T)$ is
trivial whenever $T$ has only one closed point.

Whether or not $\mc O_{\!S}\risom f_*\mc O_{\!X}$ holds universally, all
four functors have the same geometric points by Exercise~\ref{ex:Alr};
in fact, given an algebraically closed field $k$, the $k$-points of
these functors are just the elements of $\Pic(X_k)$.
 \end{ans}

\begin{ans}{lsys}
 By definition, a section in $\uH^0(X,\mc L)_{\reg}$ corresponds to an
injection $\mc L^{-1}\into\mc O_{\!X}$.  Its image is an ideal $\mc I$
such that $\mc L^{-1}\risom\mc I$.  So $\mc I$ is the ideal of an
effective divisor $D$.  Then $\mc O_{\!X}(-D)=\mc I$.  So $\mc
L^{-1}\risom\mc O_{\!X}(-D)$.  Taking inverses yields $\mc
O_{\!X}(D)\simeq\mc L$.  So $D\in|\mc L|$.  Thus we have a map
$\uH^0(X,\mc L)_{\reg}\to|\mc L|$.

If the section is multiplied by a unit in $\uH^0(X,\mc O_{\!X}^*)$, then the
injection $\mc L^{-1}\into\mc O_{\!X}$ is multiplied by the same unit, so has
the same image $\mc I$; so then $D$ is unaltered.  Conversely, if $D$
arises from a second section, corresponding to a second isomorphism $\mc
L^{-1}\risom\mc I$, then these two isomorphism differ by an automorphism
of $\mc L^{-1}$, which is given by multiplication by a unit in
$\uH^0(X,\mc O_{\!X}^*)$; so then the two sections differ by multiplication
by this unit.  Thus $\uH^0(X,\mc L)_{\reg}\big/\uH^0(X,\mc
O_{\!X}^*)\into|\mc L|$.

Finally, given $D\in|\mc L|$, by definition there exist an isomorphism
$\mc O_{\!X}(D)\simeq\mc L$.  Since $O_{\!X}(-D)$ is the ideal $\mc I$ of $D$,
the inclusion $\mc I\into\mc O_{\!X}$ yields an injection $\mc
L^{-1}\into\mc O_{\!X}$.  The latter corresponds to a section in
$\uH^0(X,\mc L)_{\reg}$, which yields $D$ via the procedure of the first
paragraph. Thus $\uH^0(X,\mc L)_{\reg}\big/\uH^0(X,\mc
O_{\!X}^*)\risom|\mc L|$.
 \end{ans}

\begin{ans}{sum}
 Let $x\in D+E$.  If $x\notin D\cap E$, then $D+E$ is a relative
effective divisor at $x$, as $D+E$ is equal to $D$ or to $E$ on a
neighborhood of $x$.  So suppose $x\in D\cap E$.  Then
Lemma~\ref{lm:ctn} says $X$ is $S$-flat at $x$, and each of $D$ and $E$
is cut out at $x$ by one element that is regular on the fiber $X_s$
through $x$.  Form the product of the two elements.  Plainly, it cuts
out $D+E$ at $x$, and it too is regular on $X_s$.  Hence $D+E$ is a
relative effective divisor at $x$ by 
Lemma~\ref{lm:ctn} again.
 \end{ans}

\begin{ans}{DivC}
 Consider a relative effective divisor $D$ on $X_T/T$.  Each fiber $D_t$
is of dimension $0$.  So its Hilbert polynomial $\chi(\mc O_{D_t}(n))$
is constant.  Its value is $\dim\tu H^0(\mc O_{D_t})$, which is just
the degree of $D_t$.

The assertions are local on $S$; so we may assume $X/S$ is projective.
Then $\Div_{X/S}$ is representable by an open subscheme
$\IDiv_{X/S}\subset \IHilb_{X/S}$ by Theorem~\ref{th:repDiv}.  And
$\IHilb_{X/S}$ is the disjoint union of open and closed subschemes of
finite type $\IHilb_{X/S}^\phi$ that parameterize the subschemes with
Hilbert polynomial $\phi$.  Set
$\IDiv_{X/S}^m:=\IDiv_{X/S}\bigcap\IHilb_{X/S}^m$.  Then the
$\IDiv_{X/S}^m$ have all the desired properties.

In general, whenever $X/S$ is separated, $X$ represents
$\Hilb_{X/S}^1$, and the diagonal subscheme $\Delta\subset X\x X$ is the
universal subscheme.  Indeed, the projection $\Delta\to X$ is an
isomorphism, so $\Delta\in\Hilb_{X/S}^1(X)$.  Now, given any $S$-map
$g:T\to X$, note $(1\x g)^{-1}\Delta=\Gamma_g$ where $\Gamma_g\subset
X\x T$ is the graph subscheme of $g$, because the $T'$-points of both
$(1\x g)^{-1}\Delta$ and $\Gamma_g$ are just the pairs $(gp,p)$ where
$p\:T'\to T$.  So $\Gamma_g\in\Hilb_{X/S}^1(T)$.

Conversely, let $\Gamma\in \Hilb_{X/S}^1(T)$.  So $\Gamma$ is a closed
subscheme of $X\x T$.  The projection $\pi\:\Gamma\to T$ is proper, and
its fibers are finite; hence, it is finite by Chevalley's
Theorem~\cite[4.4.2]{EGAIII1}.  So $\Gamma=\Spec(\pi_*\mc O_\Gamma)$.
Moreover, $\pi_*\mc O_\Gamma$ is locally free, being flat and finitely
generated over $\mc O_{\!T}$.  And forming $\pi_*\mc O_\Gamma$ commutes with
passing to the fibers, so its rank is 1.  Hence $\mc O_{\!T}\risom\pi_*\mc
O_\Gamma$.  Therefore, $\pi$ is an isomorphism.  Hence $\Gamma$ is the
graph of a map $g\:T\to X$.  So, $(1\x g)^{-1}\Delta=\Gamma$ by the
above; also, $g$ is the only map with this property, since a map is
determined by its graph.  Thus $X$ represents $\Hilb_{X/S}^1$, and
$\Delta\subset X\x X$ is the universal subscheme.

In the case at hand, $\Div_{X/S}^1$ is therefore representable by an
open subscheme $U\subset X$ by Theorem~\ref{th:repDiv}.  In fact, its
proof shows $U$ is formed by the points $x\in X$ where the fiber
$\Delta_x$ is a divisor on $X_x$.  Now, $\Delta_x$ is a $k_x$-rational
point for any $x\in X$; so $\Delta_x$ is a divisor if and only if $X_x$
is regular at $\Delta_x$.  Since $X/S$ is flat, $X_x$ is regular at
$\Delta_x$ if and only if $x\in X_0$.  Thus $X_0=\IDiv_{X/S}^1$.

Finally, set $T:=X_0^m$ and let $\Gamma_i\subset X\x T$ be the graph
subscheme of the $i$th projection.  By the above analysis, $\Gamma_i\in
\Div_{X/S}^1(T)$.  Set $\Gamma:=\sum \Gamma_i$.  Then $\Gamma\in
\Div_{X/S}^m(T)$ owing to Exercise~\ref{ex:sum} and to the additivity of
degree.  Plainly $\Gamma$ represents the desired $T$-point of
$\IDiv_{X/S}^m$.
 \end{ans}

\begin{ans}{gc&r} Let $s\in S$.  Let $K$ be the algebraically closure of
$k_s$, and set $A:=\uH^0(X_K,\mc O_{\!X_K})$.  Since $f$ is proper, $A$ is
finite dimensional as a $K$-vector space; so $A$ is an Artin ring.
Since $X_K$ is connected, $A$ is not a product of two nonzero rings by
\cite[7.8.6.1]{EGAIII2}; so $A$ is an Artin local ring.  Since $X_K$ is
reduced, $A$ is reduced; so $A$ is a field, which is a finite extension
of $K$.  Since $K$ is algebraically closed, therefore $A=K$.  Since
cohomology commutes with flat base change, consequently
$k_s\risom\uH^0(X,\mc O_{\!X_s})$.

The isomorphism $k_s\risom\uH^0(X,\mc O_{\!X_s})$ factors through $f_*(\mc
O_{\!X})\ox k_s$:
	$$k_s\to f_*(\mc O_{\!X})\ox k_s\to\uH^0(X_s,\mc O_{\!X_s}).$$ 
 So the second map is a surjection.  Hence this map is an isomorphism by
the implication (iv)$\Rightarrow$(iii) of Subsection~\ref{sb:Q} with
$\mc F:=\mc O_{\!X}$ and $\mc N:=k_s$.  Therefore, the first map is an
isomorphism too.

It follows that $\mc O_{\!S}\to f_*\mc O_{\!X}$ is surjective at $s$.
Indeed, denote its cokernel by $\mc G$.  Since tensor product is right
exact and since $k_s\to f_*(\mc O_{\!X})\ox k_s$ is an isomorphism, $\mc
G\ox k_s=0$.  So by Nakayama's lemma, the stalk $\mc G_s$ vanishes, as
claimed

Let $\mc Q$ be the $\mc O_{\!S}$-module associated to $\mc F:=\mc
O_{\!X}$ as in Subsection~\ref{sb:Q}.  Then $\mc Q$ is free at $s$ by
the implication (iv)$\Rightarrow$(i) of Subsection~\ref{sb:Q}.  And
$\rank\mc Q_s=1$ owing to the isomorphism in (\ref{eq:Q}) with $\mc N:=
k_s$.  But, with $\mc N:=\mc O_{\!S}$, the isomorphism becomes
$\SHom(\mc Q,\mc O_{\!X})\risom f_*\mc O_{\!X}$.  Hence $f_*\mc O_{\!X}$
too is free of rank 1 at $s$.  Therefore, the surjection $\mc O_{\!S}\to
f_*\mc O_{\!X}$ is an isomorphism at $s$.  Since $s$ is arbitrary, $\mc
O_{\!S}\risom f_*\mc O_{\!X}$ everywhere.

Finally, let $T$ be an arbitrary $S$-scheme.  Then $f_T\:X_T\to T$ too
is proper and flat, and its geometric fibers are reduced and connected.
Hence, by what we just proved, $\mc O_{\!T}\risom f_T*\mc O_{\!X_T}$.
 \end{ans}

\begin{ans}{LinSys}
 By Theorem~\ref{th:LinSys}, $L$ represents $\LinSys_{\mc L/X/S}$.  So
by Yoneda's Lemma~\cite[(0,1.1.4), p.~20]{EGAG}, there exists a $W\in
\LinSys_{\mc L/X/S}(L)$ possessing the required universal property.
And $W$ corresponds to the identity map $p\:L\to L$.  The proof
of  Theorem~\ref{th:LinSys} now shows 
 $\mc O_{X_L}(W)=(\mc L|X_L)\ox f_L^*\mc O_L(1)$.
 \end{ans}

\begin{ans}{0sec}
 The structure sheaf $\mc O_X$ defines a section $\sigma\:S\to
\IPic_{X/S}$.  Its image is a subscheme, which is closed if
$\IPic_{X/S}$ is separated, by \cite[Cors.\ (5.1.4), p.~275, and (5.2.4),
p.~278]{EGAG}.  Let $N\subset T$ be the pullback of this subscheme under
the map $\lambda\:T\to \IPic_{X/S}$ defined by $\mc L$.  Then the third
property holds.

Both $\mc L_N$ and $\mc O_X$ define the same map $N\to \IPic_{X/S}$.
So, since $\mc O_{\!S}\risom f_*\mc O_{\!X}$ holds universally, the
Comparison Theorem, Theorem, Theorem~\ref{th:cmp}, implies that there
exists an invertible sheaf $\mc N$ on $N$ such that the first property
holds.

Consider the second property.  Then $\mc L_{T'}\simeq f_{T'}^*{\mc
N'}$. So $\lambda t\:T'\to\IPic_{X/S}$ is also defined by $\mc
O_{X_{T'}}$; hence, $\lambda t$ factors through $\sigma\:S\to
\IPic_{X/S}$.  Therefore, $t\:T'\to T$ factors through $N$.  So, since
the first property holds, $\mc L_{T'}\simeq f_{T'}^*t^*{\mc N}$.  Hence
$\mc N'\simeq t^*\mc N$ by Lemma~(\ref{lm:fff}).  Thus the second
property holds.

Finally, suppose the pair $(N_1,\,\mc N_1)$ also possesses the first
property.  Taking $t$ to be the inclusion of $N_1$ into $T$, we conclude
that $N_1\subset N$ and $\mc N_1\simeq \mc N|N$.  Suppose $(N_1,\,\mc
N_1)$ possess the second property too.  Then, similarly, $N\subset N_1$.
Thus $N=N_1$ and  $\mc N_1\simeq \mc N$, as desired.
 \end{ans}

\begin{ans}{univshf}
 By Yoneda's Lemma~\cite[(0,1.1.4), p.~20]{EGAG}, a universal sheaf $\mc
P$ exists if and only if $\IPic_{X/S}$ represents $\Pic_{X/S}$.  Set
$P:=\IPic_{X/S}$.

Assume $\mc P$ exists.  Then, for any invertible sheaf $\mc N$ on $P$,
plainly $\mc P\ox f_P^*\mc N$ is also a universal sheaf.  Moreover, if
$\mc P'$ is also a universal sheaf, then $\mc P'\simeq \mc P\ox f_P^*\mc
N$ for some invertible sheaf $\mc N$ on $P$ by the definition with
$h:=1_P$.

Assume $\mc O_{\!S}\risom f_*\mc O_{\!X}$ holds universally.  If $\mc
P\ox f_P^*\mc N \simeq \mc P\ox f_P^*\mc N'$ for some invertible sheaves
$\mc N$ and $\mc N'$ on $P$, then $\mc N\simeq \mc N'$ by
Lemma~\ref{lm:fff}.

By Part 2 of Theorem \ref{th:cmp}, if also $f$ has a section, then
$\IPic_{X/S}$ does represent $\Pic_{X/S}$; so then $\mc P$ exists.
Furthermore, the curve $X/\bb R$ of Exercise~\ref{ex:Pfs} provides an
example where no $\mc P$ exists, because $\Pic_{(X/\bb R)\et}$ is
representable by Theorem~\ref{th:main}, but $\Pic_{X/\bb R}$ is not
since the two functors differ.
 \end{ans}

\begin{ans}{bschg}
 Say $\IPic_{X/S}$ represents $\Pic_{(X/S)\et}$.  Now, for any
$S'$-scheme $T$,
	$$\Pic_{(X_{S'}/S')\et}(T)=\Pic_{(X/S)\et}(T),$$
 which holds essentially by definition, since a map of $S'$-schemes is an
\'etale-covering if and only if it is an \'etale-covering when viewed as
a map of $S$-schemes. However,
	$$(\IPic_{X/S}\x_SS')(T)=\IPic_{X/S}(T)$$
 because the structure map $T\to S'$ is fixed.  Since the right-hand
sides of the two displayed equations are equal, so are their left-hand
sides.  Thus $\IPic_{X/S}\x_SS'$ represents $\Pic_{(X_{S'}/S')\et}$.  Of
course, a similar analysis applies when $\IPic_{X/S}$ represents one of
the other relative Picard functors.

An example is provided by the curve $X\subset \I P^2_\bb R$ of
Exercise~\ref{ex:Pfs}.  Indeed, since the functors $\Pic_{X/\bb R}$ and
$\Pic_{(X/\bb R)\et}$ differ, $\Pic_{X/\bb R}$ is not
representable. But $\Pic_{(X/\bb R)\et}$ is representable by the Main
Theorem, \ref{th:main}.  Finally, since $X_\bb C$ has a $\bb C$-point,
all its relative Picard functors are equal by the Comparison Theorem,
\ref{th:cmp}.
 \end{ans}

\begin{ans}{schpts}
 An $\mc L$ on an $X_k$ defines a map $\Spec(k)\to\IPic_{X/S}$; assign
its image to $\mc L$.  Then, given any field $k''$ containing $k$, the
pullback $\mc L|X_{k''}$ is assigned the same scheme point of
$\IPic_{X/S}$.

 Consider an $\mc L'$ on an $X_{k'}$.  If $\mc L$ and $\mc L'$ represent
the same class, then there is a $k''$ containing both $k$ and $k'$ such
that $\mc L|X_{k''}\simeq\mc L'|X_{k''}$; hence, then both $\mc L$ and
$\mc L'$ are assigned the same scheme point of $\IPic_{X/S}$.
Conversely, if $\mc L$ and $\mc L'$ are assigned the same point, take
$k''$ to be any algebraically closed field containing both $k$ and $k'$.
Then $\mc L|X_{k''}$ and $\mc L'|X_{k''}$ define the same map
$\Spec(k'')\to\IPic_{X/S}$.  Hence $\mc L|X_{k''}\simeq\mc L'|X_{k''}$
by Exercise~\ref{ex:Alr} or \ref{ex:gpts}.

Finally, given any scheme point of $\IPic_{X/S}$, let $k$ be the
algebraic closure of its residue field.  Then $\Spec(k)\to\IPic_{X/S}$
is defined by an $\mc L$ on $X_k$ by Exercise~\ref{ex:Alr} or
\ref{ex:gpts}.  So the given point is assigned to $\mc L$.  Thus the
classes of invertible sheaves on the fibers of $X/S$ correspond
bijectively to the scheme points of $\IPic_{X/S}$.
  \end{ans}

\begin{ans}{PQ-Abel}
 An $S$-map $h\:T\to \IDiv_{X/S}$ corresponds to a relative effective
divisor $D$ on $X_T$.  So the composition $\I A_{X/S}h\:T\to P$
corresponds to the invertible sheaf $\mc O_{\!X_T}(D)$.  Hence
$O_{\!X_T}(D)\simeq (1\x \I A_{X/S}h)^*\mc P \ox f_P^*\mc N$ for some
invertible sheaf $\mc N$ on $T$.  Therefore, if $T$ is viewed as a
$P$-scheme via $\I A_{X/S}h$, then $D$ defines a $T$-point $\eta$ of
$\LinSys_{\mc P/X\x P/P}$.  Plainly, the assignment $h\mapsto \eta$ is
functorial in $T$.  Thus if $\IDiv_{X/S}$ is viewed as a $P$-scheme via
$\I A_{X/S}$, then there is a natural map $\Lambda$ from its functor of
points to $\LinSys_{\mc P/X\x P/P}$.

Furthermore, $\Lambda$ is an isomorphism. Indeed, let $T$ be a
$P$-scheme.  A $T$-point $\eta$ of $\LinSys_{\mc P/X\x P/P}$ is given by
a relative effective divisor $D$ on $X_T$ such that $\mc
O_{\!X_T}(D)\simeq \mc P_T\ox f_T^*\mc N$ for some invertible sheaf $\mc
N$ on $T$.  Then $\mc O_{\!X_T}(D)$ and $\mc P_T$ define the same
$S$-map $T\to P$.  But $\mc P_T$ defines the structure map.  And $\mc
O_{\!X_T}(D)$ defines the composition $\I A_{X/S}h$ where $h\:T\to
\IDiv_{X/S}$ is the map defined by $D$.  Thus $\eta=\Lambda(h)$, and $h$
is determined by $\eta$; hence, $\Lambda$ is an isomorphism.

In other words, $\IDiv_{X/S}$ represents $\LinSys_{\mc P/X\x P/P}$.  But
$\I P(\mc Q)$ too represents $\LinSys_{\mc P/X\x P/P}$ by
Theorem~\ref{th:LinSys}.  Therefore, $\I P(\mc Q)=\IDiv_{X/S}$ as
$P$-schemes.
 \end{ans}

\begin{ans}{epi} First, suppose $F\to G$ is a surjection.  Given a map of
\'etale sheaves $\vf\:F\to H$ such that the two maps $F\x_GF\to H$ are
equal, we must show there is one and only one map $G\to H$ such that
$F\to G\to H$ is equal to $\vf$.

Let $\eta\in G(T)$.  By hypothesis, there exist an \'etale covering
$T'\to T$ and an element $\zeta'\in F(T')$ such that $\zeta'$ and $\eta$
have the same image in $G(T')$.  Set $T'':=T'\x_TT'$.  Then the two
images of $\zeta'$ in $F(T'')$ define an element $\zeta''$ of
$(F\x_GF)(T'')$.  Since the two maps $F\x_GF\to H$ are equal, the two
images of $\zeta''$ in $H(T'')$ are equal.  But these two images are
equal to those of $\vf(\zeta')\in H(T')$.  Since $H$ is a sheaf,
therefore $\vf(\zeta')$ is the image of a unique element $\theta\in
H(T)$.

Note $\theta\in H(T)$ is independent of the choice of $T'$ and
$\zeta'\in F(T')$.  Indeed, let $\zeta_1'\in F(T_1')$ be a second
choice.  Arguing as above, we find $\vf(\zeta_1')\in H(T_1')$ and
$\vf(\zeta')\in H(T')$ have the same image in $H(T_1'\x_TT')$.  So
$\zeta_1'$ also leads to $\theta$.

Define a map $G(T)\to H(T)$ by $\eta\mapsto\theta$.  Plainly this map
behaves functorially in $T$.  Thus there is a map of sheaves $G\to H$.
Plainly, $F\to G\to H$ is equal to $\vf\:F\to H$.  Finally, $G\to H$ is
the only such map, since the image of $\eta$ in $H(T)$ is determined
by the image of $\eta$ in $G(T')$, and the latter must map to
$\vf(\zeta')\in H(T')$.  Thus $G$ is the coequalizer of
$F\x_GF\rightrightarrows F$.

Conversely, suppose $G$ is the coequalizer of $F\x_GF\rightrightarrows
F$.  Form the \'etale subsheaf $H\subset G$ associated to the presheaf
whose $T$-points are the images in $G(T)$ of the elements of $F(T)$.
Then the map $F\to G$ factors through $H$.  So the two maps $F\x_GF\to
H$ are equal.  Since $G$ is the coequalizer, there is a map $G\to H$ so
that $F\to G\to H$ is equal to $F\to H$.  Hence $F\to G\to H\into G$ is
equal to $F\to G$.  So $G\to H\into G$ is equal to $1_G$ by uniqueness.
Therefore, $H=G$.  Thus $F\to G$ is a surjection.
  \end{ans}

\begin{ans}{q-proj} Theorem~\ref{th:main} implies each connected
component $Z'$ of $Z$ lies in an increasing union of open
quasi-projective subschemes of $\IPic_{X/S}$.  So $Z'$ lies in one of
them since $Z'$ is quasi-compact.  So $Z'$ is quasi-projective.  But $Z$
has only finitely many components $Z'$.  Therefore, $Z$ is
quasi-projective.
 \end{ans}

\begin{ans}{Abelpr} Set $P:=\IPic_{X/S}$, which exists by
Theorem~\ref{th:main}.  If $\mc P$ exists, then $\I A_{X/S}$
is, by Exercise~\ref{ex:PQ-Abel}, the structure map of the bundle $\I
P(\mc Q)$ where $\mc Q$ denotes the coherent sheaf on $\IPic_{X/S}$
associated to $\mc P$ as in Subsection~\ref{sb:Q}.  In particular, $\I
A_{X/S}$ is projective Zariski locally over $S$.

In general, forming $P$ commutes with extending $S$ by
Exercise~\ref{ex:bschg}.  Similarly, forming $\I A_{X/S}$ does too.  But
a map is proper if it is after an fppf base extension by
\cite[2.7.1(vii)]{EGAIV2}.

However, $f\:X\to S$ is fppf.  Moreover, $f_X\:X\x X\to X$
has a section, namely, the diagonal.  So use $f$ as a base extension.
Then, by Exercises~\ref{ex:gc&r} and \ref{ex:univshf}, a universal sheaf
$\mc P$ exists.  Therefore, $\I A_{X/S}$ is proper by the first case.
 \end{ans}

\begin{ans}{Abel} Let's use the ideas and notation of
Answer~\ref{ex:Abelpr}.  Now, $X_0$ represents $\Div_{X/S}^1$ by
Exercise~\ref{ex:DivC}.  Hence the Abel map $\I A_{X/S}$ induces a
natural map $A\:X_0\to P$, and forming $A$ commutes with extending $S$.
But a map is a closed embedding if it is after an fppf base extension by
\cite[2.7.1(xii)]{EGAIV2}.  So we may assume $\I P(\mc Q) =
\IDiv_{X/S}$.

The function $\lambda\mapsto\deg\mc P_\lambda$ is locally constant.  Let
$W\subset P$ be the open and closed subset where the function's value is
$1$.  Plainly $\I P(\mc Q_W)=\IDiv_{X/S}^1$ owing to the above.
Therefore, $X_0=\I P(\mc Q_W)$, and $A\:X_0\to P$ is equal to the
structure map of $\I P(\mc Q_W)$.  So it remains to show that this
structure map is a closed embedding.

Fix $\lambda\in W$.  Then $\dim_{k_\lambda}(\mc Q\ox k_\lambda)=
\dim_{k_\lambda} \uH^0(X_\lambda,\mc P_\lambda)$.  Suppose
$P_\lambda$ has two independent global sections.  Each defines an
effective divisor of degree 1, which is a $k_\lambda$-rational point
$x_i$.  Since neither section is a multiple of the other, the $x_i$ are
distinct.  Hence the sections generate $P_\lambda$.  So they define a
map $h\:X_\lambda\to\I P^1_{k_\lambda}$ by \cite[4.2.3]{EGAII} or
\cite[Thm.~II, 7.1, p.~150]{Ha83}.  Then $h$ is birational since each
$x_i$ is the scheme-theoretic inverse image of a $k_\lambda$-rational
point of $\I P^1_{k_\lambda}$.  Hence $h$ is an isomorphism.  But, by
hypothesis,  $X_\lambda$ is of arithmetic genus at least 1.  So there is
a contradiction.  Therefore, $\dim_{k_\lambda}(\mc
Q\ox k_\lambda)\le1$.

By Nakayama's lemma, $\mc Q$ can be generated by a single element on a
neighborhood $V\subset W$ of $\lambda$.  So there is a surjection $\mc
O_{\!V}\onto \mc Q_V$.  It defines a closed embedding $\I P(\mc Q_V )\into\I
P(\mc O_{\!V})$.  But the structure map $\I P(\mc O_{\!V})\to V$ is an
isomorphism.  Hence $\I P(\mc Q_V)\to V$ is a closed embedding.  But
$\lambda\in W$ is arbitrary.  So $\I P(\mc Q_W)\to W$ is indeed
a closed embedding.
 \end{ans}

 \begin{ans}{PE}
 Representing $\Pic_{X/S}$ is similar to representing $\Pic_{X'/S'}$ in
Example~\ref{eg:Mumford}, but simpler.  Indeed, On $X\x_S\bb Z_S$, form
an invertible sheaf $\mc P$ by placing $\mc O_{\!X}(n)$ on the $n$th copy of
$X$.  Then it suffices to show this: given any $S$-scheme $T$ and any
invertible sheaf $\mc L$ on  $X_T$, there exist a unique $S$-map
$q\:T\to\bb Z_S$ and some invertible sheaf $\mc N$ on $T$ such that
$(1\x q)^*\mc P=\mc L\ox f_T^*\mc N$.

Plainly, we may assume $T$ is connected.  Then the function
$s\mapsto\chi(X_t,\mc L_t)$ is constant on $T$ by
\cite[7.9.11]{EGAIII2}.  Now, $X_t$ is a projective space of dimension
at least 1 over the residue field $k_t$; so $\mc L_t\simeq \mc
O_{\!X_t}(n)$ for some $n$ by \cite[Prp.~6.4, p.~132, and Cors.~6.16
and 6.17, p.~145]{Ha83}.  Hence $n$ is independent of $t$.

Set $\mc M:=\mc L^{-1}(n)$.  Then $\mc M_t\simeq\mc O_{\!X_t}$ for all $t\in
T$.  Hence $\uH^1(X_t,\mc M_t)=0$ and $\uH^0(X_t,\mc M_t)=k_t$ by
Serre's explicit computation \cite[2.1.12]{EGAIII1}.  Hence
$f_{\smash{T}*}\mc M$ is invertible, and  forming it commutes with
changing the base $T$, owing to the theory in Subsection~\ref{sb:Q}.

Set $\mc N:=f_{\smash{T}*}\mc M$.  Consider the natural map
$u\:f_T^*\mc N\to \mc M$.  Forming $u$ commutes with changing $T$, since
forming $\mc N$ does.  But $u$ is an isomorphism on the fiber over each
$t\in T$.  So $u\ox k_t$ is an isomorphism.  Hence $u$ is surjective by
Nakayama's lemma. But both source and target of $u$ are invertible; so
$u$ is an isomorphism. Hence $\mc L\ox f_T^*\mc N=\mc O_{\!X_T}(n)$.

Let $q\:T\to\bb Z_S$ be the composition of the structure map $T\to S$
and the $n$th inclusion $S\into\bb Z_S$.  Plainly $(1\x q)^*\mc P=\mc
O_{\!X_T}(n)$, and $q$ is the only such $S$-map.  Thus $\bb Z_S$
represents $\Pic_{X/S}$, and $\mc P$ is a universal sheaf.
 \end{ans}

\begin{ans}{PfsCtd}
 First of all, $\IPic_{X/\bb R}$ exists by Theorem~\ref{th:main}.  Now,
$X_{\bb C}\simeq\I P^1_{\bb C}$.  Hence $\IPic_{X/\bb
R}\x_\bb R\bb C\simeq\bb Z_\bb C$ by Exercises~\ref{ex:bschg} and
\ref{ex:PE}.  The induced automorphism of $\bb Z_{\bb C\ox_\bb R\bb C}$
is the identity; indeed, a point of this scheme corresponds to an
invertible sheaf on $\I P^1_{\bb C}$, and every such sheaf is isomorphic
to its pullback under any $\bb R$-automorphism of  $\I P^1_{\bb C}$.
Hence, by descent theory, $\IPic_{X/\bb R}= \bb Z_\bb R$.

The above reasoning leads to a second proof $\Pic _{(X/\bb R)\et}$ is
representable.  Indeed, set $P:=\Pic _{(X/\bb R)\et}$.  By the above,
the pair $(P\ox_\bb R\bb C)\ox_\bb C(\bb C\ox_\bb R\bb C)
\rightrightarrows P\ox_\bb R\bb C$ is representable by the pair $\bb
Z_{\bb C\ox_\bb R\bb C} \rightrightarrows \bb Z_\bb C$, whose
coequalizer is $\bb Z_\bb R$.  On the other hand, in the category of
\'etale sheaves, the coequalizer is $P$ owing to Exercise~\ref{ex:epi}.

Notice in passing that $\IPic_{X/\bb R}= \IPic_{\I P^1_{\bb R}/\bb R}$.
However, $\Pic_{X/\bb R}$ is not representable owing to
Exercise~\ref{ex:Pfs}, where as $\Pic_{\I P^1_{\bb R}/\bb R}$ is
representable owing to Exercise~\ref{ex:PE}.
 \end{ans}

\begin{ans}{sm=>pr}
 Exercise~\ref{ex:q-proj} implies $Z$ is quasi-projective.  Hence $Z$
is projective if $Z$ is proper.  By \cite[2.7.1]{EGAIV2}, an $S$-scheme
is proper if it is so after an fppf base change, such as $f\:X\to S$.
But $f_X\:X\x X\to X$ has a section, namely, the diagonal.  Thus we may
assume $f$ has a section.

Using the Valuative Criterion for Properness \cite[Thm.~4.7,
p.~101]{Ha83}, we need only check this statement: given an $S$-scheme $T$
of the form $T=\Spec(A)$ where $A$ is a valuation ring, say with
fraction field $K$, every $S$-map $u\:\Spec(K)\to Z$ extends to an
$S$-map $T\to Z$.  We do not need to check the extension is unique if it
exists; indeed, this uniqueness holds by the Valuative Criterion for
Separatedness \cite[Thm.~4.3, p.~97]{Ha83} since $Z$ is quasi-projective,
so separated.

Since $f$ has a section, $u$ arises from an invertible sheaf $\mc L$ on
$X_K$ by Theorem~\ref{th:cmp}.  We have to extend $\mc L$ over $X_T$.
Indeed, this extension defines a map $t\:T\to \IPic_{X/S}$ extending
$u$, and $t$ factors through $Z$ because $Z$ is closed and $T$ is
integral.

Plainly it suffices to extend $\mc L(n)$ for any $n\gg0$.  So replacing
$\mc L$ if need be, we may assume $\mc L$ has a nonzero section.  It is
regular since $X_K$ is integral.  So $X_K$ has a divisor $D$ such that
$\mc O(D)=\mc L$.

  Let $D'\subset X_T$ be the closure of $D$.  Now, $X/S$ is smooth and
$T$ is regular, so $X_T$ is regular by \cite[6.5.2]{EGAIV2}, so
factorial by \cite[21.11.1]{EGAIV2}.  Hence $D'$ is a divisor.  And $\mc
O(D')$ extends $\mc L$.
 \end{ans}

\begin{ans}{q=0}
 Serre's Theorem \cite[Thm.~5.2, p.~228]{Ha83} yields $\tu
H^i(\Omega^2_{\!X}(n))=0$ for $i>0$ and $n\gg0$.  So $\phi(n)=\chi(
\Omega^2_{\!X}(n))$.  Hence
	$$q=\tu H^1(\Omega^2_{\!X})-\tu H^2(\Omega^2_{\!X})+1.$$
 Serre duality \cite[Cor.~7.13, p.~247]{Ha83} yields $\dim\tu
H^i(\Omega^i_{\!X})=\dim\tu H^{2-i}(\mc O_{\!X})$ for all $i$.  And
$\dim\tu H^0(\mc O_{\!X})=1$ since $X$ is projective and geometrically
integral.  So
	$$q=\dim\tu H^1(\mc O_X).$$
 Hence Corollary~\ref{cor:ch0} yields $\dim \IPic_{X/S}\le q$, with
equality in characteristic 0.
 \end{ans}

\begin{ans}{Enriques}
 Set $P:=\IPic_{X/S}$, which exists by Theorem~\ref{th:main}.  By
Exercises \ref{ex:gc&r} and \ref{ex:univshf}, there exists a universal
sheaf $\mc P$ on $X\x P$

Suppose $q=0$.  Then $P$ is smooth of dimension 0 everywhere by
Corollary~\ref{cor:sm}.  Let $D$ be a relative effective divisor on
$X_T/T$ where $T$ is a connected $S$-scheme.  Then $\mc O_{\!X_T}(D)$
defines a map $\tau:T\to P$, and
  $$\mc O_{\!X_T}(D)\simeq (1\x\tau)^*\mc P\ox f_T^*\mc N$$
 for some invertible sheaf $\mc N$ on $T$.  Now, $T$ is connected and
$P$ is discrete and reduced; so $\tau$ is constant.  Set $\lambda:=\tau
T$, and view $\mc P_\lambda$ as an invertible sheaf $\mc L$ on $X$.
Then $\mc L_T=(1\x\tau)^*\mc P$.  So $\mc O_{\!X_T}(D)\simeq \mc L_T\ox
f_T^*\mc N$, as required.

Consider the converse.  Again by Exercise~\ref{ex:univshf}, there is a
coherent sheaf $\mc Q$ on $P$ such that $\I P(\mc Q)=\IDiv_{X/S}$.
Furthermore, $\mc Q$ is nonzero and locally free at any closed point
$\lambda$ representing an invertible sheaf $\mc L$ on $X$ such that $\tu
H^1(\mc L)=0$ by Subsection~\ref{sb:Q}; for example, take $\mc L:=\mc
O_{\!X}(n)$ for $n\gg0$.

Let $U\subset P$ be a connected open neighborhood of $\lambda$ on which
$Q$ is free.  Let $T\subset \I P(\mc Q)$ be the preimage of $U$, and let
$D$ be the universal relative effective divisor on $X_T/T$.  Then the
natural map $A\:T\to U$ is smooth with irreducible fibers.  So $T$ is
connected.  Moreover, $A$ is the map defined by $\mc O_{\!X_T}(D)$.

Suppose $\mc O_{\!X_T}(D)\simeq \mc M_T\ox f_T^*\mc N$ for some
invertible sheaves $\mc M$ on $X$ and $\mc N$ on $T$.  Then $A\:T\to U$
is also defined by $\mc M_T$.  Say $\mu\in P$ represents $\mc M$.  Then
$A$ factors through the inclusion of the closed point $\mu$.  Hence
$\mu=\lambda$; moreover, since $A$ is smooth and surjective, its image,
the open set $U$, is just the reduced closed point $\lambda$.  Now, there
is an automorphism of $P$ that carries 0 to $\lambda$, namely,
``multiplication'' by $\lambda$.  So  $P$  is smooth of dimension 0 at
0.  Therefore, $q=0$ by Corollary~\ref{cor:sm}.

In characteristic 0, a priori $P$ is smooth by Corollary~\ref{cor:ch0}.
Now, $A\:T\to U$ is smooth.  Hence, $T$ is smooth too.  But the
preceding argument shows that, if the condition holds for this $T$, then
$q=0$, as required.
 \end{ans}

\begin{ans}{jac}
 By hypothesis, $\dim X_s=1$ for $s\in S$; so $\tu H^2(\mc O_{\!X_s})=0$.
Hence the $\IPicz_{X_s/k_s}$ are smooth by Proposition~\ref{prp:H2}, so
of dimension $p_a$ by Proposition~\ref{cor:sm}.  Hence, by
Proposition~\ref{prp:P0}, the $\IPicz_{X_s/k_s}$ form a family of finite
type, whose total space is the open subscheme $\IPicz_{X/S}$ of
$\IPic_{X/S}$.  And $\IPic_{X/S}$ is smooth over $S$ again by
Proposition~\ref{prp:H2}.

Hence $\IPicz_{X/S}$ is quasi-projective by Exercise~\ref{ex:q-proj}.

If $X/S$ is smooth, then $\IPicz_{X/S}$ is projective over $S$ by
Exercise~\ref{ex:sm=>pr}.  Alternatively, use Theorem~\ref{th:qpp&p} and 
Proposition~\ref{prp:P0} again to conclude $\IPicz_{X/S}$ is proper, so
projective since it is quasi-projective.

Conversely, assume $\IPicz_{X/S}$ is proper, and let us prove $X/S$ is
smooth, Since $X/S$ is flat, we need only prove each $X_s$ is smooth.
So we may replace $S$ by the spectrum of the algebraic closure of $k_s$.
If $p_a=0$, then $X$ is smooth, indeed $X=\I P^1$, by \cite[Ex.~1.8(b),
p.~298]{Ha83}.

Suppose $p_a>0$.  Let $X_0$ be the open subscheme where $X$ is smooth.
Then there is a closed embedding $A\:X_0\into\IPic_{X/S}$ by
Exercise~\ref{ex:Abel}.  Its image consists of points $\lambda$
representing invertible sheaves of degree $1$.  Fix a rational point
$\lambda$, and define an automorphism $\beta$ of $\IPic_{X/S}$ by
$\beta(\kappa):=\kappa\lambda^{-1}$.  Then $\beta A$ is a closed
embedding of $X_0$ in $\IPicz_{X/S}$.

By assumption, $\IPicz_{X/S}$ is proper.  So $X_0$ is proper.  Hence
$X_0\into X$ is proper since $X$ is separated.  Hence $X_0$ is closed in
$X$.  But $X_0$ is dense in $X$ since $X$ is integral and the ground
field is algebraically closed.  Hence $X_0=X$; in other words, $X$ is
smooth.
 \end{ans}

\begin{ans}{Altpf}
 As before, by Lemma~\ref{lem:hp}, there is an $m$ such that every $\mc
N(m)$ is generated by its global sections.  So there is a section that
does not vanish at any given associated point of $X$; since these points
are finite in number, if $\sigma$ is a general linear combination of the
corresponding sections, then $\sigma$ vanishes at no associated point.
So $\sigma$ is regular, whence defines an effective divisor $D$ such
that $\mc O_{\!X}(-D)=\mc N^{-1}(-m)$.

Plainly $\mc N^{-1}$ is numerically equivalent to $\mc O_{\!X}$ too.  So
$\chi(\mc N^{-1}(n))=\chi(\mc O_{\!X}(n))$ by Lemma~\ref{lem:hp}.  Hence
the sequence $0\to\mc O_{\!X}(-D)\to\mc O_{\!X}\to\mc O_{\!D}\to0$
yields
	$$\chi(\mc O_{\!D}(n))=\psi(n) \text{ where }
	  \psi(n) := \chi(\mc O_{\!X}(n))-\chi(\mc O_{\!X}(n-m)).$$

 Let $T\subset \IDiv_{X/k}$ be the open and closed subscheme
parameterizing the effective divisors with Hilbert polynomial $\psi(n)$.
Then $T$ is a $k$-scheme of finite type.  Let $\mc M'$ be the invertible
sheaf on $X_T$ associated to the universal divisor; set $\mc M:=\mc
M'(-n)$.  Then there exists a rational point $t\in T$ such that $\mc
N=\mc M_t$.  Thus the $\mc N$ numerically equivalent to $\mc O_{\!X}$
form a bounded family.
 \end{ans}

\begin{ans}{h0L}
 Suppose $a<a_r$.  Suppose $\mc L(-1)$ has a nonzero section.  It
defines an effective divisor $D$, possibly 0.  Hence
  $$\textstyle0\le\int h^{r-1}[D]=\int h^{r-1}\ell-\int h^r =a-a_r<0,$$
 which is absurd.  Thus $\tu H^0(\mc L(-1))=0$.

Let $H$ be a hyperplane section of $X$.  Then there is an exact sequence
 $$0\to\mc L(n-1)\to\mc L(n)\to \mc L_H(n)\to0.$$
 It yields the following bound:
 \begin{equation}\label{equation:bd}
 \dim \tu H^0(\mc L(n)) -\dim \tu H^0(\mc L(n-1)) \le
 \dim \tu H^0(\mc L_H(n)).
 \end{equation}
 Since $\binom{n+i}i-\binom{n-1+i}i=\binom{n+i-1}{i-1}$, the sequence
also yields the following formula:
 $$\chi(\mc L_H(n))=\textstyle\sum_{0\le i\le r-1}a_{i+1}\binom{n+i}i.$$
 
Suppose $r=1$.  Then $\dim \tu H^0(\mc L_H(n))=\chi(\mc L_H(n))=a_1$.
Therefore, owing to Equation~(\ref{equation:bd}), induction on $n$ yields
$\dim \tu H^0(\mc L(n))\le a_1(n+1)$, as desired.

Furthermore, $\mc L_H$ is $0$-regular.  Set $m:=\dim \tu H^1(\mc
L(-1))$.  Then $\mc L$ is $m$-regular by Mumford's conclusion at the
bottom of \cite[p.~102]{Mm66}.  But
	$$m=\dim \tu H^0(\mc L(-1))-\chi(\mc L(-1))
	=0-a_1(-1+1)-a_0=-a_0.$$
 Thus we may take $\Phi_1(u_0):=-u_0$ where $u_0$ is an indeterminate.

Suppose $r\ge2$.  Then we may take $H$  irreducible by Bertini's
Theorem \cite[Thm.~12, p.~374]{Se50} or \cite[Cor.~6.7, p.~80]{Jo79}.
Set $h_1:=c_1O_{\!H}(1)$ and  $\ell_1:=c_1\mc L_H$.  Then
  $\int\ell_1 h_1^{r-2}=\int \ell h^{r-2}[H]=a<a_r$.
 So by induction on $r$, we may assume
 $$\dim \tu H^0(\mc L_H(n))\le a_r\tbinom{n+r-1}{r-1}.$$
  Therefore, owing to Equation~(\ref{equation:bd}), induction on $n$
yields the desired bound.

Furthermore, we may assume $\mc L_H$ is $m_1$-regular where
$m_1:=\Phi_{r-1}(a_1\dotsc,a_{r-1})$.  Set $m:=m_1+\dim \tu H^1(\mc
L(m_1-1))$.   By Mumford's same work, $\mc L$ is $m$-regular.
But
 \begin{align*}
 m&=m_1+\dim \tu H^0(\mc L(m_1-1))-\chi(\mc L(m_1-1))\\ &\le\textstyle
    m_1+a_r\binom{m_1-1+r}r-\sum_{0\le i\le r}a_i\binom{m_1-1+i}i.
 \end{align*}
 The latter expression is a polynomial in $a_0,\dotsc,a_{r-1}$ and
$m_1$.  So it is a polynomial $\Phi_r$ in $a_0,\dotsc,a_{r-1}$ alone, as
desired.

In general, consider $\mc N:=\mc L(-a)$.  Then
	$$\chi(\mc N(n))=\textstyle\sum_{0\le i\le r}b_i\binom{n+i}i
  \text{ where } b_i:=\sum_{j=0}^{r-i}a_{i+j}(-1)^j\binom{a-i-j}j.$$
 Set  $\nu:=c_1\mc N$ and  $b:=\int \nu h^{r-1}$.  Then 
$$\textstyle b = \int \ell h^{r-1}-a\int h^r=a-aa_r\le0<a_r.$$
 Hence $\mc N$ is $m$-regular where $m:=\Phi_r(b_0,\dotsc,b_{r-1})$.
But the $b_i$ are polynomials in $a_0,\dotsc,a_r$ and $a$.  Hence there
is a polynomial $\Psi_r$ depending only on $r$ such that
$m:=\Psi_r(a_0,\dotsc,a_r;a)$, as desired.
 \end{ans}

\begin{ans}{Gtauk}
 Let $k'$ be the algebraic closure of $k$.  If $H\ox k'\subset G^\tau\ox
k'$, then $H\subset G^\tau$.  But $G^\tau\ox k'=(G\ox k')^\tau$ by
Lemma~\ref{ex:Gtauk}.  Thus we may assume $k=k'$.

Then $H\subset\bigcup_{h\in H(k)}hG^0$.  But $G^0$ is open, so $hG^0$ is
too.  And $H$ is quasi-compact.  So $H$ lies in finitely many 
$hG^0$.  So $G^0(k)$ has finite index in $H(k)G^0(k)$, say $n$.  Then
$h^n\in G^0(k)$ for every $h\in H(k)$.  So $\vf_n(H)\subset G^0$.
Thus $H\subset G^\tau$.
 \end{ans}

\begin{ans}{teq}
 For any $n$, plainly $\mc L^{\ox n}$ corresponds to $\vf_n\lambda$.
And $\mc L^{\ox n}$ is algebraically equivalent to $\mc O_{\!X}$ if and only
if $\vf_n\lambda\in\IPicz_{X/k}$ by Proposition~\ref{prp:algeq}.  So
$\mc L$ is $\tau$-equivalent to $\mc O_{\!X}$ if and only if
$\lambda\in\IPict_{X/k}$ by Definitions~\ref{dfn:numeq}
and~\ref{dfn:Gtau}.
 \end{ans}

\begin{ans}{ptauproj}
 Theorem~\ref{th:main} implies $\IPic_{X/k}$ exists and represents
$\Pic_{(X/S)\et}$.  So $\IPict_{X/k}$ is of finite type by
Proposition~\ref{prp:pictk}.  Hence $\IPict_{X/k}$ is quasi-projective
by Exercise~\ref{ex:q-proj}.

Suppose $X$ is also geometrically normal.  Since $\IPict_{X/k}$ is
quasi-projective, to prove it is projective, it suffices to prove it is
complete.  By Proposition~\ref{prp:pictk}, forming $\IPicz_{X/k}$
commutes with extending $k$.  And by \cite[2.7.1(vii)]{EGAIV2}, a
$k$-scheme is complete if (and only if) it is after extending $k$.  So
assume $k$ is algebraically closed.

As $\lambda$ ranges over the $k$-points of $\IPict_{X/k}$, the cosets
$\lambda \IPicz_{X/k}$ cover $\IPict_{X/k}$.  So finitely many cosets
cover, since $\IPicz_{X/k}$ is an open by Proposition~\ref{prp:pic0} and
since $\IPict_{X/k}$ is quasi-compact, Now, $\IPicz_{X/k}$ is projective
by Theorem~\ref{th:qpp&p}, so complete, And $\IPict_{X/k}$ is closed
again by Proposition~\ref{prp:pictk}.  Therefore, $\IPict_{X/k}$ is
complete.
 \end{ans}

\begin{ans}{bddsbset}
 Suppose $L$ is bounded.  Then $\mc M$ defines a map $\theta\:T\to
\IPic_{X/S}$, and $\theta(T) \supset \Lambda$.  Since $T$ is Noetherian,
plainly so is $\theta(T)$; whence, plainly so is any subspace of
$\theta(T)$.  Thus $\Lambda$ is quasi-compact.

Conversely, suppose $\Lambda$ is quasi-compact.  Since $\IPic_{X/S}$ is
locally of finite type by Proposition~\ref{prp:lft}, there is an open
subscheme of finite type containing any given point of $\Lambda$.  So
finitely many of the subschemes cover $\Lambda$.  Denote their union by
$U$.

The inclusion $U\into \IPic_{X/S}$ is defined by an invertible sheaf
$\mc M$ on $X_T$ for some fppf covering $T\to U$.  Replace $T$ be an
open subscheme so that $T\to U$ is of finite type and surjective.  Since
$U$ is of finite type, so is $T$.  Given $\lambda\in\Lambda$, let $t\in
T$ map to $\lambda$.  Then $\lambda$ corresponds to the class of $\mc
M_t$.
 \end{ans}

\begin{ans}{sm=>Pictpr}
 Theorem~\ref{th:Ptaufin} asserts $\IPict_{X/S}$ is of finite type.  So
it is projective by Exercise~\ref{ex:sm=>pr}. 
 \end{ans}

\begin{ans}{curves}
 Plainly, replacing $S$ by an open subset, we may assume $X/S$ is
projective and $S$ is connected.  Given $s\in S$, set $\psi(n):=\chi(\mc
O_{\!X_s}(n))$.  Then $\psi(n)$ is independent of $s$.  Given $m$, set
$\phi(n):=m+\psi(n)$.

Let $\lambda\in\IPic_{X/S}$.  Then $\lambda\in\IPicm_{X/S}$ if and only
if $\lambda$ represents an invertible sheaf $\mc L$ of degree $m$.  And
$\lambda\in\IPicp_{X/S}$ if and only if $\chi(\mc L(n))=\phi(n)$.  But,
	$$\chi(\mc L(n))=\deg(\mc L(n))+\psi(0)=\deg(\mc L)+\psi(n)$$
 by Riemann's Theorem and the additivity of $\deg(\bullet)$. Hence
$\IPicm_{X/S}=\IPicp_{X/S}$.  So Theorem~\ref{thm:Pphifin} yields all
the assertions, except for the two middle about $\IPicz_{X/S}$.

To show $\IPicz_{X/S}=\IPict_{X/S}$, similarly we need only show
$\deg\mc L=0$ if and only if $\mc L$ is $\tau$-equivalent to $\mc O_{\!X}$,
for, by Exercise~\ref{ex:teq}, the latter holds if and only if
$\lambda\in\IPict_{X/S}$.  Plainly, we may assume $\mc L$ lives on a
geometric fiber of $X/S$.  Then the two conditions on $\mc L$ are
equivalent by Theorem ~\ref{thm:numeq}.

Since $\deg$ is additive, multiplication carries $\IPicz_{X/S}\x
\IPicm_{X/S}$ set-theoretically into $\IPicm_{X/S}$.  So $\IPicz_{X/S}$
acts on $\IPicm_{X/S}$ since these two sets are open in $\IPic_{X/S}$.

Since $X/S$ is flat with integral geometric fibers, its smooth locus
$X_0$ provides an fppf covering of $S$.  Temporarily, make the base
change $X_0\to S$.  After it, the new map $X_0\to S$ has a section.  Its
image is a relative effective divisor $D$, and tensoring with $\mc
O_{\!X}(mD)$ defines the desired isomorphism from $\IPicz_{X/S}$ to
$\IPicm_{X/S}$.

Finally, to show there is no abuse of notation, we must show the fiber
$(\IPicz_{X/S})_s$ is connected.  To do so, we instead make the base
change to the spectrum of an algebraically closed field $k\supset k_s$.
Then $X_0$ has a $k$-rational point $D$, and again tensoring with $\mc
O_{\!X}(mD)$ defines an isomorphism from $\IPicz_{X/k}$ to
$\IPicm_{X/k}$.  So it suffices to show $\IPicm_{X/k}$ is connected for
some $m\ge1$.

Let $\beta\:X_0^m\to\IDiv_{X/k}^m\to\IPicm_{X/S}$ be composition of the
map $\alpha$ of Exercise~\ref{ex:DivC} and the Abel map.  Since $X$ is
integral, so is the $m$-fold product $X_0^m$.  Hence it suffices to show
$\beta$ is surjective for some $m\ge1$.

By Exercise~\ref{ex:h0L}, there is an $m_0\ge1$ such that every
invertible sheaf on $X$ of degree $0$ is $m_0$-regular.  Set
$m:=\deg(\mc O_{\!X}(m_0))$.  Then every invertible sheaf $\mc L$ on $X$ of
degree $m$ is $0$-regular, so generated by its global sections.

In particular, for each singular point of $X$, there is a global section
that does not vanish at it.  So, since $k$ is infinite, a general linear
combination of these sections vanishes at no singular point.  This
combination defines an effective divisor $E$ such that $\mc O_{\!X}(E)=\mc
L$.  It follows that $\beta$ is surjective, as desired.
 \end{ans}

\begin{ans}{conv}
 Suppose $\Lambda$ is quasi-compact.  Then, owing to
Exercise~\ref{ex:bddsbset}, there exist an $S$-scheme $T$ of finite type
and an invertible sheaf $\mc M$ on $X_T$ such that every polynomial
$\phi\in\Pi$ is of the form $\phi(n)=\chi(\mc M_t(n))$ for some $t\in
T$.  Hence, by  \cite[7.9.4]{EGAIII2}, the number of $\phi$ is at most
the number of connected components of $T$.  Thus $\Pi$ is finite.

Suppose $\Lambda$ is connected.  Then its closure is too.  So we may
assume $\Lambda$ is closed.  Give $\Lambda$ its reduced subscheme
structure.  Then the inclusion $\Lambda\into\IPic_{X/S}$ is defined by
an invertible sheaf $\mc M$ on $X_T$ for some fppf covering $T\to
\Lambda$.  Fix $t\in T$ and set $\phi(n):=\chi(\mc M_t(n))$.  Fix $n$,
and form the set $T'$ of points $t'$ of $T$ such that $\chi(\mc
M_{t'}(n)) = \phi(n)$.  By \cite[7.9.4]{EGAIII2}, the set $T'$ is open,
and so is its complement.  Hence their images are open in $\Lambda$, and
plainly these images are disjoint.  But $\Lambda$ is connected.  Hence
$T'=T$. Thus $\Pi=\{\phi\}$.
 \end{ans}

\section{Basic intersection theory}\label{sc:B}
 This appendix contains an elementary treatment of basic intersection
theory, which is more than sufficient for many purposes, including the
needs of Section~\ref{sc:Pictau}.  The approach was originated in
1959--60 by Snapper.  His results were generalized and his proofs were
simplified immediately afterward by Cartier \cite{Ca60}.  Their work was
developed further in fits and starts by the author.
 
The index theorem was proved by Hodge in 1937.  Immediately afterward,
B. Segre~\cite[\S~1]{Se37} gave an algebraic proof for surfaces, and
this proof was rediscovered by Grothendieck in 1958.  Their work was
generalized a tad in \cite[p.~662]{Kl71}, and a variation appears below
in Theorem~\ref{thm:HIT}.  From the index theorem,
Segre~\cite[\S6]{Se37} derived a connectedness statement like
Corollary~\ref{cor:conn} for surfaces, and the proof below is basically
his.

\begin{dfn}\label{dfn:K(X)}
 Let $\I F(X/S)$ or $\I F$ denote the Abelian category of coherent
sheaves $\mc F$ on $X$ whose support $\Supp\mc F$ is proper over an
Artin subscheme of $S$, that is, a 0-dimensional Noetherian closed
subscheme.  For each $r\ge0$, let $\I F_r$ denote the full subcategory
of those $\mc F$ such that $\dim\Supp\mc F\le r$.

Let $\I K(X/S)$ or $\I K$ denote the ``Grothendieck group'' of $\I F$,
namely, the free Abelian group on the $\mc F$, modulo short exact
sequences.  Abusing notation, let $\mc F$ also denote its class.  And if
$\mc F=\mc O_{\!Y}$ where $Y\subset X$ is a subscheme, then let $[Y]$ also
denote the class.  Let $\I K_r$ denote the subgroup generated by $\I
F_r$.

Let $\chi\:\I K\to\bb Z$ denote the homomorphism induced by the Euler
characteristic, which is just the alternating sum of the lengths of the
cohomology groups.

Given $\mc L\in \Pic(X)$, let $c_1(\mc L)$ denote the endomorphism of
$\I K$ defined by the following formula:
	$$c_1(\mc L)\mc F:=\mc F-\mc L^{-1}\ox\mc F.$$
 Note that  $c_1(\mc L)$ is well defined since tensoring with  $\mc L^{-1}$
preserves exact sequences.
 \end{dfn}

\begin{lem}\label{lem:div}
 Let $\mc L\in \Pic(X)$.  Let $Y\subset X$ be a closed subscheme with
$\mc O_{\!Y}\in\I F$.  Let $D\subset Y$ be an effective divisor such that
$\mc O_{\!Y}(D)\simeq\mc L_Y$.  Then
 	$$c_1(\mc L)\cdot[Y]=[D].$$
 \end{lem}
\begin{proof}
 The left side is defined since $\mc O_{\!Y}\in\I F$.  The equation
results from the sequence $0\to\mc O_{\!Y}(-D)\to\mc O_{\!Y}\to \mc
O_{\!D}\to 0$ since $\mc O_{\!Y}(-D) \simeq \mc L^{-1}\ox\mc O_{\!Y}$.
 \end{proof}

\begin{lem}\label{lem:add}
 Let $\mc L,\,\mc M\in \Pic(X)$.  Then the following relations hold:
\begin{gather*}
 c_1(\mc L)c_1(\mc M) = c_1(\mc L)+c_1(\mc M)-c_1(\mc L\ox\mc M);\\
 c_1(\mc L)c_1(\mc L^{-1}) = c_1(\mc L) + c_1(\mc L^{-1});\\
 c_1(\mc O_{\!X})=0.
 \end{gather*} 
Furthermore,  $c_1(\mc L)$ and $c_1(\mc M)$ commute.
 \end{lem}\begin{proof}
 Let $\mc F\in\I K$.  By definition, $c_1(\mc O_{\!X})\mc F=0$; thus the
third relation holds.
Plainly, each side of the first relation carries
 $\mc F$  into
	$$\mc F-\mc L^{-1}\ox\mc F-\mc M^{-1}\ox\mc F
	+\mc L^{-1}\ox\mc M^{-1}\ox\mc F.$$
 Thus the first relation holds.  It and the third relation imply the
second.  Furthermore, the first relation implies that $c_1(\mc L)$ and
$c_1(\mc M)$ commute.
 \end{proof}

\begin{lem}\label{lem:cyc}
  Given $\mc F\in\I F_r$,  let $Y_1$, \dots, $Y_s$ be the $r$-dimensional
irreducible components of\/ $\Supp\mc F$ equipped with their induced
reduced structure, and let $l_i$ be the length of the stalk of $\mc F$
at the generic point of $Y_i$.  Then, in $\I K_r$,
         $$\mc F \equiv \sum l_i\cdot [Y_i] \bmod \I K_{r-1}.$$
 \end{lem}
\begin{proof}
 The assertion holds if it does after we replace $S$ by a neighborhood of
the image of\/ $\Supp\mc F$.  So we may assume $S$ is Noetherian.

Let $\I F'\subset\I F_r$ denote the family of $\mc F$ for which the
assertion holds.  Since length$(\bullet)$ is an additive function, $\I
F'$ is ``exact'' in the following sense: for any short exact sequence
$0\to\mc F'\to\mc F\to\mc F''\to0$ such that two of the $\mc F$s belong
to $\I F'$, then the third does too.  Trivially, $\mc O_{\!Y}\in\I F'$
for any closed integral subscheme $Y\subset X$ such that $\mc
O_{\!Y}\in\I F_r$.  Hence $\I F'=\I F_r$ by the ``Lemma of
D\'evissage,'' \cite[Thm.~3.1.2]{EGAIII1}.
 \end{proof}

\begin{lem}\label{lem:deg}
 Let $\mc L\in \Pic(X)$.  Then $c_1(\mc L)\I K_r\subset\I K_{r-1}$ for
all $r$.
 \end{lem}
\begin{proof}
 Let $\mc F\in\I F_r$.  Then $\mc F$ and $\mc L^{-1}\ox\mc F$ are
isomorphic at the generic point of each component of\/
$\Supp\mc F$.  So Lemma~\ref{lem:cyc} implies $c_1(\mc L)\mc F\in\I
K_{r-1}$.
 \end{proof}

\begin{lem}\label{lem:AA}
 Let $\mc L\in \Pic(X)$, let $\mc F\in\I K_r$, and let $m\in\bb Z$.
Then
 $$\mc L^{\ox m}\ox\mc F=\textstyle
 \sum_{i=0}^r\binom{m+i-1}ic_1(\mc L)^i\mc F.$$
 \end{lem}
\begin{proof}
 Let $x$ be an indeterminate, and consider the formal identity
	    $$(1-x)^n=\textstyle\sum_{i\ge0}(-1)^i \binom ni x^i.$$
 Replace $x$ by $1-y^{-1}$, set $n:=-m$, and use the familiar identity
		$$\textstyle(-1)^i\binom ni=\binom{m+i-1}i,$$
 to obtain the formal identity
	       $$y^m=\textstyle\sum\binom{m+i-1}i(1-y^{-1})^i.$$
 It yields the assertion, because
$c_1(L)^i\mc F=0$ for $i>r$ owing to Lemma~\ref{lem:deg}.
 \end{proof}

\begin{thm}[Snapper]\label{thm:Snapper}
 Let $\mc L_1, \dotsc, \mc L_n\in \Pic(X)$, let $m_1,\dotsc, m_n\in\bb
Z$, and let $\mc F\in\I K_r$.  Then the Euler characteristic
$\chi\bigl(\mc L_1^{\ox m_1}\ox\cdots\ox\mc L_n^{\ox m_n}\ox\mc
F\bigr)$ is given by a polynomial in the $m_i$ of degree at most $r$.
In fact,
 $$\chi\bigl(\mc L_1^{\ox m_1}\ox\dotsm\ox\mc  L_n^{\ox m_n}\ox\mc F\bigr)
 \textstyle=\sum a(i_1,\dotsc,i_n)\binom{m_1+i_1-1}{i_1}\cdots
	 \binom {m_n+i_n-1}{i_n}$$
 where $i_j\ge0$ and   $\sum i_j\le r$ and where
   $a(i_1,\dots,i_n):=\chi\bigl(c_1(\mc L_1)^{i_1}\dotsm c_1(\mc
L_n)^{i_n}\mc F\bigr)$.
 \end{thm}
\begin{proof}
 The theorem follows from Lemmas~\ref{lem:AA} and \ref{lem:deg}.
 \end{proof}

\begin{dfn}\label{dfn:intnos}
  Let $\mc L_1, \dotsc, \mc L_r\in \Pic(X)$, repetitions allowed.  Let
$\mc F\in\I K_r$.  Define the {\it intersection
number\/} or {\it intersection symbol\/} by the formula
 $$\textstyle\int c_1(\mc L_1)\dotsm c_1(\mc L_r)\mc F
 :=\chi\bigl(c_1(\mc L_1)\dotsm c_1(\mc L_r)\mc F\bigr)
 \in\bb Z.$$
 If $\mc F=\mc O_{\!X}$, then also write  $\int c_1(\mc
L_1)\dotsm c_1(\mc L_r)$ for the number.
 If $\mc L_j=\mc O_{\!X}(D_j)$ for a divisor $D_j$, then also write
$(D_1\dotsm D_r\cdot \mc F)$, or just $(D_1\dotsm D_r)$
if $\mc F=\mc O_{\!X}$.
 \end{dfn}

\begin{thm}\label{thm:1stprpts}
 Let $\mc L_1, \dotsc, \mc L_r\in \Pic(X)$ and $\mc F\in\I K_r$.

\tu{(1)} If $\mc F\in \I F_{r-1}$, then $\int c_1(\mc L_1)\dotsm
c_1(\mc L_r) \mc F = 0$.

\tu{(2) (symmetry and additivity)} The symbol $\int c_1(\mc L_1)\dotsm
c_1(\mc L_r) \mc F$ is symmetric in the $\mc L_j$.  Furthermore, it is a
homomorphism separately in each $\mc L_j$ and in $\mc F$.

\tu{(3)} Set $\mc E:=\mc L_1^{-1}\oplus\dotsb\oplus\mc L_r^{-1}$.  Then
  $$\textstyle\int c_1(\mc L_1)\dotsm c_1(\mc L_r)\mc F =
\sum_{i=0}^r(-1)^i\chi\bigl(\bigl(\bigwedge^i\mc E\bigr)\ox\mc F\bigr).$$
 \end{thm} \begin{proof}
 Part (1) results from Lemma~\ref{lem:deg}.  So the symbol is a
homomorphism in each $\mc L_j$ owing to the relations asserted in
Lemma~\ref{lem:add}.  Furthermore, the symbol is symmetric owing to the
commutativity asserted in Lemma~\ref{lem:add}.  Part (3) results from
the definitions.
 \end{proof}

\begin{cor}\label{cor:r=2}
 Let $\mc L_1,\, \mc L_2\in \Pic(X)$ and $\mc F\in\I K_2$.  Then
  $$\textstyle\int c_1(\mc L_1)c_1(\mc L_2)\mc F = \chi(\mc F)
 - \chi(\mc L_1^{-1}\ox \mc F) -\chi(\mc L_2^{-1}\ox\mc F)
 +\chi(\mc L_1^{-1}\ox\mc L_2^{-1}\ox\mc F).$$
 \end{cor}\begin{proof}
 The assertion is a special case of Part (3) of
Proposition~\ref{thm:1stprpts}.
 \end{proof}

\begin{cor}\label{cor:CM}
 Let $D_1,\dotsc,D_r$ be effective divisors on $X$, and $\mc F\in\I
 F_r$.  Set $Z:=D_1\cap\dotsm\cap D_r$.  Suppose $Z\cap\Supp F$ is
finite, and at each of its points, $F$ is Cohen--Macaulay.  Then
	$$(D_1\dotsm D_r\cdot \mc F)=\length \tu H^0(\mc F_Z)
	 \text{ where }\mc F_Z:=\mc F\ox\mc O_{\!Z}.$$
  \end{cor}\begin{proof}
 For each $j$, set $\mc L_j:=\mc O_{\!X}(D_j)$ and let $\sigma_j \in \tu
H^0(\mc L_j)$ be the section defining $D_j$.   Set $\mc E:=\mc
L_1^{-1}\oplus\dotsb\oplus\mc L_r^{-1}$.  Form the corresponding Koszul
complex $(\bigwedge^\bullet\mc E\bigr)\ox\mc F$ and its cohomology
sheaves $\mc H^i\bigl(\bigl( \bigwedge^\bullet\mc E\bigr)\ox\mc
F\bigr)$.  Then
  $$\textstyle\int c_1(\mc L_1)\dotsm c_1(\mc L_r)\mc F =
 \sum_{i=0}^r(-1)^i\chi\bigl( \mc H^i\bigl(\bigl(
 \bigwedge^\bullet\mc E\bigr)\ox\mc F\bigr)\bigr).$$
 owing to Part (3) of Proposition~\ref{thm:1stprpts} and to the
additivity of $\chi$.  Furthermore, essentially by definition,
$\mc H^0\bigl(\bigl( \bigwedge^\bullet\mc E\bigr)\ox\mc F\bigr)=\mc
F_Z$.  And by standard  local algebra, the higher $\mc H^i$ vanish.
Thus the assertion holds.
 \end{proof}

\begin{lem}\label{lem:icyc}
 Let $\mc L_1, \dotsc, \mc L_r\in \Pic(X)$ and $\mc F\in\I F_r$.  Let $Y_1$,
\dots, $Y_s$ be the $r$-dimensional irreducible components of\/
$\Supp\mc F$ given their induced reduced structure, and let
$l_i$ be the length of the stalk of $\mc F$ at the generic point of
$Y_i$.  Then
  $$\textstyle\int c_1(\mc L_1)\dotsm c_1(\mc L_r)\mc F =
	\sum_i l_i\int c_1(\mc L_1)\dotsm c_1(\mc L_r)[Y_i].$$
 \end{lem}
\begin{proof}
 Apply Lemma~\ref{lem:cyc} and Parts (1) and (2) of
Proposition~\ref{thm:1stprpts}.
 \end{proof}

\begin{lem}\label{lem:idiv}
 Let $\mc L_1, \dotsc, \mc L_r\in \Pic(X)$ and $Y\subset X$ a closed
subscheme with $\mc O_{\!Y}\in\I F$.  Let $D\subset Y$ be an effective
divisor such that $\mc O_{\!Y}(D)\simeq\mc L_r|Y$.  Then
  $$\textstyle\int c_1(\mc L_1)\dotsm c_1(\mc L_r)[Y] =
	\int c_1(\mc L_1)\dotsm c_1(\mc L_{r-1})[D].$$
 \end{lem}\begin{proof}
 Apply Lemma~\ref{lem:div}.
 \end{proof}

\begin{prp}\label{prp:amp}
 Let $\mc L_1, \dotsc, \mc L_r\in \Pic(X)$ and $\mc F\in\I F_r$.  If all
the $\mc L_j$ are relatively ample and if $\mc F\notin \I K_{r-1}$, then
 $$\textstyle\int c_1(\mc L_1)\dotsm c_1(\mc L_r)\mc F>0.$$
 \end{prp}\begin{proof}
  Proceed by induction on $r$.  If $r=0$, then $\int \mc F=\dim\tu
H^0(\mc F)$ essentially by definition, and $H^0(\mc F)\neq0$ since $\mc
F\notin \I K_{r-1}$ by hypothesis.

Suppose $r\ge1$.  Owing to Proposition~\ref{lem:icyc}, we may assume
$\mc F=\mc O_{\!Y}$ where $Y$ is integral.  Owing to Part (2) of
Theorem~\ref{thm:1stprpts}, we may replace $\mc L_r$ by a multiple, and
so assume it is very ample.  Then, for the corresponding embedding of
$Y$, a hyperplane section $D$ is a nonempty effective divisor such that
$\mc O_{\!Y}(D)\simeq\mc L_r|Y$.  Hence the assertion results from
Proposition~\ref{lem:idiv} and the induction hypothesis.
 \end{proof}

\begin{lem}\label{lem:prjfm}
 Let $g\:X'\to X$ be an $S$-map.  Let $\mc L_1, \dotsc, \mc L_r\in
\Pic(X)$ and let $\mc F\in\I F_r(X'/S)$.  Then
       $$\textstyle\int c_1(g^*\mc L_1)\dotsm c_1(g^*\mc L_r)\mc F=
	    \int c_1(\mc L_1)\dotsm c_1(\mc L_r)g_*\mc F .$$
 \end{lem}\begin{proof}
 Let $\mc G\in\I F_r(X'/S)$.  Then, by hypothesis, $\Supp\mc G$ is
proper over an Artin subscheme of $S$, and $\dim\Supp\mc G\le r$;
furthermore, $X/S$ is separated.  Hence, the restriction $g|\Supp\mc G$
is proper; so $g(\Supp\mc G)$ is closed.  And by the dimension
theory of schemes of finite type over Artin schemes, $\dim g(\Supp\mc G)\le
r$.  Therefore, $\tu R^ig_*\mc G\in\I F_r(X/S)$ for all $i$.

 Define a map $\I F_r(X'/S)\to\I K_r(X/S)$ by $\mc G \mapsto \sum_{i=0}^r
(-1)^i\tu R^ig_*\mc G$.  It induces a homomorphism $\tu Rg_*\:\I
K_r(X'/S)\to\I K_r(X/S)$.  And $\chi(Rg_*(\mc G)) = \chi(\mc G)$ owing to
the Leray Spectral Sequence \cite[{\bf0}-12.2.4]{EGAIII1} and to the
additivity of $\chi$ \cite[{\bf0}-11.10.3]{EGAIII1}.  Furthermore, $\mc
L\ox \tu R^ig_*(\mc G)\risom\tu R^ig_*(g^*\mc L\ox \mc G)$ for any $\mc
L\in\Pic(X)$ by \cite[{\bf0}-12.2.3.1]{EGAIII1}.  Hence $c_1(\mc
L)Rg_*(\mc G) = Rg_*(c_1(g^*\mc L)\mc G)$.  Therefore,
    $$\textstyle\int c_1(g^*\mc L_1)\dotsm c_1(g^*\mc L_r)\mc F=
	  \int c_1(\mc L_1)\dotsm c_1(\mc L_r)\tu Rg_*\mc F.$$

Finally, $\tu Rg_*\mc F \equiv \tu g_*\mc F \bmod\I K_{r-1}(X/S)$,
because $\tu R^ig_*\mc F\in\I F_{r-1}$ for $i\ge1$ since, if $W\subset
X'$ is the locus where $\Supp\mc F\to X$ has fibers of dimension at least
1, then $\dim g(W)\le r-1$.  So Part (1) of Theorem~\ref{thm:1stprpts}
yields the asserted formula.
 \end{proof}

\begin{prp}[Projection Formula]\label{prp:prjfm} Let $g\:X'\to X$ be an
$S$-map.  Let $\mc L_1, \dotsc, \mc L_r\in \Pic(X)$.  Let $Y'\subset X'$
be an integral subscheme with $\mc O_{Y'}\in\I F_r(X'/S)$. Set
$Y:=gY'\subset X$, give $Y$ its induced reduced structure, and let
$\deg(Y'/Y)$ be the degree of the function field extension if finite and
be $0$ if not.  Then
       $$\textstyle\int c_1(g^*\mc L_1)\dotsm c_1(g^*\mc L_r)[Y']=
	 \deg(Y'/Y)\int c_1(\mc L_1)\dotsm c_1(\mc L_r)[Y].$$
 \end{prp}\begin{proof}
  Lemma~\ref{lem:cyc} yields $g_*\mc O_{Y'}\equiv \deg(Y'/Y)[Y]\bmod \I
K_{r-1}(X/S)$.  So the assertion results from Lemma~\ref{lem:prjfm}
and from Part (1) of Theorem ~\ref{thm:1stprpts}.
 \end{proof}

\begin{prp}\label{prp:bschg}
 Assume $S$ is the spectrum of a field, and let $T$ be the spectrum of
an extension field.  Let $\mc L_1, \dotsc, \mc L_r\in
\Pic(X)$ and $\mc F\in\I F_r(X/S)$.  Then
   $$\textstyle\int c_1(\mc L_{1,T})\dotsm c_1(\mc L_{r,T})\mc F_T =
	\int c_1(\mc L_1)\dotsm c_1(\mc L_r)\mc F.$$
 \end{prp}\begin{proof}
 The base change $T\to S$ preserves short exact sequences.  So it
induces a homomorphism $\kappa\:\I K_r(X/S)\to\I K_r(X_T/T)$.  Plainly
$\kappa$ preserves the Euler characteristic.  The assertion now follows.
 \end{proof}

\begin{prp}\label{prp:aleq}
 Let $\mc L_1, \dotsc, \mc L_r\in \Pic(X)$.  Let $\mc F$ be a flat
coherent sheaf on $X$.  Assume  $\Supp\mc F$  is proper and of relative
dimension $r$.  Then the function
   $$y\mapsto\textstyle\int c_1(\mc L_1)\dotsm c_1(\mc L_r)\mc F_y$$
  is locally constant.
 \end{prp}\begin{proof}
 The assertion results from Definition~\ref{dfn:intnos} and
\cite[Cor., top p.~50]{Mm70}.
 \end{proof}

\begin{dfn}\label{dfn:neq}
  Let $\mc L,\,\mc N\in\Pic(X)$.  Call them {\it numerically equivalent}
if $\int c_1(\mc L)[Y]=\int c_1(\mc N)[Y]$ for all closed integral
curves $Y\subset X$ with $\mc O_{\!Y}\in \I F_1$.
 \end{dfn}

\begin{prp}\label{prp:neq}
 Let $\mc L_1,\dotsc,\mc L_r;\,\mc N_1,\dotsc,\mc N_r\in \Pic(X)$ and
$\mc F\in\I K_r$.  If  $\mc L_j$ and $\mc N_j$ are numerically equivalent
for each  $j$, then 
 $$\textstyle\int c_1(\mc L_1)\dotsm c_1(\mc L_r)\mc F =
  \int c_1(\mc N_1)\dotsm c_1(\mc N_r)\mc F.$$
 \end{prp}\begin{proof}
 If $r=1$, then the assertion results from Lemma~\ref{lem:icyc}.
Suppose $r\ge2$.  Then $c_1(\mc L_2)\dotsm c_1(\mc L_r)\mc F\in\I K_1$
by Lemma~\ref{lem:deg}.  Hence
 $$\textstyle\int c_1(\mc L_1)c_1(\mc L_2)\dotsm c_1(\mc L_r)\mc F =
  \int c_1(\mc N_1)c_1(\mc L_2)\dotsm c_1(\mc L_r)\mc F.$$
 Similarly, $c_1(\mc N_1)c_1(\mc L_3)\dotsm c_1(\mc L_r)\mc F\in\I K_1$,
 and so
 $$\textstyle
  \int c_1(\mc N_1)c_1(\mc L_2)c_1(\mc L_3) \dotsm c_1(\mc L_r)\mc F =
  \int c_1(\mc N_1)c_1(\mc N_2)c_1(\mc L_3)\dotsm c_1(\mc L_r)\mc F.$$
 Continuing in this fashion yields the assertion.
\end{proof}

\begin{prp}\label{prp:pbnq}
   Let $g\:X'\to X$ be an $S$-map.  Let $\mc L,\,\mc N\in\Pic(X)$.

\tu{(1)} If $\mc L$ and $\mc N$ are numerically equivalent, then so are
$g^*\mc L$ and $g^*\mc N$.

\tu{(2)} Conversely, when $g$ is proper and surjective, if $g^*\mc L$
and $g^*\mc N$ are numerically equivalent, then so are $\mc L$ and $\mc
N$.
 \end{prp}\begin{proof}
 Let $Y'\subset X'$ be a closed integral curve with $\mc O_{Y'}\in \I
F_1(X'/S)$.  Set $Y:=g(Y')$ and give $Y$ its induced reduced
structure.  Then Proposition~\ref{prp:prjfm} yields
 \begin{align*}
 \textstyle\int c_1(g^*\mc L)[Y']&=\deg(Y'/Y)
 \textstyle\int c_1(\mc L)[Y] \text{ and }\\
 \textstyle\int c_1(g^*\mc N)[Y']&=\deg(Y'/Y)\textstyle\int c_1(\mc N)[Y].
 \end{align*}
 Part (1) follows.

Conversely, suppose $g$ is proper and surjective.  Let $Y\subset X$ be a
closed integral curve with $\mc O_{Y}\in \I F_1(X/S)$.  Then $Y$ is a
complete curve in the fiber $X_s$ over a closed point $s\in S$.  Hence,
since $g$ is proper, there exists a complete curve $Y'$ in $X'_s$
mapping onto $Y$.  Indeed, let $y\in Y$ be the generic point, and $y'\in
g^{-1}Y$ a closed point; let $Y'$ be the closure of $y'$ given $Y'$ its
induced reduced structure.  Plainly $\mc O_{Y'}\in \I F_1(X'/S)$ and
$\deg(Y'/Y)\neq0$.  The two equations displayed above now yield Part (2).
 \end{proof}

\begin{dfn}\label{dfn:deg}
 Assume $S$ is Artin, and $X$ a proper curve.  Let $\mc
L\in\Pic(X)$.  Define its {\it degree\/}  $\deg(\mc L)$ by the formula
	$$\deg(\mc L):=\textstyle\int c_1(\mc L).$$
 Let $D$ be a divisor on $X$.  Define its {\it degree\/} $\deg(D)$ by
$\deg(D):=\deg(\mc O_{\!X}(D))$.
 \end{dfn}

\begin{prp}\label{prp:Riemann}
 Assume $S$ is Artin, and $X$ a proper curve.

\tu{(1)} The map $\deg\:\Pic(X)\to\bb Z$ is a homomorphism.

\tu{(2)} Let $D\subset X$ be an effective divisor. Then
 $$\deg(D) = \dim \tu H^0(\mc O_{\!D}).$$

\tu{(3) (Riemann's Theorem)} Let $\mc L\in\Pic(X)$.  Then
 $$\chi(\mc L) = \deg(\mc L)+\chi(\mc O_{\!X}).$$

\tu{(4)} Suppose $X$ is integral, and let $g\:X'\to X$ be the
normalization map.  Then
	$$\deg(\mc L)=\deg(g^*\mc L).$$
 \end{prp}
\begin{proof}
  Part (1) results from Theorem~\ref{thm:1stprpts}~(2).  And Part
(2) results from Lemma \ref{lem:idiv}.  As to Part (3), note $\deg(\mc
L^{-1})=-\deg(\mc L)$ by Part (1).  And the definitions yield $\deg(\mc
L^{-1})=\chi(\mc O_{\!X})-\chi(\mc L)$.  Thus Part (3) holds.  Finally,
Part (3) results from the definition and the Projection Formula.
 \end{proof}

\begin{dfn}\label{dfn:pa}
 Assume $S$ is Artin, and $X$ a proper surface.  Given a divisor $D$ on
$X$, set
	$$p_a(D):=1-\chi\bigl(c_1(\mc O_{\!X}(D))\,\mc O_{\!X}\bigr).$$
 \end{dfn}

\begin{prp}\label{prp:pa}
 Assume $S$ is Artin, and $X$ a proper surface.  Let $D$ and $E$ be
divisors on $X$.  Then
	$$p_a(D+E)=p_a(D)+p_a(E)+(D\cdot E)-1.$$

Furthermore, if  $D$  is effective, then 
	$$p_a(D)=1-\chi(\mc O_{\!D});$$
 in other words, $p_a(D)$ is equal to the \textup{arithmetic genus} of $D$. 
 \end{prp}\begin{proof}
 The assertions result from Lemmas~\ref{lem:add} and \ref{lem:div}.
 \end{proof}

\begin{prp}[Riemann--Roch for surfaces]\label{prp:RR}
Assume $S$ is the spectrum of a field, and $X$ is a reduced, projective,
equidimensional, Cohen--Macaulay surface. Let $\bigomega$ be a dualizing
sheaf, and set $\mc K:=\bigomega-\mc O_{\!X}$.  Let $D$ be a divisor on
$X$.  Then $\mc K\in \I K_1$; furthermore,
   $$p_a(D) = \frac{(D^2)+(D\cdot\mc K)}2 + 1 \quad \text{and}\quad
 \chi(\mc O_{\!X}(D)) = \frac{(D^2)-(D\cdot\mc K)}2+\chi(\mc
 O_{\!X}).$$
 If $X/S$ is Gorenstein, that is, $\bigomega=\mc O_{\!X}(K)$ for some
``canonical'' divisor $K$, then
 $$p_a(D) = \frac{\bigl(D\cdot(D+K)\bigr)}2 + 1 \quad \text{and}\quad
 \chi(\mc O_{\!X}(D))
  = \frac{\bigl(D\cdot(D-K)\bigr)}{2}+\chi(\mc O_{\!X}).$$
 \end{prp}\begin{proof}
  Since $X$ is reduced, $\bigomega$ is isomorphic to $\mc O_{\!X}$ on a
dense open subset of $X$ by \cite[ (2.8), p.~8]{AK70}.  Hence $\mc K\in
\I K_1$.

 Set $\mc L:=\mc O_{\!X}(D)$.  Then $(D^2):=\int c_1(\mc L)^2=-\int
c_1(\mc L)c_1(\mc L^{-1})$ by Parts (1) and (2) of
Theorem~\ref{thm:1stprpts}.
  Now, the definitions yield
 \begin{align*}
  c_1(\mc L)(-c_1(\mc L^{-1})\mc O_{\!X}+\mc K)
  &=c_1(\mc L)(\mc L-2\mc O_{\!X}+\bigomega)\\
  &=\mc L+\bigomega-3\mc O_{\!X}+2\mc L^{-1}-\mc L^{-1}\ox\bigomega
  \text{\quad and} \\
 c_1(\mc L)(-c_1(\mc L^{-1})\mc O_{\!X}-\mc K)
  &=c_1(\mc L)(\mc L-\bigomega)
  =\mc L-\bigomega-\mc O_{\!X} + \mc L^{-1}\ox\bigomega.
  \end{align*}
 But, $\tu H^i(\mc L)$ is dual to $\tu H^{2-i}(\mc L^{-1}\ox\bigomega)$
by duality theory; see \cite[Cor.~7.7, p.~244]{Ha83}, where $k$ needn't
be taken algebraically closed.  So $\chi(\mc L)=\chi(\mc
L^{-1}\ox\bigomega)$.  Similarly, $\chi(\mc O_{\!X})=\chi(\bigomega)$.
Therefore,
  $$(D^2)+(D\cdot\mc K)=2(\chi(\mc L^{-1})-\chi(\mc O_{\!X}))
 \quad \text{and}\quad
 (D^2)-(D\cdot\mc K)=2(\chi(\mc L) -\chi(\mc O_{\!X})).$$
 Now, $-c_1(\mc O_{\!X}(D))\,\mc O_{\!X}=\mc L^{-1}-\mc O_{\!X}$.  The
first assertion follows.

Suppose $\bigomega$ is invertible.  Then $-c_1(\bigomega^{-1})\mc
O_{\!X}=\mc K$ owing to the definitions.  And $-\int c_1(\mc
L)c_1(\bigomega^{-1})= \int c_1(\mc L)c_1(\bigomega)$ by Part (2) of
Theorem~\ref{thm:1stprpts}.  Therefore, $(D\cdot\mc K)=(D\cdot K)$.
Hence Part (2) of Theorem~\ref{thm:1stprpts} yields the second
assertion.
 \end{proof}

\begin{thm}[Hodge Index]\label{thm:HIT}
 Assume $S$ is the spectrum of a field, and $X$ is a geometrically
irreducible complete surface.  Assume there is an $\mc H\in \Pic(X)$ such
that $\int c_1(\mc H)^2>0$.  Let $\mc L\in\Pic(X)$.  Assume $\int
c_1(\mc L)c_1(\mc H)=0$ and $\int c_1(\mc L)^2\ge0$.  Then $\mc L$ is
numerically equivalent to $\mc O_{\!X}$.
 \end{thm}\begin{proof}
 We may extend the ground field to its algebraic closure owing to
Proposition~\ref{prp:bschg}.  Furthermore, we may replace $X$ by its
reduction; indeed, the hypotheses are preserved due to
Lemma~\ref{lem:icyc}, and the conclusion is preserved due to
Definition~\ref{dfn:neq}.

By Chow's Lemma, there is a surjective map $g\:X'\to X$ where $X'$ is an
integral projective surface.  Furthermore, we may replace $X'$ by its
normalization.  Now, we may replace $X$ by $X'$ and $\mc H$ and $\mc L$
by $g^*\mc H$ and $g^*\mc L$.  Indeed, the hypotheses are preserved
due to the Projection Formula, Proposition~\ref{prp:prjfm}.  And the
conclusion is preserved due to Part (2) of Proposition~\ref{prp:pbnq}.

By way of contradiction, assume that there exists a closed integral
subscheme $Y\subset X$ such that $\int c_1(\mc L)\mc O_{\!Y}\neq0$.  Let
$g\:X'\to X$ be the blowing-up along $Y$, and $E:=g^{-1}Y\subset X$ the
exceptional divisor.  Let $E_1,\dotsc,E_s$ be the irreducible components
of $E$, and give them their induced reduced structure.

Since $X$ is normal, it has only finitely may singular points.  Off
them, $Y$ is a divisor, and $g$ is an isomorphism.  Hence one of the
$E_i$, say $E_1$ maps onto $Y$, and the remaining $E_i$ map onto points.
Therefore, $\int c_1(g^*\mc L)[E_1]=\int c_1(g\mc L)[Y]$ and $\int
c_1(g^*\mc L)[E_i]=0$ for $i\ge2$ by the Projection Formula.  Hence
Lemma~\ref{lem:icyc} yields $\int c_1(g^*\mc L)[E] =\int c_1(\mc L)[Y]$.
The latter is nonzero by the new assumption, and the former is equal to
$\int c_1(g^*\mc L) c_1(\mc O_{\!X'}(E))$ by Lemma~\ref{lem:idiv}.

Set $\mc M:=\mc O_{\!X'}(E)$.  Then $\int c_1(g^*\mc L) c_1(\mc M)\neq0$.
Moreover, by the Projection Formula, $\int c_1(g^*\mc H)>0$ and $\int
c_1(g^*\mc L)c_1(g^*\mc H)=0$ and $\int c_1(g^*\mc L)^2\ge0$.  Let's
prove this situation is absurd.  First, replace $X$ by $X'$ and $\mc H$
and $\mc L$ by $g^*\mc H$ and $g^*\mc L$.

Let $\mc G$ be an ample invertible sheaf on $X$.  Set $\mc H_1:=\mc
G^{\ox m}\ox \mc M$.  Then 
	$$\textstyle\int c_1(\mc L)c_1(\mc H_1)
	=m\int c_1(\mc L)c_1(\mc G)+\int c_1(\mc L)c_1(\mc M)$$
 by additivity (see Part (2) of Theorem~\ref{thm:1stprpts}).  Now,
$\int c_1(\mc L)c_1(\mc M)\neq0$.  Hence there is an $m>0$ so that $\int
c_1(\mc L)c_1(\mc H_1)\neq 0$ and so that $\mc H_1$ is ample.

Set $\mc L_1:=\mc L^{\ox p}\ox \mc H^{\ox q}$.  Since $\int c_1(\mc
L)c_1(\mc H)=0$, additivity yields
 \begin{align*}
 \textstyle\int c_1(\mc L_1)^2 &= \textstyle
	p^2\int c_1(\mc L)^2+q^2\int c_1(\mc H)^2,\text{ and}\\
 \textstyle\int c_1(\mc L_1)c_1(\mc H_1) &= \textstyle
	p\int c_1(\mc L)c_1(\mc H_1)+q\int c_1(\mc H)c_1(\mc H_1).
 \end{align*}
 Since $\int c_1(\mc L)c_1(\mc H_1)\neq 0$, there are $p,\,q$ with
$q\neq0$ so that $\int c_1(\mc L_1)c_1(\mc H_1)=0$.  Then $\int c_1(\mc
L_1)^2 >0$ since $\int c_1(\mc L)^2\ge0$ and $\int c_1(\mc H)^2>0$.
Replace $\mc L$ by $\mc L_1$ and $\mc H$ by $\mc H_1$.  Then $\mc H$ is
ample, $\int c_1(\mc L)c_1(\mc H)=0$ and $\int c_1(\mc L)^2 >0$.

Set $\mc N:=\mc L^{\ox n}\ox \mc H^{-1}$ and $\mc H_1:=\mc L\ox\mc H^a$.
Take $a>0$ so that $\mc H_1$ is ample.  By additivity,
	$$\textstyle\int c_1(\mc N)c_1(\mc H_1)
	= n\int c_1(\mc L)^2-a\int c_1(\mc H)^2.$$
 Take $n>0$  so that  $\int c_1(\mc N)c_1(\mc H_1)>0$.  Then additivity
and Proposition~\ref{prp:amp} yield 
 \begin{align*}
 \textstyle\int c_1(\mc N)c_1(\mc H)&=-\int c_1(\mc H)^2<0,\\
 \int c_1(\mc N)^2&= n^2\int c_1(\mc L)^2+\int c_1(\mc H)^2>0.
 \end{align*}
 But this situation stands in contradiction to the next lemma.
 \end{proof}

\begin{lem}\label{lem:final}
 Assume $S$ is the spectrum of a field, and  $X$ is an integral
surface.   Let $\mc N\in\Pic(X)$, and assume $\int c_1(\mc N)^2>0$.  Then
these conditions are equivalent:\smallskip

\tu{\hphantom{i}(i)} For every ample sheaf $\mc H$, we have $\int c_1(\mc N)
c_1(\mc H)>0$.

\tu{(i$'$)} For some ample sheaf $\mc H$, we have $\int c_1(\mc
N)c_1(\mc H)>0$. 

\tu{(ii)} For some $n>0$, we have $\tu H^0(\mc N^{\ox n})\neq0$.  
 \end{lem}\begin{proof}
 Suppose (ii) holds.  Then there exists an effective divisor $D$ such
that $\mc N^{\ox n}\simeq \mc O_{\!X}(D)$.  And $D\neq0$ since $\int c_1(\mc
N)^2>0$.  Hence (i) results as follows:
	$$\textstyle\int c_1(\mc N) c_1(\mc H)
 =\int c_1(\mc H)c_1(\mc N) =\int c_1(\mc H)[D]>0$$
 by symmetry, by Lemma~\ref{lem:idiv}, and by Proposition~\ref{prp:amp}.

Trivially, (i) implies (i$'$).  Finally, assume (i$'$), and let's prove
(ii).  Let $\bigomega$ be a dualizing sheaf for $X$; then $\bigomega$ is
torsion free of rank 1, and $\tu H^2(\mc L)$ is dual to $\Hom(\mc L,
\bigomega)$ for any coherent sheaf $\mc F$ on $X$; see
\cite[p.~149-17]{FGA}, \cite[(1.3), p.~5, and (2.8), p.~8]{AK70}, and
\cite[Prp,~7.2, p.~241]{Ha83}.  Set $\mc K:=\bigomega-\mc O_{\!X}\in \I K_1$.

Suppose $\mc L$ is invertible and $\tu H^2(\mc L)$ is nonzero.  Then
there is a nonzero map $\mc L\to\bigomega$, and it is injective since
$X$ is integral.  Let $\mc F$ be its cokernel.  Then
     $$\mc K =  \mc F -c_1(\mc L^{-1})\mc O_{\!X} \text{ in } \I K_1.$$
 Hence Proposition~\ref{prp:amp}, symmetry, and additivity yield
    $$\textstyle\int c_1(\mc H)\mc K\ge \int c_1(\mc L)c_1(\mc H).$$
 
Take $\mc L:=\mc N^{\ox n}$.  Then $\int c_1(\mc H)c_1(\mc L)=n\int
c_1(\mc H)c_1(\mc N)$ by additivity.  But $\int c_1(\mc N)c_1(\mc H)>0$
by hypothesis.  Hence $\tu H^2(\mc N^{\ox n})$ vanishes for $n\gg0$.  Now,
	$$\chi(\mc N^{\ox n})=\textstyle\int c_1(\mc
	N)^2\,\binom{n+1}2+a_1n +a_0$$
 for some $a_1,\,a_0$ by Snapper's Theorem, Theorem~\ref{thm:Snapper}.
But  $\int c_1(\mc N)^2>0$ by hypothesis.  Therefore, (ii) holds.
 \end{proof}

\begin{cor}\label{cor:conn} Assume $S$ is the spectrum of a field, and
$X$ a geometrically irreducible projective $r$-fold with $r\ge2$.  Let
$D,\,E$ be effective divisors, with $E$ possibly trivial.  Assume $D$ is
ample.  Then $D+E$ is connected.
 \end{cor}
\begin{proof}
 Plainly we may assume the ground field is algebraically closed and $X$
is reduced.  Fix $n>0$ so that $nD+E$ is ample; plainly we may replace
$D$ and $E$ by $nD+E$ and 0.
Proceeding by way of contradiction, assume $D$ is the
disjoint union of two closed subschemes $D_1$ and $D_2$.  Plainly $D_1$
and $D_2$ are divisors; so $D=D_1+D_2$.

Proceed by induction on $r$.  Suppose $r=2$.  Then, since $D_1$ and $D_2$
are disjoint,   $(D_1\cdot D_2)=0$ by Lemma~\ref{lem:idiv}. 
 Now, $D$ is ample.  Therefore, Proposition
 \ref{prp:amp} yields
	$$(D_1^2)=(D.D_1)>0\text{ and }(D_2^2)=(D.D_2)>0.$$
 These conclusions contradict  Theorem~\ref{thm:HIT}
with $\mc H:=\mc O_X(D_1)$ and $\mc L:=\mc O_X(D_2)$.

Finally, suppose $r\ge3$.  Let $H$ be a general hyperplane section of
$X$.  Then $H$ is integral by Bertini's Theorem \cite[Thm.~12,
p.~374]{Se50}.  And $H$ is not a component of $D$.  Set $D':=D\cap H$
and $D_i':=D_i\cap H$.  Plainly $D_1'$ and $D_2'$ are disjoint, and
$D'=D_1'+D_2'$; also, $D'$ is ample.  So induction yields the desired
contradiction.
 \end{proof}

\end{document}